\documentclass{amsart}
\ifdefined\pdfminorversion\pdfminorversion=7\relax\fi
\usepackage[utf8]{inputenc}
\usepackage{amsfonts, amssymb}
\usepackage{amssymb,amsmath,amscd,euscript,verbatim,array}
\usepackage{geometry}
\usepackage{anysize}
\usepackage{fancyhdr}
\usepackage{indentfirst}
\usepackage{graphicx}
\usepackage{color}
\usepackage{ifpdf}
\usepackage{marginnote}
\usepackage{enumitem}
\usepackage{mathrsfs}
\usepackage[hidelinks]{hyperref}
\relax {}
\usepackage{bm}
\setlist[itemize]{leftmargin=2em}

\newtheorem{Theorem}{Theorem}
\newtheorem{MainTheorem}{Theorem}

\newtheorem{Definition}{Definition}[section]
\newtheorem{Proposition}[Definition]{Proposition}
\newtheorem{Corollary}{Corollary}
\newtheorem{Lemma}[Definition]{Lemma}

\theoremstyle{definition}
\newtheorem{Remark}[Definition]{Remark}

\newcommand{\dist}{\operatorname{dist}}
\newcommand{\bthm}{\begin{Theorem}}
	\newcommand{\ethm}{\end{Theorem}}
\newcommand{\bpr}{\begin{Proposition}}
	\newcommand{\epr}{\end{Proposition}}
\newcommand{\blm}{\begin{Lemma}}
	\newcommand{\elm}{\end{Lemma}}
\newcommand{\bex}{\begin{Exercise}}
	\newcommand{\eex}{\end{Exercise}}
\newcommand{\be}{\begin{equation}}
	\newcommand{\ee}{\end{equation}}
\newcommand{\beal}{\begin{aligned}}
	\newcommand{\enal}{\end{aligned}}
\newcommand{\brm}{\begin{Remark}}
	\newcommand{\erm}{\end{Remark}}

\begin{document}
\raggedbottom

	\title[Quantitative Accumulation of Invariant Circles Near a Non-splitting Separatrix]{Quantitative Accumulation of Invariant Circles Near a Non-splitting Separatrix}
	
	\author{Deyu Kong}
	\address{School of Mathematics and Statistics, Beijing Institute of Technology, Beijing 100081, China}
	\email{dykong@bit.edu.cn}

	\author{Lin Wang}
	\address{School of Mathematics and Statistics, Beijing Institute of Technology, Beijing 100081, China}
	\email{lwang@bit.edu.cn}

	

		\begin{abstract}
For a smooth area-preserving cylinder map with a non-splitting separatrix and an exact symplectic linearization at its saddle, we prove accumulation by smooth essential invariant circles when the boundary transition is sufficiently $C^4$-close to a translation. As an application, we give an answer to Herman's open problem in 1983.
    \end{abstract}
	\maketitle
	\vspace*{-1pt}
	{\fontsize{9}{9.5}\selectfont\tableofcontents}
	\clearpage
	
	\section{Introduction}\label{sec:intro}

Invariant circles are central to the dynamics of area-preserving cylinder maps. Classical KAM theory and its finite-regularity and geometric refinements give persistence under suitable arithmetic and smallness assumptions \cite{Kol54,Arn63,Mos62,Rus83,Her86,Pos82,Sal04,TL25,dLGJV05,GEdL08}. Aubry--Mather and converse KAM theory describe ordered invariant sets and mechanisms that destroy invariant circles \cite{Mat82,Mat86,Mat87,Mat88,For94,Wan12,Wan22,Wan27,Bes00,CW13}. A non-splitting separatrix lies at the intersection of these questions: its hyperbolic dynamics obstruct a direct perturbative argument, but its invariant geometry may still support accumulating quasiperiodic circles.

Let $f:\mathbb A\to\mathbb A$ be a smooth area-preserving diffeomorphism isotopic to the identity. A \emph{non-splitting separatrix} $\Sigma$ is an essential closed curve formed by coinciding stable and unstable branches of a hyperbolic saddle. Two oppositely oriented connections between a lifted saddle $o_0$ and its horizontal translate $o_1=\tau(o_0)$ may form a heteroclinic cycle in the universal cover $\mathbb R^2$. The bounded invariant Jordan disk $D$ enclosed by this cycle is called an \emph{island}; this terminology does not presume an elliptic fixed point. Following \cite{KK21}, invariant circles \emph{accumulate} on $\Sigma$ or $\partial D$ if their Hausdorff distance $d_H$ from that boundary tends to zero. We seek quantitative conditions for essential circles to accumulate from the exterior and for contractible circles to accumulate on $\partial D$ from within.

Herman's hyperbolic perturbation provides a concrete instance of this problem. In \cite[Chapter~II, \S4]{Her83}, he constructed a smooth family $f_\lambda$ with two non-splitting separatrices and destroyed small-rotation invariant circles by adding a localized perturbation \cite[Chapter~II, Theorem~4.6]{Her83}. He then asked whether the unmodified family already has this destruction property \cite[Chapter~II, Remark~4.7(3)]{Her83}. We give a negative answer: for each sufficiently small fixed $\lambda$, the original map admits smooth invariant graphs with arbitrarily small nonzero rotation numbers. We also obtain contractible invariant circles inside its island, even though the original cutoff gives an unstable segment of fixed points rather than an elliptic center.

Two earlier results explain both the possibility of accumulation and the need for additional hypotheses. Fayad--Krikorian \cite{FK09} proved that a smooth Diophantine invariant circle is accumulated by a positive-area union of invariant circles, without a twist assumption. Katok--Krikorian \cite{KK21} developed a logarithmic return-map construction near a non-splitting separatrix and proved accumulation for sufficiently small perturbations of a fixed integrable model. Their counterexamples \cite[Theorem~B]{KK21} show that the existence of a non-splitting separatrix alone is insufficient. The present criterion controls the boundary transition along the separatrix directly, without requiring closeness of the full map to a reference integrable model.

\subsection{Difficulties}\label{subsec:difficulties}
Three issues must be addressed.
\begin{enumerate}[label=\textup{(\arabic*)},leftmargin=2em,topsep=4pt,itemsep=5pt,parsep=0pt]
\item \emph{Degenerating rotation numbers near the saddle.}
Passage times through a saddle neighborhood diverge, while nearby rotation numbers tend to zero. For the twist map
\[
f_\varepsilon(x,y)=(x+y,\,y+\varepsilon F(x+y)),
\]
an invariant parameterization $x(t)=t+u(t)$ of rotation number $\omega$ satisfies
\[
u(t+\omega)-2u(t)+u(t-\omega)=\varepsilon F(t+u(t)).
\]
Linearization at $u=0$ gives, for $m\ne0$,
\[
\widehat u_m=-\frac{\varepsilon\widehat F_m}{4\sin^2(\pi m\omega)}.
\]
If $\widehat F_1\ne0$, the first correction has size comparable to $\varepsilon\omega^{-2}|\widehat F_1|$ as $\omega\to0$. Thus a perturbative expansion about horizontal circles does not remain uniformly small in the original coordinates \cite{Mos62,Her83,Sal04}.

\item \emph{Global passage and weak hyperbolicity.}
Local hyperbolicity does not control the transition along the entire separatrix, as the examples of \cite{KK21} demonstrate. In Herman's family, the one-step hyperbolic exponent is $\log(1+\lambda)\sim\lambda$, and passage along a compact portion of the separatrix takes $\mathcal O(\lambda^{-1})$ iterates. The transverse derivatives of the resulting two-dimensional transition need not be uniformly bounded as $\lambda\to0$. A useful criterion must separate the boundary distortion from these map-dependent derivatives and control the accumulated boundary error.

\item \emph{Reconstruction and the topology of the island.}
A KAM graph for a normalized return map is not yet an embedded circle invariant under the original map. Reconstruction must retain the integer part of the return frequency, close the curve, and recover invariance under $f_0$ from invariance under an iterate $f_0^M$. Inside an island, the orbit passes two saddle neighborhoods, and its rotation about an interior fixed point differs from its horizontal rotation, which is zero. Moreover, neither a separatrix containing a saddle nor a degenerate fixed segment supplies the smooth Diophantine circle required to start the theorem of \cite{FK09}.
\end{enumerate}

\subsection{Main results}\label{subsec:general-overview}
The hypotheses are \textbf{(G1)}, a non-splitting separatrix, and \textbf{(G2)}, an exact symplectic linearization at the saddle; see Section~\ref{subsec:geometric-assumptions}. The finite-step passage along a connection defines a boundary transition map in prescribed logarithmic coordinates. Its $C^4$ distance from translations is specified in~\eqref{eq:transition-distance-translations}. Throughout, smooth means $C^\infty$, and the quasiperiodic circles below carry dynamics smoothly conjugate to irrational rotations.

\begin{MainTheorem}[A single separatrix]\label{thm:intro-criterion}
Let $f_0$ be a smooth area-preserving cylinder diffeomorphism isotopic to the identity and fixing a hyperbolic saddle. Set $G=f_0^M$, $M\in\mathbb N$, with the lift $\widetilde G=\widetilde f_0^M$ fixing that saddle. Suppose $G$ satisfies \textbf{(G1)--(G2)} on a chosen side of $\Sigma$ and its boundary transition is sufficiently $C^4$-close to a translation. Then:
\begin{itemize}[topsep=3pt,itemsep=3pt,parsep=0pt]
\item For some constant-type $\alpha$ and all sufficiently large $k$, there is a smooth quasiperiodic essential $f_0$-invariant circle $\Gamma_k$ on that side, with
\[
\omega_k=\frac{1}{M(k-\alpha)},\qquad
d_H(\Gamma_k,\Sigma)\le C_Ge^{-k\nu/2},
\]
where $\nu$ is the expanding hyperbolic exponent of $G$, with rotation oriented as in \textbf{(G1)}. Every collar on this side contains a positive-area union of smooth quasiperiodic essential invariant circles.
\item Under positive twist bounded away from zero, $\Gamma_k$ are smooth graphs. If $\Sigma=\operatorname{graph}\psi$ is Lipschitz, they converge uniformly to $\psi$ from the chosen side at an exponential rate.
\end{itemize}
\end{MainTheorem}

Theorem~\ref{thm:intro-criterion} treats one side of a single connection. When two opposite connections enclose an island, each boundary transition must be controlled. An exterior orbit follows one connection; an interior orbit passes both saddles before completing a turn. This gives the following distinction between essential and contractible invariant circles.

\begin{MainTheorem}[An island bounded by two connections]\label{thm:intro-island}
Let $f_0$, $\widetilde f_0$, and $G$ be as in Theorem~\ref{thm:intro-criterion}. Suppose the connections $S_\pm$ bound an island $D$ as in Definition~\ref{def:island-geometry}, with compatible \textbf{(G2)} charts in all four sectors. Assume both boundary transitions $b_\pm$ are sufficiently $C^4$-close to translations. Then:
\begin{itemize}[topsep=3pt,itemsep=3pt,parsep=0pt]
\item \emph{Exterior circles.} Theorem~\ref{thm:intro-criterion} applies above $\Sigma_+=\pi(S_+)$ and below $\Sigma_-=\pi(S_-)$, with positive and negative horizontal rotation numbers, respectively.
\item \emph{Interior circles.} The lift $\widetilde f_0$ has a fixed point in $D$. For any such point $p$, there are smooth quasiperiodic invariant Jordan curves $C_k\subset D$ enclosing $p$, contractible on the cylinder, such that
\[
\rho_{\mathrm{in}}(C_k)=\frac{1}{M(2k-\alpha_{\mathrm{in}})},\qquad
d_H(C_k,\partial D)\le C_Ge^{-k\nu/2}
\]
for all sufficiently large $k$. Here $\alpha_{\mathrm{in}}$ and $\rho_{\mathrm{in}}(C_k)$ are of constant type, and rotation about $p$ is oriented by $S_+$ followed by $S_-$. Every interior collar of $\partial D$ contains a positive-area union of smooth quasiperiodic invariant circles enclosing $p$.
\item \emph{Horizontal rotation.} Every orbit in $\overline D$ has horizontal rotation number zero. Under positive twist bounded away from zero, $\rho_x^->0$ throughout the upper exterior and $\rho_x^+<0$ throughout the lower exterior; see~\eqref{eq:lower-upper-rotation}.
\end{itemize}
\end{MainTheorem}

No ellipticity, stability, isolation, or uniqueness is required of the fixed point $p$ in Theorem~\ref{thm:intro-island}; the starting index of $C_k$ may depend on $p$. The positive-area assertions concern families in each collar, not the displayed countable sequences alone. The smallness thresholds in Theorems~\ref{thm:intro-criterion} and~\ref{thm:intro-island} are uniform for fixed arithmetic data and $\nu$ in a compact subset of $(0,\infty)$; starting indices and distance constants may depend on the map. Theorem~\ref{thm:boundary-transition-criterion}, Corollaries~\ref{cor:f0-invariant} and~\ref{cor:graph-corollary}, and Proposition~\ref{prop:island} give the precise statements. Lemma~\ref{lem:one-step-transition} supplies a sufficient condition in terms of one step of the original map.

\paragraph{Application to Herman's question.}\label{subsec:herman-overview}
For Herman's family, the forcing and the two separatrix functions are
\[
\varphi_{\lambda^2}=g_\lambda+g_\lambda^{-1}-2\operatorname{id},\qquad
\psi_+=g_\lambda-\operatorname{id},\qquad
\psi_-=g_\lambda^{-1}-\operatorname{id};
\]
see Section~\ref{subsec:auto-1-2}. Let $D_\lambda$ be the island between their graphs over $(0,1)$. Theorems~\ref{thm:intro-criterion} and~\ref{thm:intro-island} apply for every sufficiently small parameter, rather than only along a parameter subsequence.

\begin{MainTheorem}[Herman's hyperbolic family]\label{thm:intro-herman}
There exist $\lambda_*>0$ and $C>0$ such that, for every $0<\lambda<\lambda_*$, the Herman map $f_\lambda$ has the following properties.
\begin{itemize}[topsep=3pt,itemsep=3pt,parsep=0pt]
\item \emph{Exterior graphs at arbitrarily small rotation.} For every sufficiently large $k$, there are smooth quasiperiodic invariant graphs $u_{\lambda,k}^+>\psi_+$ and $u_{\lambda,k}^-<\psi_-$ with constant-type rotation numbers $\pm\omega_{\lambda,k}$, where $\omega_{\lambda,k}\to0^+$, such that
\[
\|u_{\lambda,k}^\pm-\psi_\pm\|_\infty\le C_\lambda e^{-k\nu_\lambda/2},
\qquad \|(u_{\lambda,k}^\pm)'\|_\infty\le C\lambda.
\]
Here $M_\lambda=\lfloor\lambda^{-1}\rfloor$ and $\nu_\lambda=M_\lambda\log(1+\lambda)$.
\item \emph{Island and exterior dynamics.} All conclusions of Theorem~\ref{thm:intro-island} hold, including positive-area families in each exterior and interior collar, zero horizontal rotation in $\overline D_\lambda$, and the signs of exterior rotation. For any fixed point $p\in D_\lambda$ of the lift fixing the saddle, the interior curves have clockwise rotation numbers
\[
\frac{1}{M_\lambda(2k-\alpha_{\lambda,\mathrm{in}})},
\]
where $\alpha_{\lambda,\mathrm{in}}$ is of constant type and independent of $k$.
\item \emph{The original cutoff.} The point $(1/2,0)$ lies in the interior of a horizontal segment of fixed points. The map is a shear near $(1/2,0)$, which is Lyapunov unstable.
\end{itemize}
The starting indices and $C_\lambda$ may depend on $\lambda$. In particular, for each such fixed $\lambda$ and every $r>0$, an exterior invariant graph has rotation number in $(0,r)$.
\end{MainTheorem}

\begin{figure}[htbp]
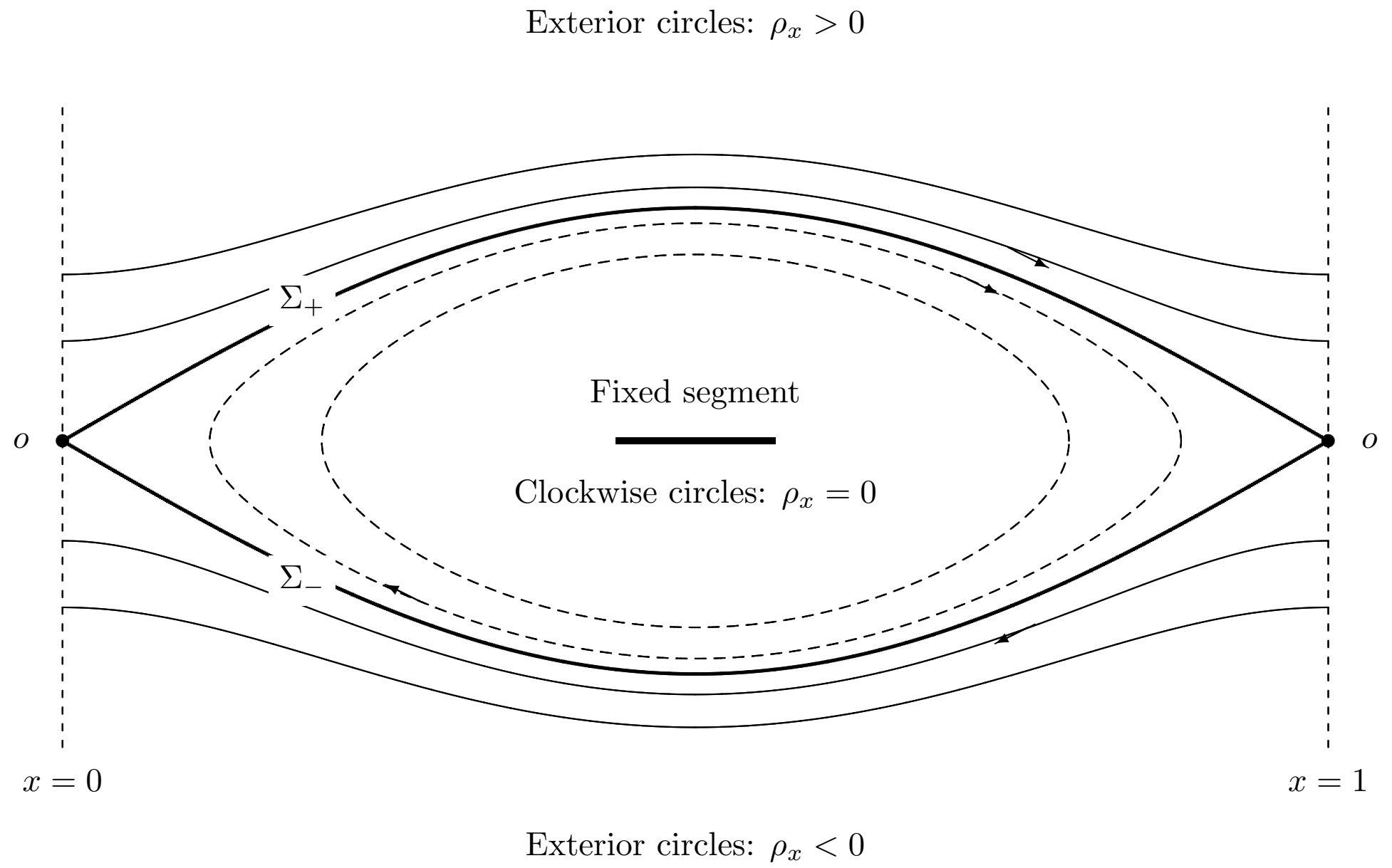

\centering
\begingroup
\setlength{\unitlength}{1pt}
\normalsize\color{black}
\scalebox{0.88}{%
}
\end{picture}%
}
\endgroup

\caption{Theorem~\ref{thm:intro-herman} on a cylinder cut along $x=0$; the vertical edges are identified. Thick solid curves are the separatrices $\Sigma_\pm$, thin solid curves are exterior invariant circles, and dashed curves are contractible interior invariant circles. Arrows indicate the direction of motion. The interior circles surround the fixed segment and have horizontal rotation number $\rho_x=0$. }
\label{fig:herman-island-portrait}
\end{figure}

For Herman's  parameters, this last assertion has the following quantifiers. Given $a>0$, there is $N_*(a)\in\mathbb N$, independent of the radii $r_N$, such that
\[
N\ge N_*(a),\quad \lambda_N=N^{-a/2},\quad r_N>0
\quad\Longrightarrow\quad
\exists\,\omega_N\in(0,r_N)
\]
with a smooth $f_{\lambda_N}$-invariant graph of rotation number $\omega_N$. Taking $r_N=N^{-a/2-\varepsilon}$ gives the negative answer to Herman's question; Theorem~\ref{cor:herman-window} proves this statement by first fixing $N$ and then choosing the accumulation index. Theorem~\ref{thm:herman-hyperbolic-family} gives the exterior estimates, and Section~\ref{subsec:herman-island-application} proves the island assertions and identifies the fixed segment.

\subsection{Method and comparison with earlier work}\label{subsec:comparison}\label{subsec:auto-1-3}
The argument addresses the three difficulties above in the same order.
\begin{enumerate}[label=\textup{(\arabic*)},leftmargin=2em,topsep=4pt,itemsep=5pt,parsep=0pt]
\item \emph{A return map with nondegenerate logarithmic twist.}
We use the logarithmic renormalization of Katok--Krikorian \cite[Sections~4--8]{KK21}. Near the saddle, logarithmic coordinates turn hyperbolic scaling into phase translation. On the rescaled annuli
\[
v=e^{-k\nu}s,\qquad s\in[1/2,1],
\]
the return map has a nondegenerate logarithmic twist, although the rotation number of the reconstructed circle tends to zero. Our exact linearization hypothesis \textbf{(G2)} is stronger than the saddle normal form used in \cite{KK21}, where nonlinear integrable terms remain. It permits extension of the transition map across open overlaps on the universal cover while retaining the exact symplectic form and the integer part of the return frequency. The resulting return lift is
\[
\mathcal R=\mathcal P\circ\mathcal C_\nu.
\]
This construction separates the vanishing rotation number in the original cylinder from the fixed twist and arithmetic bounds on the return annulus.

\item \emph{Boundary control separated from transverse derivatives.}
Write the boundary map as $b(\theta)=\theta+\sigma+\zeta(\theta)$, with $\zeta$ of zero mean. Expanding $\mathcal P$ at $v=0$ and using area preservation gives
\[
\partial_v\mathcal P_2(\theta,0)=\frac1{b'(\theta)},\qquad
\|\mathcal F_k-\mathcal T_{\Lambda_\nu}\|_{C^3}
\le C_\nu\|\zeta\|_{C^4(\mathbb T)}+C_{f,\nu}e^{-k\nu}.
\]
The persistent term depends only on the one-dimensional boundary displacement. Transverse Taylor remainders are multiplied by $e^{-k\nu}$ and vanish once the map is fixed. Unlike the perturbative hypothesis of \cite[Theorem~A and Remark~1.2]{KK21}, this criterion does not require uniform smallness of the entire two-dimensional transition relative to a fixed integrable model.

For Herman's family, vertical rescaling and $M_\lambda=\lfloor\lambda^{-1}\rfloor$ iterates keep $\nu_\lambda$ in $[1/2,1]$. Since $g_\lambda$ is affine near the endpoints, a linear symplectic change of coordinates gives \textbf{(G2)}. Along the separatrix, the integral coordinate
\[
t=h_\lambda(x)=\int_{x_*}^x\frac{\mathrm du}{W_\lambda(u)}
\]
reduces one step to $t\mapsto t+\lambda+\lambda^2E_\lambda(t)$ with uniform derivative bounds. Composition through $n=\mathcal O(\lambda^{-1})$ steps then gives
\[
\|\zeta_\lambda\|_{C^4(\mathbb T)}\le Cn\lambda^2\le C\lambda.
\]
The reversor gives the corresponding bound on the second connection. No uniform estimate on the transverse derivatives is needed to choose the initial accumulation index for a fixed $\lambda$.

\item \emph{Reconstruction for the original map and for the island.}
Herman's invariant curve theorem \cite[Chapter~V, \S7.10]{Her86} gives a curve for the normalized return map. An orbit extension reconstructs a smooth embedded circle in the original cylinder. Minimality and the order of essential circles recover $f_0$-invariance from $f_0^M$-invariance and give the rotation number $1/[M(k-\alpha)]$. Under positive twist, Birkhoff's theorem first gives a Lipschitz graph; its smooth regular parameterization then gives a smooth graph \cite{Arn09}.

For the island, the two passages compose:
\[
\mathcal R_{\mathrm{in}}=\mathcal R_-\circ\mathcal R_+,
\qquad
\Lambda_{\mathrm{in}}(s)=\sigma_++\sigma_-+\frac2\nu\log s.
\]
Passing to the universal cover of $\mathbb R^2\setminus\{p\}$ records a full turn about $p$. Reconstruction yields contractible circles with rotation $1/[M(2k-\alpha_{\mathrm{in}})]$ about $p$ and zero horizontal rotation. This argument requires no elliptic center. The constructed smooth Diophantine circles then provide the starting points for Fayad--Krikorian \cite[Theorem~1]{FK09}, which yields positive-area families in every exterior and interior collar; see Remark~\ref{rem:scope}.
\end{enumerate}

The distinction from the sine perturbation studied by Liu--Wang \cite{LW25} is also relevant. Their map has separatrix splitting and transverse homoclinic intersections, whereas the stable and unstable branches considered here coincide. Their persistence theorem selects frequencies along a parameter subsequence \cite[Theorem~2]{LW25}. Here each sufficiently small fixed parameter in Herman's original family admits infinitely many smooth invariant graphs with rotation numbers tending to zero.

\subsection{Organization of the paper}\label{subsec:organization}
Section~\ref{sec:criterion} gives the notation, hypotheses, and precise statements. Sections~\ref{sec:transition} and~\ref{sec:kam} construct and normalize the return map. Section~\ref{sec:reconstruction} reconstructs the invariant circles, proves Theorem~\ref{thm:intro-criterion}, and treats the return around the heteroclinic cycle to prove Theorem~\ref{thm:intro-island}. Section~\ref{sec:family} establishes the uniform boundary estimates for Herman's family and proves Theorem~\ref{thm:intro-herman}. Appendix~\ref{sec:appendix} collects the chain-rule estimates and smooth cutoff constructions.
\vspace{2em}

\noindent\textbf{Acknowledgement.}
This work was partially supported by the National Natural Science Foundation of China (Grant No. 12671217 and 12122109).

\section{Notations and precise statements}\label{sec:criterion}
This section gives the conventions and hypotheses used throughout the paper. Theorem~\ref{thm:boundary-transition-criterion} and Corollaries~\ref{cor:f0-invariant} and~\ref{cor:graph-corollary} give the quantitative form of Theorem~\ref{thm:intro-criterion}. Definition~\ref{def:island-geometry} and Proposition~\ref{prop:island} specify the island geometry and the interior conclusion of Theorem~\ref{thm:intro-island}. For Theorem~\ref{thm:intro-herman}, the estimates on the upper branch are stated in Theorems~\ref{thm:herman-hyperbolic-family} and~\ref{cor:herman-window}; Section~\ref{subsec:herman-island-application} treats the lower branch and the island.

\subsection{Notation and conventions}\label{subsec:auto-1-1}
	Write $\mathbb N=\{1,2,\ldots\}$, $\mathbb T=\mathbb R/\mathbb Z$, and $\mathbb A=\mathbb T\times\mathbb R$, with the standard area form $\mathrm dx\wedge\mathrm dy$. Throughout, \emph{smooth} means $C^\infty$, and the same notation is used for a function on the circle and its $1$-periodic lift. The unit horizontal translation on the universal cover is denoted by $\tau(x,y)=(x+1,y)$. When a map has a hyperbolic saddle fixed point, we fix a lift of this saddle to the universal cover $\mathbb R^2$ and choose the unique lift of the map that fixes it. Horizontal rotation numbers are computed with respect to this lift, with their integer parts retained. For a contractible invariant circle enclosing an interior fixed point $p$, its \emph{rotation number about $p$} is the rotation number in the punctured plane, with the orientation and lift specified in Proposition~\ref{prop:island}. We write $f$ for a map on the cylinder and $\widetilde f$ for its chosen lift, and use the same convention for $f_\lambda$, $f_0$, $f_\varphi$, $F_\lambda$, and $G_\lambda$. These lifts satisfy $\widetilde f\circ\tau=\tau\circ\widetilde f$.

    For a lifted orbit $(x_n,y_n)=\widetilde f_0^{\,n}(z)$, we use the lower and upper horizontal rotation numbers
    \begin{equation}\label{eq:lower-upper-rotation}
        \rho_x^-(z):=\liminf_{n\to\infty}\frac{x_n-x_0}{n},
        \qquad \rho_x^+(z):=\limsup_{n\to\infty}\frac{x_n-x_0}{n}.
    \end{equation}
    When they coincide and are finite, their common value is the horizontal rotation number $\rho_x(z)$.
    \par

    The rigid rotation of the circle is denoted by $R_\omega(t)=t+\omega\pmod{1}$. Horizontal phase translations in logarithmic coordinates are denoted by $T_c(\theta,v)=(\theta+c,v)$. Here $T_1$ corresponds to one local iterate, whereas $\tau(x,y)=(x+1,y)$ translates the original spatial coordinate. Their geometric and dynamical roles are different. An irrational number $\omega$ is of \emph{constant type} if there exists $\gamma>0$ such that
	\begin{equation}\label{eq:constant-type}
		\left|\omega-\frac pq\right|\ge\frac\gamma{q^2}
		\qquad(p\in\mathbb Z,\ q\in\mathbb N).
	\end{equation}
	Throughout the paper, we fix a constant-type pair $\xi\in(0,1)$ and $\gamma_\xi>0$ satisfying
	\begin{equation}\label{eq:fix-constant-type}
		\left|\xi-\frac pq\right|\ge\frac{\gamma_\xi}{q^2}
		\qquad(p\in\mathbb Z,\ q\in\mathbb N).
	\end{equation}
	All subsequent smallness thresholds depend on the arithmetic data through $\gamma_\xi$. An \emph{essential circle} in the cylinder is an embedded closed curve whose horizontal projection has absolute degree $1$. For a cylinder map $f$, write its lift as $\widetilde f=(X,Y)$. We call $f$ \emph{exact symplectic} if the $1$-form $Y\,\mathrm dX-y\,\mathrm dx$ has a single-valued primitive on $\mathbb A$. The positive twist condition is $\partial_yX>0$.
	
	We use the standard fact that an area-preserving cylinder diffeomorphism isotopic to the identity is exact symplectic if it has at least one invariant essential circle; equivalently, its flux vanishes. See Stark \cite{Sta88}.
	
	We use the standard metric on $\mathbb T$ and the product metric on $\mathbb A$:
	\[
	d_{\mathbb T}(x,x')=\min_{m\in\mathbb Z}|x-x'-m|,\qquad
	d_{\mathbb A}((x,y),(x',y'))
	=\sqrt{d_{\mathbb T}(x,x')^2+|y-y'|^2},
	\]
	where arbitrary real representatives may be used on the right-hand side defining $d_{\mathbb T}$. For nonempty compact sets, the Hausdorff distance is denoted by
	\[
	d_H(A,B)=\max\Bigl\{\sup_{a\in A}\dist(a,B),\;\sup_{b\in B}\dist(b,A)\Bigr\}.
	\]
	
	Let $\pi_x:\mathbb R^2\to\mathbb R$ be the projection onto the first coordinate; the corresponding projection in logarithmic phase coordinates is denoted by $\pi_1$. If an oriented lifted curve $\gamma:\mathbb R\to\mathbb R^2$ satisfies $\widetilde f\bigl(\gamma(t)\bigr)=\gamma(t+\omega)$ and $\gamma(t+1)=\tau\bigl(\gamma(t)\bigr)$, write $\pi_x\bigl(\gamma(t)\bigr)=t+u(t)$, where $u$ is $1$-periodic. Its rotation number is uniquely determined by the limit
	\[
	\frac{\pi_x\bigl(\widetilde f^{n}(\gamma(t))\bigr)-\pi_x\bigl(\gamma(t)\bigr)}{n}
	=\omega+\frac{u(t+n\omega)-u(t)}n\longrightarrow\omega \qquad (n\to\infty)
	.
	\]

\subsection{Geometric assumptions and boundary transition}\label{subsec:geometric-assumptions}
Let $f : \mathbb{A} \to \mathbb{A}$ be a smooth area-preserving diffeomorphism isotopic to the identity. Denote the universal cover of the cylinder by $\widetilde{\mathbb{A}} = \mathbb{R}^2$, its natural projection by $\pi : \mathbb{R}^2 \to \mathbb{A}$, and the unit horizontal translation by $\tau(x,y) = (x+1, y)$. Assume the following:
	\begin{enumerate}
		\item [\textbf{(G1)}] There is a hyperbolic fixed point $o \in \mathbb{A}$. Fix a lift of $o$ to $\mathbb{R}^2$, still denoted by $o$, and choose the unique lift $\widetilde{f}$ satisfying $\widetilde{f}(o) = o$. There is an embedded arc $\widetilde{\Sigma} \subset \mathbb{R}^2$ joining $o$ to $\tau(o)$ whose interior $\widetilde{\Sigma}\setminus\{o,\tau(o)\}$ is both a chosen unstable branch of $o$ and a chosen stable branch of $\tau(o)$, oriented from $o$ toward $\tau(o)$ under forward iteration. The integer translates $\{\tau^n(\widetilde{\Sigma})\}_{n \in \mathbb{Z}}$ meet only at adjacent endpoints. Their projection $\Sigma := \pi(\widetilde{\Sigma})$ is thus a simple essential closed curve in $\mathbb{A}$, smooth away from the saddle and satisfying $f(\Sigma) = \Sigma$. We call $\Sigma$ a \emph{non-splitting separatrix} of $f$, or simply a separatrix.
		
			\item [\textbf{(G2)}]
		There are an open neighborhood $U(o)\subset\mathbb R^2$ of $o$, an open neighborhood $V\subseteq\mathbb R^2_{(Q,P)}$ of the origin, and a smooth symplectic diffeomorphism $\Phi_o:U(o)\to V$ satisfying $\Phi_o(o)=(0,0)$ and $\Phi_o^*(\mathrm dQ\wedge\mathrm dP)=\mathrm dx\wedge\mathrm dy$. Set $U(\tau^k(o)):=\tau^k(U(o))$ and $\Phi_{\tau^k(o)}:=\Phi_o\circ\tau^{-k}$. In these coordinates:
		\begin{enumerate}[label=(\roman*), font=\normalfont]
			\item Whenever $z=\Phi_o^{-1}(Q,P)\in U(o)$ and $\widetilde f(z)\in U(o)$, the map is exactly linearized as
			\begin{equation}\label{eq:local-linearization}
				(\Phi_o\circ\widetilde f\circ\Phi_o^{-1})(Q,P)
				=(e^\nu Q,\,e^{-\nu}P),\qquad \nu>0.
			\end{equation}
			\item The local outgoing branch $\widetilde\Sigma\cap U(o)$ and incoming branch $\tau^{-1}(\widetilde\Sigma)\cap U(o)$ correspond to $V\cap\{P=0,\,Q\ge0\}$ and $V\cap\{Q=0,\,P\ge0\}$, respectively, and the chosen side of the cylinder corresponds to $V\cap\{Q>0,\,P>0\}$.
		\end{enumerate}
	
	\end{enumerate}
	By continuity and $\widetilde f(o)=o$, we fix $Q_{\mathrm{box}},P_{\mathrm{box}}>0$ sufficiently small so that the closed box $\mathcal B:=[0,Q_{\mathrm{box}}]\times[0,P_{\mathrm{box}}]\Subset V$ and its preimages
	\begin{equation}\label{eq:2.2}
	\begin{gathered}
		U_{\mathrm{box}}(o):=\Phi_o^{-1}(\mathcal B)\Subset U(o),\\
		U_{\mathrm{box}}(\tau^k(o)):=\Phi_{\tau^k(o)}^{-1}(\mathcal B)
		=\tau^k\bigl(U_{\mathrm{box}}(o)\bigr)\Subset U(\tau^k(o))
	\end{gathered}
	\end{equation}
	satisfy $\widetilde f\bigl(U_{\mathrm{box}}(o)\bigr)\Subset U(o)$. Thus~\eqref{eq:local-linearization} applies to every $z\in U_{\mathrm{box}}(o)$, and we use it only for orbit segments remaining in $U_{\mathrm{box}}(o)$. Figure~\ref{fig:local-symplectic-chart} shows the local symplectic chart and its action on the two separatrix branches.

	\begin{figure}[tbp]
\centering
\begingroup
\setlength{\unitlength}{1pt}
\normalsize\color{black}
\scalebox{0.88}{%
\begin{picture}(450,247)
\linethickness{0.45pt}
\put(0,37){\begin{picture}(450,0)
\put(95,195){\makebox(0,0)[c]{(a) Universal cover}}
\put(350,195){\makebox(0,0)[c]{(b) Local coordinates}}
\qbezier[25](15,100)(15,165)(90,165)
\qbezier[25](90,165)(170,165)(170,100)
\qbezier[25](170,100)(170,28)(90,28)
\qbezier[25](90,28)(15,28)(15,100)
\qbezier[333](85,48)(115,69)(145,90)
\qbezier[336](145,90)(115,111.5)(85,133)
\qbezier[336](85,133)(55,111.5)(25,90)
\qbezier[333](25,90)(55,69)(85,48)
\linethickness{0.9pt}
\qbezier[406](12,99.1)(48.5,73.55)(85,48)
\qbezier[418](85,48)(122.5,74.5)(160,101)
\linethickness{0.45pt}
\qbezier(37,86)(50,65)(77,64)
\qbezier(77,64)(100,63)(133,87)
\put(37,86){\circle*{2.6}}
\put(43,77){\circle*{2.6}}
\put(51,70){\circle*{2.6}}
\put(61,66){\circle*{2.6}}
\put(73,64){\circle*{2.6}}
\put(82,64){\circle*{2.6}}
\put(101,68){\circle*{2.6}}
\put(123,80){\circle*{2.6}}
\put(85,48){\circle*{3.5}}
\put(84,37){\makebox(0,0)[c]{$o$}}
\put(77,115){\makebox(0,0)[c]{$U_{\mathrm{box}}(o)$}}
\put(150,40){\makebox(0,0)[c]{$U(o)$}}
\put(20,124){\makebox(0,0)[c]{$\tau^{-1}(\widetilde\Sigma)$}}
\put(160,116){\makebox(0,0)[c]{$\widetilde\Sigma$}}
\put(188,99){\vector(1,0){54}}
\put(215,116){\makebox(0,0)[c]{$\Phi_o$}}
\qbezier[25](275,103)(275,175)(350,175)
\qbezier[25](350,175)(433,175)(433,103)
\qbezier[25](433,103)(433,25)(350,25)
\qbezier[25](350,25)(275,25)(275,103)
\put(280,45){\vector(1,0){160}}
\put(300,28){\vector(0,1){153}}
\put(443,45){\makebox(0,0)[l]{$Q$}}
\put(300,187){\makebox(0,0)[c]{$P$}}
\qbezier[500](300,45)(355,45)(410,45)
\qbezier[432](410,45)(410,92.5)(410,140)
\qbezier[500](410,140)(355,140)(300,140)
\qbezier[432](300,140)(300,92.5)(300,45)
\put(300,35){\makebox(0,0)[c]{$0$}}
\put(410,33){\makebox(0,0)[c]{$Q_{\mathrm{box}}$}}
\put(288,140){\makebox(0,0)[r]{$P_{\mathrm{box}}$}}
\put(387,122){\makebox(0,0)[c]{$\mathcal B$}}
\put(420,20){\makebox(0,0)[c]{$V$}}
\qbezier[76](312,126.67)(313.5,118.5)(315,110.33)
\qbezier[66](315,110.33)(317,103.46)(319,96.58)
\qbezier[54](319,96.58)(321.5,91.21)(324,85.83)
\qbezier[53](324,85.83)(327.5,81.22)(331,76.61)
\qbezier[58](331,76.61)(336,72.76)(341,68.9)
\qbezier[70](341,68.9)(348,65.86)(355,62.82)
\qbezier[94](355,62.82)(365,60.44)(375,58.07)
\qbezier[92](375,58.07)(385,56.69)(395,55.32)
\put(312,126.67){\circle*{2.7}}
\put(315,110.33){\circle*{2.7}}
\put(319,96.58){\circle*{2.7}}
\put(324,85.83){\circle*{2.7}}
\put(331,76.61){\circle*{2.7}}
\put(341,68.9){\circle*{2.7}}
\put(355,62.82){\circle*{2.7}}
\put(375,58.07){\circle*{2.7}}
\put(395,55.32){\circle*{2.7}}
\put(365,94){\makebox(0,0)[c]{$QP=v>0$}}
\put(365,59){\vector(4,-1){15}}
\put(28,82){\makebox(0,0)[c]{$z$}}
\put(337,129){\makebox(0,0)[c]{$\Phi_o(z)$}}
\put(225,4){\makebox(0,0)[c]{$\Phi_o\circ\widetilde f\circ\Phi_o^{-1}(Q,P)=(e^\nu Q,e^{-\nu}P)$}}
\end{picture}}
\put(225,24){\makebox(0,0)[c]{$\Phi_o(U(o))=V,\qquad \Phi_o(U_{\mathrm{box}}(o))=\mathcal B\Subset V$}}
\put(225,6){\makebox(0,0)[c]{$\Phi_o^*(\mathrm dQ\wedge\mathrm dP)=\mathrm dx\wedge\mathrm dy$}}
\end{picture}%
}
\endgroup

\caption{Local symplectic coordinates. (a) The open saddle neighborhood $U(o)\subset\mathbb R^2_{(x,y)}$ has a dotted boundary; the solid quadrilateral is $U_{\mathrm{box}}(o)$. Dots show a local orbit. (b) The chart $\Phi_o:U(o)\to V\subset\mathbb R^2_{(Q,P)}$ sends this box to the rectangle $\mathcal B$ and the incoming and outgoing branches to the positive $P$ and $Q$ semiaxes. Orbit points satisfy $QP=v>0$.}
\label{fig:local-symplectic-chart}
\end{figure}

	\begin{Remark}
		An \emph{exact symplectic linearization} is a smooth symplectic change of coordinates for which~\eqref{eq:local-linearization} holds as an exact identity, with no nonlinear remainder. Here \emph{exact} refers to the linearization identity; exact symplecticity of a cylinder map is the one-form condition defined in Section~\ref{subsec:auto-1-1}. This is a nontrivial structural assumption and does not follow merely from hyperbolicity in \textbf{(G1)}. The Hartman--Grobman theorem provides only a local $C^0$ topological conjugacy at a general hyperbolic fixed point. The eigenvalue relation $e^\nu \cdot e^{-\nu} = 1$ creates resonances, so the conjugacy need not be smooth and higher-order nonlinear terms cannot in general be removed. Moreover, the chart $\Phi_o$ must satisfy the pullback identity $\Phi_o^*(\mathrm{d}Q \wedge \mathrm{d}P) = \mathrm{d}x \wedge \mathrm{d}y$, making the local linearization compatible with the area-preserving geometry.
	\end{Remark}
	
	To convert local hyperbolic scaling into phase translation, write
	\[
	\mathcal B_{\mathrm{out}}:=(0,Q_{\mathrm{box}}]\times[0,P_{\mathrm{box}}],\qquad
	\mathcal B_{\mathrm{in}}:=[0,Q_{\mathrm{box}}]\times(0,P_{\mathrm{box}}],
	\]
	These sets cover $\mathcal B\setminus\{(0,0)\}$. Restricting the smooth diffeomorphisms on the open half-planes $\{Q>0\}$ and $\{P>0\}$ to $\mathcal B_{\mathrm{out}}$ and $\mathcal B_{\mathrm{in}}$, respectively, defines the \emph{outgoing and incoming logarithmic coordinates}:
	\begin{equation}\label{eq:2.3}
		\begin{aligned}
			\psi_{\mathrm{out}}(Q,P)&=(\theta_{\mathrm{out}},v):=(\nu^{-1}\log Q,\,QP),\\
			\psi_{\mathrm{in}}(Q,P)&=(\theta_{\mathrm{in}},v):=(-\nu^{-1}\log P,\,QP).
		\end{aligned}
	\end{equation}
	Here the phase $\theta$ is the rescaled logarithm of the expanding coordinate $Q$, or the negative rescaled logarithm of the contracting coordinate $P$, and $v=QP$ is invariant under the linearized map. The logarithm converts the multiplicative change of $Q$ or $P$ into an additive phase change.
	The outgoing semiaxis $\{P=0,\,Q>0\}$ and the incoming semiaxis $\{Q=0,\,P>0\}$ are each mapped to the invariant boundary $\{v=0\}$. The inverse maps are $\psi_{\mathrm{out}}^{-1}(\theta,v)=(e^{\nu\theta},ve^{-\nu\theta})$ and $\psi_{\mathrm{in}}^{-1}(\theta,v)=(ve^{\nu\theta},e^{-\nu\theta})$. Composing with the inverse symplectic chart $\Phi_o^{-1}$ defines the inverse local maps from logarithmic coordinates $(\theta, v)$ to $U_{\mathrm{box}}(o)$:
	\begin{equation}\label{eq:log-charts}
		\begin{aligned}
			J_{\mathrm{out}}(\theta,v)
			&:=(\Phi_o^{-1}\circ\psi_{\mathrm{out}}^{-1})(\theta,v)
			=\Phi_o^{-1}(e^{\nu\theta},ve^{-\nu\theta}),\\
			J_{\mathrm{in}}(\theta,v)
			&:=(\Phi_o^{-1}\circ\psi_{\mathrm{in}}^{-1})(\theta,v)
			=\Phi_o^{-1}(ve^{\nu\theta},e^{-\nu\theta}).
		\end{aligned}
	\end{equation}
	In either chart, the lift $\widetilde{f}$ is locally conjugate to the horizontal phase translation $T_1(\theta, v) = (\theta+1, v)$; that is, $J_{\mathrm{out}} \circ T_1 = \widetilde{f} \circ J_{\mathrm{out}}$ and $J_{\mathrm{in}} \circ T_1 = \widetilde{f} \circ J_{\mathrm{in}}$, where in general we write $T_c(\theta, v) := (\theta+c, v)$ for $c \in \mathbb{R}$.
	
	Choose open intervals $B_{\mathrm{out}},B_{\mathrm{in}}\subset\mathbb R$ of length greater than $1$, called the \emph{outgoing basic phase interval} and the \emph{incoming basic phase interval}, respectively. For $h>0$ such that the following maps are defined, $B_{\mathrm{out}}\times[0,h)$ and $B_{\mathrm{in}}\times[0,h)$ are called the outgoing and incoming \emph{basic rectangles}. The set $\mathbb R\times[0,h)$ and its quotient $\mathbb T\times[0,h)$ are one-sided collar neighborhoods, or simply collars, of $v=0$.
	
	In $(Q,P)$ coordinates, the sets $\psi_{\mathrm{out}}^{-1}(B_{\mathrm{out}}\times\{0\})$ and $\psi_{\mathrm{in}}^{-1}(B_{\mathrm{in}}\times\{0\})$ lie on the outgoing and incoming semiaxes, respectively. The images $\psi_{\mathrm{out}}^{-1}(B_{\mathrm{out}}\times[0,h))$ and $\psi_{\mathrm{in}}^{-1}(B_{\mathrm{in}}\times[0,h))$ of the basic rectangles are called the outgoing and incoming \emph{basic strips}. The maps $J_{\mathrm{out}},J_{\mathrm{in}}$ still take values in the original covering space. When no confusion can arise, their images $J_{\mathrm{out}}(B_{\mathrm{out}}\times[0,h))$ and $J_{\mathrm{in}}(B_{\mathrm{in}}\times[0,h))$ in the universal cover under $\Phi_o^{-1}$ are also called the outgoing and incoming basic strips.
	
	We use $0<v_*<v_+$ to specify the domain and range of the transition map; see Lemma~\ref{lem:finite_transition}. Since $T_1$ corresponds to one iterate of $\widetilde{f}$ inside the saddle box, the quotient $\mathbb{T} \times [0, v_*)$ induced by the integer translations $T_1$ is the space of equivalence classes of local orbits passing through this region, that is, the local orbit space. It is distinct from a neighborhood of $\Sigma$ in the original cylinder. Near the corner of the separatrix, a transverse coordinate $\rho=\min(Q,P)$ can be chosen such that the logarithmic coordinate satisfies $v=QP$; see Lemma~\ref{lem:corner-collar}.
	
	The finite-step transition from the outgoing basic rectangle to the incoming basic rectangle at the adjacent saddle, after subtraction of the transition time, gives a normalized transition map $\mathcal{P}$ on $\mathbb{R}\times[0,v_*)$; its explicit construction is given in Section~\ref{sec:transition}. It satisfies
	\[
	\mathcal{P}(\theta + 1, v) = \mathcal{P}(\theta, v) + (1, 0), \qquad \mathcal{P}(\theta, 0) = (b(\theta), 0), \qquad b'(\theta) > 0.
	\]
	
	The \emph{boundary transition map} is the phase map $b(\theta):=\pi_1(\mathcal P(\theta,0))$ obtained by restricting the normalized transition to the invariant boundary $v=0$. It records the incoming logarithmic phase after passage to the adjacent saddle, identification by $\tau^{-1}$, and subtraction of the integer transition time. Since $b(\theta+1)=b(\theta)+1$ and $b'>0$, it is a lift of an orientation-preserving diffeomorphism of $\mathbb T$. Its boundary displacement $b(\theta) - \theta$ is decomposed into its mean $\sigma \in \mathbb{R}$ and a periodic function $\zeta \in C^\infty(\mathbb{T})$ of zero mean:
	\begin{equation}
		b(\theta) = \theta + \sigma + \zeta(\theta), \qquad \sigma := \int_0^1 \big(b(\theta) - \theta\big)\,\mathrm{d}\theta, \qquad \int_0^1 \zeta(\theta)\,\mathrm{d}\theta = 0.
	\end{equation}
	Since $\zeta$ is $1$-periodic with zero mean, each of $\zeta,\zeta',\zeta'',\zeta'''$ has zero mean and hence a zero in each period. The fundamental theorem of calculus gives
	\[
	\begin{gathered}
	\|\zeta^{(j)}\|_\infty\le\|\zeta^{(j+1)}\|_\infty\qquad(0\le j\le3),\\
	\|\zeta\|_{C^4(\mathbb T)}
	=\max_{0\le j\le4}\|\zeta^{(j)}\|_\infty
	=\|\zeta^{(4)}\|_\infty.
	\end{gathered}
	\]
	In particular, the distance from $b$ to the translations, computed for lifts in the maximum norm above, satisfies
    \begin{equation}\label{eq:transition-distance-translations}
        \inf_{c\in\mathbb R}\|b-(\operatorname{id}+c)\|_{C^4(\mathbb T)}
        =\|\zeta\|_{C^4(\mathbb T)}=\|b^{(4)}\|_\infty .
    \end{equation}
    Indeed, the fourth derivative gives a lower bound for every $c$, and $c=\sigma$ attains it. This is the meaning of $C^4$-closeness to a translation in Theorems~\ref{thm:intro-criterion} and~\ref{thm:intro-island}.
	The smallness condition in~\eqref{eq:zeta-C4-assumption} constrains only the nontranslational part $\zeta$ of the boundary displacement, not its mean $\sigma$. The logarithmic coordinates~\eqref{eq:2.3} are part of the normalization. If the outgoing and incoming phase origins are translated by $\theta_{\mathrm{out}}' = \theta_{\mathrm{out}} + a$ and $\theta_{\mathrm{in}}' = \theta_{\mathrm{in}} + c$, with $a, c \in \mathbb{R}$ (including the coordinate scalings $(Q, P) \mapsto (e^{\nu a}Q, e^{-\nu a}P)$ that preserve the symplectic form and the hyperbolic linearization~\eqref{eq:local-linearization}), then the boundary transition map $b(\theta)$ becomes $c + b(\theta - a)$, its mean displacement becomes $\sigma + c - a$, and its zero-mean periodic part $\zeta(\theta)$ becomes the translate $\zeta(\theta - a)$. Its $C^4(\mathbb{T})$ norm $\|\zeta\|_{C^4(\mathbb{T})}$ is unchanged. A general nonlinear reparameterization of the separatrix, however, changes the local passage formula near the saddle and need not preserve this norm. We therefore state the criterion in the logarithmic coordinates~\eqref{eq:2.3}.

\subsection{The general accumulation criterion}\label{subsec:criterion-statements}
\begin{Theorem}\label{thm:boundary-transition-criterion}
		Fix constant-type data \eqref{eq:fix-constant-type}. For every $\nu > 0$, there exist an integer $Q_* \ge 1$ and a positive number $\delta_{\nu, \gamma_\xi} > 0$, depending only on $\nu$ and $\gamma_\xi$, with the following property. If $f$ satisfies (G1)--(G2) and
		\begin{equation}\label{eq:zeta-C4-assumption}
			\|\zeta\|_{C^4(\mathbb{T})} < \delta_{\nu, \gamma_\xi},
		\end{equation}
		
		then there are integers $j$ and $k_f$ and an irrational number $\alpha = (j + \xi)/Q_*$ of constant type such that, for every integer $k\ge k_f$, the map $f$ has a smooth essential invariant circle $\Gamma_k$. Its lift has a smooth parameterization $\gamma_k:\mathbb R\to\mathbb R^2$ satisfying $\gamma_k(t+1)=\tau\bigl(\gamma_k(t)\bigr)$ and, for the chosen lift,
		\begin{equation}\label{eq:rotation-conjugacy}
			\widetilde f\bigl(\gamma_k(t)\bigr) = \gamma_k(t + \omega_k), \qquad \omega_k = \frac{1}{k - \alpha} > 0.
		\end{equation}
		Moreover,
		\begin{equation}\label{eq:omegak-estimate}
			\left\vert{} \omega_k - \frac{p}{q} \right\vert{} \ge \frac{\omega_k}{2 q^2} \qquad (p \in \mathbb{Z},\; q \in \mathbb{N}),
		\end{equation}
		
		and the invariant circle $\Gamma_k$ is disjoint from $\Sigma$ and approaches $\Sigma$ exponentially:
		\begin{equation}\label{eq:Gamma-Sigma-dist}
			\Gamma_k \cap \Sigma = \varnothing, \qquad \sup_{z \in \Gamma_k} \operatorname{dist}(z, \Sigma) \le d_H(\Gamma_k, \Sigma) \le C_f e^{-k\nu/2}.
		\end{equation}
		
		For fixed $\gamma_\xi$, the integer $Q_*$ and threshold $\delta_{\nu, \gamma_\xi}$ can be chosen uniformly when $\nu$ ranges over a compact subset of $(0, \infty)$. The initial index $k_f$ and constant $C_f > 0$ may depend on the map $f$.
	\end{Theorem}
	No twist condition is imposed on $f$, so these invariant circles need not be graphs over $x$.

	\begin{Corollary}\label{cor:f0-invariant}
		Let $f_0$ be a smooth area-preserving cylinder diffeomorphism isotopic to the identity, with a fixed saddle and a lift fixing that saddle. Set $G=f_0^M$ and $\widetilde G=\widetilde f_0^M$, where $M\in\mathbb N$, and suppose that $G$ and this lift satisfy the assumptions of Theorem~\ref{thm:boundary-transition-criterion} at the saddle. Let $k_G$ and $C_G$ denote the threshold and distance constant obtained by applying that theorem to $G$. Then for every $k\ge k_G$, the invariant circle of $G$ obtained from Theorem~\ref{thm:boundary-transition-criterion} is also invariant under $f_0$, and its rotation number with respect to $\widetilde f_0$ is
		\begin{equation}\label{eq:f0-rotation-number}
			\omega_k = \frac{1}{M(k - \alpha)}.
		\end{equation}
		
		With the same threshold $k_G$, the constant-type estimate in~\eqref{eq:omegak-estimate} remains valid; the requirement on $k$ in Lemma~\ref{lem:reciprocal-frequency} is independent of $M$.
	\end{Corollary}
	
	\begin{Corollary}\label{cor:graph-corollary}
		Under the assumptions of Theorem~\ref{thm:boundary-transition-criterion} (with $M=1, G=f_0$) or Corollary~\ref{cor:f0-invariant} (with $M \ge 1, G=f_0^M$), let $\nu>0$ be the logarithm of the expanding hyperbolic multiplier of $G$, and retain the threshold $k_G$ and distance constant $C_G$ from that result. If $f_0$ satisfies the positive twist condition $\partial_y(\pi_x\circ\widetilde{f}_0) \ge b_* > 0$, then for every $k \ge k_G$:
		\begin{enumerate}
			\item $\Gamma_k$ is a smooth graph $\Gamma_k = \operatorname{graph} u_k$ in the original coordinates $(x, y)$, with $u_k \in C^\infty(\mathbb{T})$, and the circle map induced by $f_0$ on this graph is smoothly conjugate to a rigid rotation with rotation number $\omega_k = \frac{1}{M(k - \alpha)}$;
			\item If the separatrix is also a Lipschitz graph $\Sigma = \operatorname{graph}\psi$, with $L_\psi := \operatorname{Lip}(\psi)$, and the first quadrant $\{Q > 0, P > 0\}$ of the local saddle chart lies above $\Sigma$, then
			\begin{equation}\label{eq:graph-distance-estimate}
				0 < u_k(x) - \psi(x) \le (1 + L_\psi) C_G e^{-k\nu/2} \qquad (\forall x \in \mathbb{T}).
			\end{equation}
			
		\end{enumerate}
	\end{Corollary}
	
	\begin{Remark}
		Lemma~\ref{lem:herman-C3} gives an invariant graph for the return map in the $(\theta, s)$ coordinates used for normalization. The orbit extension map $\iota$ and the identification of points related by $\tau$ in the reconstruction of Section~\ref{sec:reconstruction}, however, do not preserve graph structure. Even an exact symplectic coordinate change can destroy it. For a constant $L, \epsilon > 0$, consider the exact symplectic shear
		\[ H_L(\theta, s) = (\theta + Ls, \, s), \qquad H_L^*(s\,\mathrm{d}\theta) - s\,\mathrm{d}\theta = \mathrm{d}\!\left(\frac{1}{2}Ls^2\right).\]
		The image of the smooth graph $s = s_0 + \epsilon\sin(2\pi\theta)$ under $H_L$ is still an embedded essential circle, but the derivative of its first coordinate with respect to $\theta$ is $1 + 2\pi L\epsilon\cos(2\pi\theta)$. If $2\pi L\epsilon > 1$, the first-coordinate projection is no longer injective, so the curve is not a graph over the original angular coordinate. Without a twist assumption on the original map, Theorem~\ref{thm:boundary-transition-criterion} therefore asserts only the existence of smooth embedded essential circles in the original cylinder. Corollaries~\ref{cor:f0-invariant} and~\ref{cor:graph-corollary} have two separate roles, independent of graph structure in the strip coordinates. Invariance under $f_0^M$ does not directly imply invariance under the one-step map $f_0$; Corollary~\ref{cor:f0-invariant} uses the total order of essential circles to recover invariance under $f_0$ and the exact rotation number. Corollary~\ref{cor:graph-corollary} uses the positive twist condition for $f_0$ and Birkhoff's graph theorem to obtain a Lipschitz graph in the original coordinates, and then uses the smooth parameterization to obtain a smooth graph.
	\end{Remark}

\subsection{An island bounded by a heteroclinic cycle}\label{subsec:island-data}

\begin{Definition}\label{def:island-geometry}
	For $G=f_0^M$ and its chosen lift $\widetilde G=\widetilde f_0^M$ as in Theorem~\ref{thm:intro-island}, let $o_0=o$ and $o_1=\tau(o_0)$. In addition to the non-splitting arc $S_+$ from $o_0$ to $o_1$ in \textbf{(G1)}, suppose there is a non-splitting arc $S_-$ oriented by $\widetilde G$ from $o_1$ to $o_0$, whose interior is an unstable branch of $o_1$ and a stable branch of $o_0$. Assume that $S_+$ and $S_-$ meet only at $\{o_0,o_1\}$, so that $S_+\cup S_-$ forms a \emph{heteroclinic cycle} in $\mathbb R^2$ enclosing a bounded open Jordan disk $D$, called an \emph{island}, {
whose integer translates are pairwise disjoint and whose closures meet only at consecutive saddle vertices. More precisely,
\[
\tau^m(D)\cap\tau^n(D)=\varnothing\quad(m\ne n),\qquad
\overline{\tau^n(D)}\cap\overline{\tau^{n+1}(D)}=\{\tau^{n+1}(o_0)\},
\]
and nonadjacent closures are disjoint
}. Write $\Sigma_\pm=\pi(S_\pm)$ with $S_+$ the upper arc. In the compatible symplectic charts of \textbf{(G2)}, we assume that $D$ corresponds to $\{Q>0,\,P<0\}$ near $o_0$ and to $\{Q<0,\,P>0\}$ near $o_1$, while the exterior sectors along $S_\pm$ satisfy $QP>0$.
\end{Definition}

In all four local sectors at $o_0$ and $o_1$, the outgoing and incoming logarithmic phases are defined by $\theta_{\mathrm{out}}=\nu^{-1}\log|Q|$ and $\theta_{\mathrm{in}}=-\nu^{-1}\log|P|$. Let $b_+$ and $b_-$ denote the normalized boundary transition maps of $G$ in the forward direction along $S_+$ and $S_-$, respectively, constructed as in Section~\ref{sec:transition}. Since the translates $\tau^n(D)$ are disjoint, $\pi|_D$ is injective and $\Sigma_\pm$ form two homoclinic loops to $o$ in $\mathbb A$. By an \emph{interior one-sided neighborhood} (or \emph{interior collar}) of $\partial D$ we mean a set of the form $U\cap D$, where $U\subset\mathbb R^2$ is an open neighborhood of the entire Jordan curve $\partial D=S_+\cup S_-$.

The chosen lift $\widetilde f_0$ already preserves $S_+$ and $S_-$: it fixes $o_0$, commutes with $\widetilde G$, and cannot interchange the two unstable half-branches at $o_0$ because their forward $\widetilde G$-orbits converge to distinct fixed lifts $\tau(o_0)$ and $\tau^{-1}(o_0)$. Translation gives the same property at $o_1$, so $\widetilde f_0(D)=D$; in particular, every lifted orbit in $\overline D$ is bounded and has horizontal rotation number zero. Moreover, since $D$ is homeomorphic to $\mathbb R^2$ and $\widetilde f_0|_D$ preserves a finite area of full support, no open disk in $D$ has pairwise disjoint forward iterates, and Brouwer's plane translation theorem \cite{Fra92} guarantees that $\widetilde f_0$ has at least one fixed point $p\in D$.

\begin{Proposition}\label{prop:island}
	Let $f_0$ and $G=f_0^M$ satisfy \textbf{(G1)--(G2)} and Definition~\ref{def:island-geometry}, and fix the arithmetic data~\eqref{eq:fix-constant-type}. There exists a threshold $\delta^{\rm island}_{\nu,\gamma_\xi}>0$, uniform for $\nu$ in compact subintervals of $(0,\infty)$, such that if
	\[
	\max_{j\in\{+,-\}}\inf_{c\in\mathbb R}
	\|b_j-(\operatorname{id}+c)\|_{C^4(\mathbb T)}
	<\delta^{\rm island}_{\nu,\gamma_\xi},
	\]
	then for any fixed point $p\in D$ of $\widetilde f_0$ the following hold. Measure angles in $\mathbb R^2\setminus\{p\}$ about $p$, positively along the boundary loop $S_+,S_-$, and use the angular lift fixing a lift of $o_0$. Then for some constant-type number $\alpha_{\rm in}$ and every sufficiently large $k$, $\widetilde f_0$ has a smooth invariant Jordan curve $C_k\subset D$ surrounding $p$ on which $\widetilde f_0$ is smoothly conjugate to a rigid rotation with rotation number $1/[M(2k-\alpha_{\rm in})]$ about $p$, and $d_H(C_k,\partial D)\le C_G e^{-k\nu/2}$ (the starting index may depend on $p$). Moreover, every interior one-sided neighborhood of $\partial D$ contains smooth quasiperiodic invariant circles surrounding $p$ whose union has positive area, and their projections to $\mathbb A$ are contractible.
\end{Proposition}

\subsection{Herman's hyperbolic family and rotation-number windows}\label{subsec:auto-1-2}
	We now define the hyperbolic perturbation family introduced by Herman in \cite[Chapter~II, \S4]{Her83} and state the quantitative results underlying Theorem~\ref{thm:intro-herman}. The coordinate rescaling, iteration, and verification of the criterion are carried out in Section~\ref{sec:family}.
	
	In \cite[Chapter~II, Proposition~4.2]{Her83} (with parameters $\delta = \lambda^2$ and $\lambda_N = N^{-a/2}$ in his notation), Herman constructed a family of increasing smooth diffeomorphisms $g_\lambda : [0, 1] \to [0, 1]$ ($0 < \lambda \le \lambda_0$) via explicit smooth cutoffs. For all sufficiently small $\lambda$, the map $g_\lambda$ is affine near both endpoints:
	\begin{equation}\label{eq:endpoint-affine-structure}
		\begin{aligned}
			g_\lambda(x) &= (1 + \lambda)x && (0 \le x \le r_0), \\
			g_\lambda(x) &= 1 + \frac{x - 1}{1 + \lambda} && (1 - r_0 \le x \le 1),
		\end{aligned}
	\end{equation}
	where $r_0 \in (0, 1/4)$ is independent of $\lambda$. It satisfies the uniform $C^5$ bound on the closed interval $[0, 1]$:
	\begin{equation}\label{eq:C5-interval-norm}
		\|g_\lambda - \mathrm{id}\|_{C^5([0, 1])} \le C\lambda,
	\end{equation}
	for a constant $C > 0$ independent of $\lambda$. Moreover, for every compact set $K \Subset (0, 1)$, there is a constant $c_K > 0$ independent of $\lambda$ such that
	\begin{equation}\label{eq:positive-drift}
		\min_{x \in K} \bigl(g_\lambda(x) - x\bigr) \ge c_K \lambda.
	\end{equation}
	
	Set $\psi_\lambda = g_\lambda - \mathrm{id}$ and define
	\begin{equation}\label{eq:phi-def}
		\varphi_{\lambda^2}(x) = g_\lambda(x) + g_\lambda^{-1}(x) - 2x.
	\end{equation}
	By the affine endpoint condition~\eqref{eq:endpoint-affine-structure}, the function $\varphi_{\lambda^2}(x)$ equals $\lambda^2x/(1+\lambda)$ on $[0, r_0]$ and $\lambda^2(x-1)/(1+\lambda)$ on $[1 - \frac{r_0}{1+\lambda} , 1]$. Its values at $x = 0$ and $x = 1$ are both zero, and all one-sided derivatives match: the first derivatives equal $\lambda^2/(1+\lambda)$ and all higher derivatives vanish. Hence $\varphi_{\lambda^2}$ defines a smooth function $\varphi_{\lambda^2} \in C^\infty(\mathbb{T})$ on $\mathbb{T}$. {For completeness, the cancellation used in \cite[Chapter~II, Section~4.3]{Her83} follows directly from~\eqref{eq:C5-interval-norm}. For small $\lambda$, $g_\lambda'\ge1/2$. Successively differentiating $g_\lambda\circ g_\lambda^{-1}=\operatorname{id}$ through order five gives a uniform $C^5$ bound for $g_\lambda^{-1}$. The identity $g_\lambda^{-1}-\operatorname{id}=-\psi_\lambda\circ g_\lambda^{-1}$ consequently gives
\[
\|g_\lambda^{-1}-\operatorname{id}\|_{C^5([0,1])}\le C\lambda.
\]
If $z=g_\lambda^{-1}(x)$, then $x=z+\psi_\lambda(z)$ and
\[
\begin{aligned}
\varphi_{\lambda^2}(x)
&=\psi_\lambda(x)-\psi_\lambda(z)\\
&=\psi_\lambda(z)\int_0^1
\psi_\lambda'\bigl(z+t\psi_\lambda(z)\bigr)\,\mathrm dt.
\end{aligned}
\]
The product and chain rules bound the last expression in $C^4$ by $C\|\psi_\lambda\|_{C^5}^2$. The affine endpoint formulas identify these interval estimates with the periodic ones. Therefore there are constants $C_j>0$, independent of $\lambda$, such that, for $0\le j\le4$,}
	\begin{equation}\label{eq:phi-C4-bound}
		\|\varphi_{\lambda^2}\|_{C^j(\mathbb T)} \le C_j\lambda^2.
	\end{equation}
	Define the cylinder map $f_\lambda : \mathbb{A} \to \mathbb{A}$ by its lift fixing the origin:
	\begin{equation}\label{eq:f-original}
		\widetilde{f}_\lambda(x, y) = \bigl(x + y,\ y + \varphi_{\lambda^2}(x + y)\bigr).
	\end{equation}
	This is a smooth area-preserving positive twist map isotopic to the identity. The continuous $1$-periodic extension of $\psi_\lambda$ defines a non-splitting separatrix $\Sigma_\lambda = \operatorname{graph}\psi_\lambda$ of $f_\lambda$, with a corner at the saddle $x \in \mathbb{Z}$ and smooth elsewhere. The multipliers at the saddle are $1 + \lambda$ and $(1 + \lambda)^{-1}$.
	
	After the vertical rescaling $y = \lambda r$, we take $M_\lambda = \lfloor\lambda^{-1}\rfloor$ iterates of the resulting map $F_\lambda$ to obtain $G_\lambda = F_\lambda^{M_\lambda}$. Applying Theorem~\ref{thm:boundary-transition-criterion} and Corollaries~\ref{cor:f0-invariant} and~\ref{cor:graph-corollary} yields the following result for Herman's hyperbolic family.
	
\begin{figure}[tbp]
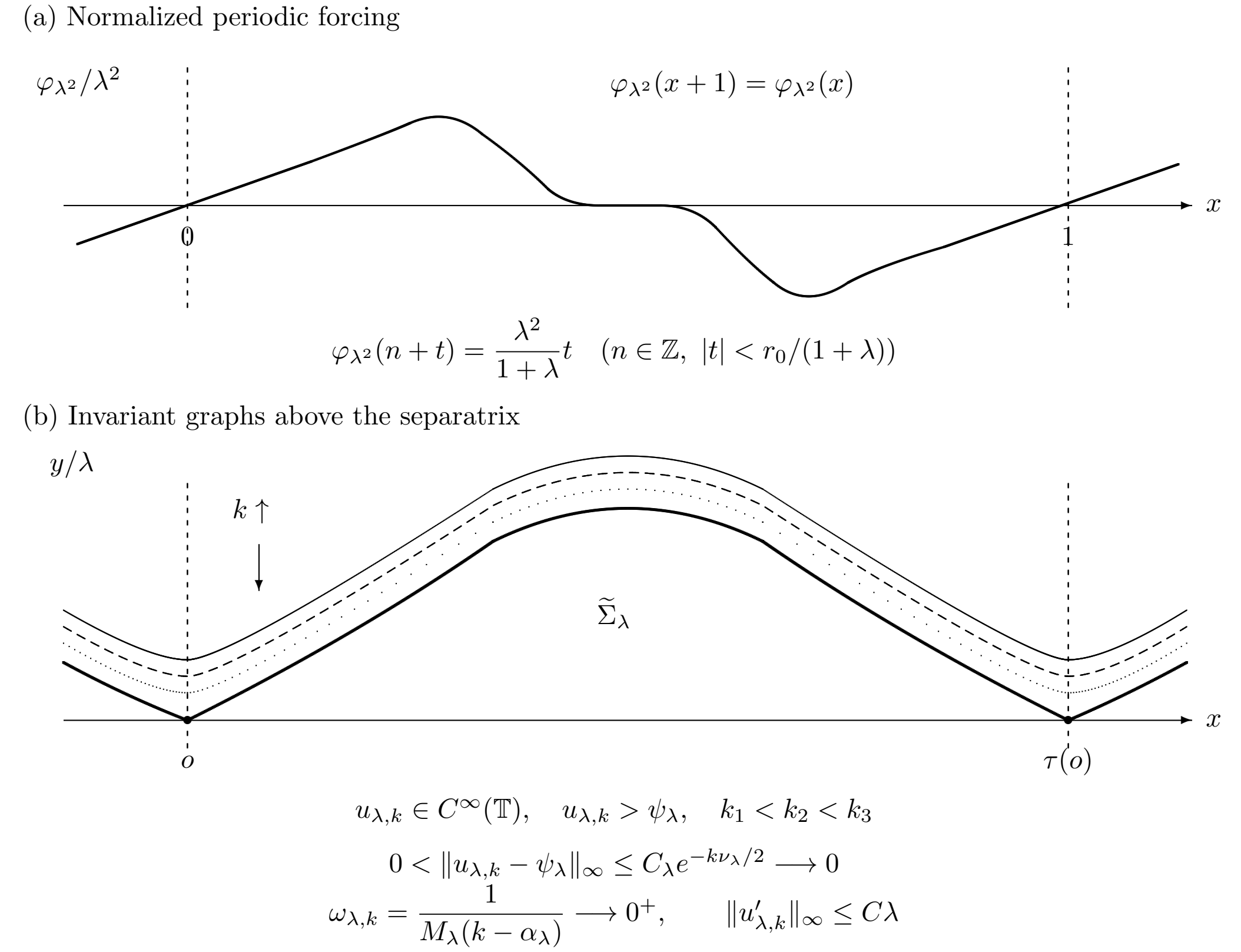

\centering
\begingroup
\setlength{\unitlength}{1pt}
\normalsize\color{black}
\scalebox{0.88}{%
}
\put(225,25){\makebox(0,0)[c]{$0<\|u_{\lambda,k}-\psi_\lambda\|_\infty\le C_\lambda e^{-k\nu_\lambda/2}\longrightarrow0$}}
\put(225,6){\makebox(0,0)[c]{$\omega_{\lambda,k}=\dfrac{1}{M_\lambda(k-\alpha_\lambda)}\longrightarrow0^+,\qquad \|u_{\lambda,k}^\prime\|_\infty\le C\lambda$}}
\end{picture}%
}
\endgroup

\caption{Smooth forcing and invariant graphs (schematic). (a) The normalized forcing is affine near each integer; all endpoint derivatives match, so its periodic extension is smooth. The zero plateau corresponds to the Herman cutoff. (b) At the vertical scale $y/\lambda$, thin solid, dashed, and dotted curves represent $\widetilde\Gamma_{\lambda,k_1}$, $\widetilde\Gamma_{\lambda,k_2}$, and $\widetilde\Gamma_{\lambda,k_3}$, respectively; the thick curve is $\widetilde\Sigma_\lambda$, with corners at $o$ and $\tau(o)$. Increasing $k$ brings the invariant graph closer to the separatrix.}
\label{fig:smoothness-and-accumulation}
\end{figure}

	\begin{Theorem}\label{thm:herman-hyperbolic-family}
		Fix a family of increasing smooth diffeomorphisms $g_\lambda$ satisfying~\eqref{eq:endpoint-affine-structure}--\eqref{eq:positive-drift}. For every fixed constant-type pair $(\xi, \gamma_\xi)$ (see~\eqref{eq:fix-constant-type}), there exist a number $\lambda_* = \lambda_*(\gamma_\xi) > 0$, depending only on $\gamma_\xi$ and the fixed family $g_\lambda$, and a positive integer $Q_* \ge 1$ independent of $\lambda$, such that for each $0 < \lambda < \lambda_*$ there are integers $j_\lambda, k_\lambda$ and a constant-type irrational frequency $\alpha_\lambda = (j_\lambda + \xi)/Q_*$ with the following properties. Set
		\begin{equation}\label{eq:thm1.1-scale-parameters}
			M_\lambda = \lfloor \lambda^{-1} \rfloor, \qquad \nu_\lambda = M_\lambda \log(1 + \lambda).
		\end{equation}
		For every integer $k \ge k_\lambda$, the map $f_\lambda$ has a smooth invariant graph $\Gamma_{\lambda,k} = \operatorname{graph}u_{\lambda,k}$, where $u_{\lambda,k} \in C^\infty(\mathbb{T})$, and its induced circle map is smoothly conjugate to the rigid rotation $R_{\omega_{\lambda,k}}$, where
		\begin{equation}\label{eq:thm1.1-rotation-number}
			\omega_{\lambda,k} = \frac{1}{M_\lambda(k - \alpha_\lambda)} > 0.
		\end{equation}
		This rotation number satisfies the constant-type Diophantine estimate
		\begin{equation}\label{eq:thm1.1-diophantine}
			\left|\omega_{\lambda,k} - \frac{p}{q}\right| \ge \frac{\omega_{\lambda,k}}{2q^2} \qquad (p \in \mathbb{Z},\; q \in \mathbb{N}).
		\end{equation}
		Moreover, the smooth graph function $u_{\lambda,k}$ satisfies
		\begin{align}
			u_{\lambda,k}(x) &> \psi_\lambda(x) \qquad (\forall x \in \mathbb{T}), \label{eq:thm1.1-u-above-psi} \\
			\|u_{\lambda,k}'\|_\infty &\le C\lambda, \label{eq:thm1.1-slope-bound} \\
			0 < \|u_{\lambda,k} - \psi_\lambda\|_\infty &\le C_\lambda e^{-k\nu_\lambda / 2} \longrightarrow 0 \qquad (k \to \infty). \label{eq:thm1.1-decay-rate}
		\end{align}
		The constant $C > 0$ in~\eqref{eq:thm1.1-slope-bound} is independent of both $\lambda$ and $k\ge k_\lambda$. The constant $C_\lambda < \infty$ in~\eqref{eq:thm1.1-decay-rate} is finite for each fixed $\lambda$ but may depend on $\lambda$ (see Figure~\ref{fig:smoothness-and-accumulation}).
	\end{Theorem}
	In \cite[Chapter II, Remark 4.7(3)]{Her83}, Herman wrote:
	\begin{quote}\small
		\emph{je ne sais pas si le diff\'eomorphisme $f_{n,a}$ v\'erifie le Th\'eor\`eme de 4.6}
	\end{quote}
	That is, \textquotedblleft I do not know whether the diffeomorphism $f_{n,a}$ satisfies Theorem 4.6.\textquotedblright\ His subsequent discussion concerns nonzero rotation numbers. In Theorem 4.6, the map modified by the additional perturbation has no essential invariant circle with rotation number in a specified neighborhood of zero. For the unmodified original hyperbolic family $f_{N,a} := f_{\lambda_N}$, where $\lambda_N = N^{-a/2}$, the question is as follows.\medskip
	
	\noindent\textbf{Question (Herman).}
	\emph{For fixed $a > 1$ and $\varepsilon > 0$, is it true that, for every sufficiently large integer $N$, the map $f_{N,a}$ has no essential invariant circle satisfying}
	\begin{equation}\label{eq:herman-question}
		0 < |\omega| < N^{-a/2 - \varepsilon}\,?
	\end{equation}
	Here and below, a rotation-number window of radius $r > 0$ means the set $\{\omega \in \mathbb{R} : 0 < |\omega| < r\}$; after orienting the chosen branch, we use the corresponding half-interval. The answer is \textbf{negative}. By Theorem~\ref{thm:herman-hyperbolic-family}, for every sufficiently large $N$, the map $f_{N,a}$ has infinitely many smooth invariant graphs of constant type whose rotation numbers tend to zero, and hence enter any prescribed punctured neighborhood of zero, including~\eqref{eq:herman-question}. The reason is that, even when the invariant circles spend arbitrarily many iterates near the saddle, the nonconstant boundary displacement in the normalized return map remains small.

	\begin{Theorem}\label{cor:herman-window}
        Fix a family $g_\lambda$ as in Theorem~\ref{thm:herman-hyperbolic-family}, and define $f_\lambda$ by~\eqref{eq:f-original}.
		Given any $a > 0$, set $\lambda_N = N^{-a/2}$ and $f_{N,a} := f_{\lambda_N}$. Fix a pair $(\xi, \gamma_\xi)$ satisfying~\eqref{eq:fix-constant-type}. There exists a positive integer $N_*(a, \gamma_\xi)$ such that, for every integer $N\ge N_*$ and every $r_N>0$, the map $f_{N,a}$ has a smooth invariant graph on which the induced circle map is smoothly conjugate to a rigid rotation $R_{\omega_N}$, where $\omega_N$ is irrational of constant type and satisfies $0<\omega_N<r_N$. In particular, if $a > 1$ and $\varepsilon > 0$, choosing $r_N=N^{-a/2-\varepsilon}$ gives a smooth invariant graph satisfying~\eqref{eq:herman-question}.
	\end{Theorem}

\begin{Remark}\label{rem:scope}
	Here and below, positive area means that the union of smooth quasiperiodic invariant circles contained in the given collar has positive two-dimensional Lebesgue measure. The explicit constant-type sequence $\{\Gamma_k\}$ constructed in Theorem~\ref{thm:boundary-transition-criterion} and Corollary~\ref{cor:f0-invariant} is countable and has measure zero. Given an exterior collar of $\Sigma$, fix $k$ sufficiently large so that $\Gamma_k$ and a closed tubular neighborhood of $\Gamma_k$ lie in the collar. Since $f_0$ preserves area and leaves the essential circle $\Gamma_k$ invariant, it has the intersection property on that tubular neighborhood, and the local theorem of Fayad--Krikorian \cite[Theorem~1; Theorem~3 and Remark~1]{FK09} yields a positive-area family of smooth quasiperiodic invariant circles inside the collar. In Section~\ref{subsec:island-proof}, the same theorem is applied in the {punctured lifted island} $D\setminus\{p\}$ around $C_k$ to obtain the interior positive-area family of Proposition~\ref{prop:island}. We do not assert density one at the boundary or a measure lower bound uniform in $\lambda$.
\end{Remark}

\section{Finite-step transition and the return map}\label{sec:transition}
	We construct a Poincar\'e return lift in a one-sided collar of the separatrix, following the approach of \cite[Sections 4--5]{KK21}. To track the winding number precisely, we work on the covering space and extend by integer translations, retaining explicitly the chosen lift of the map and the number of transition iterates.
	\subsection{Logarithmic coordinates}\label{subsec:auto-3-1}
	Consider the inverse local maps $J_{\mathrm{out}} = \Phi_o^{-1} \circ \psi_{\mathrm{out}}^{-1}$ and $J_{\mathrm{in}} = \Phi_o^{-1} \circ \psi_{\mathrm{in}}^{-1}$ defined in Section~\ref{sec:criterion}. In the local symplectic chart $\Phi_o$, they are represented by the inverse logarithmic charts $\psi_{\mathrm{out}}^{-1}$ and $\psi_{\mathrm{in}}^{-1}$:
	\begin{equation*}
		\begin{aligned}
			\psi_{\mathrm{out}}^{-1}(\theta, v) &\colon \quad Q = e^{\nu\theta}, \quad P = v e^{-\nu\theta};\\
			\psi_{\mathrm{in}}^{-1}(\theta, v) &\colon \quad P = e^{-\nu\theta}, \quad Q = v e^{\nu\theta}.
		\end{aligned}
	\end{equation*}
	The pullbacks of the standard one-form $P\,\mathrm{d}Q$ on the linearizing coordinate space are
	\begin{equation}\label{eq:3.1}
		(\psi_{\mathrm{out}}^{-1})^*(P\,\mathrm{d}Q) = \nu v\,\mathrm{d}\theta_{\mathrm{out}}, \qquad (\psi_{\mathrm{in}}^{-1})^*(P\,\mathrm{d}Q) = \mathrm{d}v + \nu v\,\mathrm{d}\theta_{\mathrm{in}}.
	\end{equation}
	Indeed, direct differentiation gives $\mathrm{d}Q_{\mathrm{out}} = \nu e^{\nu\theta}\,\mathrm{d}\theta, \, P_{\mathrm{out}} = v e^{-\nu\theta}$ and $\mathrm{d}Q_{\mathrm{in}} = e^{\nu\theta}\,\mathrm{d}v + \nu v e^{\nu\theta}\,\mathrm{d}\theta, \, P_{\mathrm{in}} = e^{-\nu\theta}$. Taking exterior derivatives in~\eqref{eq:3.1} and using the symplectic chart identity $(\Phi_o^{-1})^*(\mathrm{d}x \wedge \mathrm{d}y) = \mathrm{d}Q \wedge \mathrm{d}P$, we see that $J_{\mathrm{out}}$ and $J_{\mathrm{in}}$ both pull back the symplectic area form $\mathrm{d}x \wedge \mathrm{d}y$ on the universal cover to the same constant multiple:
	\begin{equation*}
		\begin{aligned}
			J_{\mathrm{out}}^*(\mathrm{d}x \wedge \mathrm{d}y)
			&= (\psi_{\mathrm{out}}^{-1})^*(\mathrm{d}Q \wedge \mathrm{d}P) = \nu\,\mathrm{d}\theta \wedge \mathrm{d}v,\\
			J_{\mathrm{in}}^*(\mathrm{d}x \wedge \mathrm{d}y)
			&= (\psi_{\mathrm{in}}^{-1})^*(\mathrm{d}Q \wedge \mathrm{d}P) = \nu\,\mathrm{d}\theta \wedge \mathrm{d}v.
		\end{aligned}
	\end{equation*}
	Thus the local coordinate change $\mathcal{C}_\nu = \psi_{\mathrm{out}} \circ \psi_{\mathrm{in}}^{-1}$ between the charts, as well as the unnormalized finite-step transition map induced by $\widetilde{f}^L$ below, preserves the area form $\mathrm{d}\theta \wedge \mathrm{d}v$.
	
	For points in the common domain $\mathcal{B}_{\mathrm{out}} \cap \mathcal{B}_{\mathrm{in}} = (0, Q_{\mathrm{box}}] \times (0, P_{\mathrm{box}}]$, that is, $v = QP > 0$, the identity $Q = e^{\nu\theta_{\mathrm{out}}} = v e^{\nu\theta_{\mathrm{in}}}$ gives the local change from incoming to outgoing coordinates:
	\begin{equation}
		\mathcal{C}_\nu := \psi_{\mathrm{out}} \circ \psi_{\mathrm{in}}^{-1} \colon \psi_{\mathrm{in}}(\mathcal{B}_{\mathrm{out}} \cap \mathcal{B}_{\mathrm{in}}) \longrightarrow \psi_{\mathrm{out}}(\mathcal{B}_{\mathrm{out}} \cap \mathcal{B}_{\mathrm{in}}),
	\end{equation}
	Explicitly,
	\begin{equation}
		\mathcal{C}_\nu(\theta, v) = \left( \theta + \frac{\log v}{\nu}, \; v \right).
	\end{equation}
	We continue to write $\mathcal C_\nu$ for the smooth diffeomorphism $\mathbb R\times(0,\infty)\to\mathbb R\times(0,\infty)$ defined by the same formula.
	
	\subsection{The finite-step transition map}\label{subsec:auto-3-2}
	Fix $\epsilon_0 \in (0, 1)$ and choose an outgoing basic phase interval $B_{\mathrm{out}} = (\theta_-, \theta_+)$ satisfying $\vert{}B_{\mathrm{out}}\vert{} = 1 + \epsilon_0 < 2$ and
	\begin{equation}\label{eq:theta_plus_bound}
		e^{\nu(\theta_+ + 1)} < Q_{\mathrm{box}} \quad \Longleftrightarrow \quad \theta_+ < \nu^{-1}\log Q_{\mathrm{box}} - 1.
	\end{equation}
	Its overlap with its unit translate is the nonempty open interval
	\begin{equation*}
		O_{\mathrm{out}} := B_{\mathrm{out}} \cap (B_{\mathrm{out}} - 1) = (\theta_-, \, \theta_+ - 1) \neq \varnothing.
	\end{equation*}
	Taking the closure of the basic interval on the local outgoing semiaxis, set
	\begin{equation*}
		E_{\mathrm{out}} := J_{\mathrm{out}}\big(\overline{B_{\mathrm{out}}} \times \{0\}\big) = \Phi_o^{-1}\big([e^{\nu\theta_-}, e^{\nu\theta_+}] \times \{0\}\big) \subset \widetilde{\Sigma}\setminus\{o,\tau(o)\}.
	\end{equation*}
	By (G1), $\widetilde{\Sigma}\setminus\{o,\tau(o)\} \subset W^s(\tau(o))$, so $E_{\mathrm{out}}$ is also a compact subarc of the stable branch of the adjacent saddle $\tau(o)$.
\begin{Lemma}\label{lem:finite_transition}
		Under \textbf{(G1)}--\textbf{(G2)}, there exist a positive integer $L$, an open interval $B_{\mathrm{in}}=(\theta_{\mathrm{in}}^-,\theta_{\mathrm{in}}^+)$ with $1<|B_{\mathrm{in}}|<2$, and constants $0<v_*<v_+$ with $v_*<e^{-\nu(\theta_{\mathrm{in}}^+-\theta_-)}$ such that $J_{\mathrm{out}}(\overline{B_{\mathrm{out}}}\times[0,v_*])$, $J_{\mathrm{in}}(\overline{B_{\mathrm{in}}}\times[0,v_+])$, and their one-step forward images under $\widetilde f$ are contained in $U_{\mathrm{box}}(o)$, and the unnormalized finite-step transition map
		\begin{equation}\label{eq:P-hat-exact-mapping}
			\widehat{\mathcal P}:=J_{\mathrm{in}}^{-1}\circ\tau^{-1}\circ\widetilde f^L\circ J_{\mathrm{out}}
			:\overline{B_{\mathrm{out}}}\times[0,v_*]
			\longrightarrow B_{\mathrm{in}}\times[0,v_+)
		\end{equation}
		is {
the restriction of a smooth diffeomorphism on an open neighborhood of the closed rectangle onto its image
}. Moreover, $\widehat{\mathcal P}(\theta,0)=\bigl(b^{\mathrm{raw}}(\theta),0\bigr)$, where $(b^{\mathrm{raw}})'>0$, and
		\[
		\widehat{\mathcal P}\bigl(\overline{B_{\mathrm{out}}}\times(0,v_*]\bigr)\subset B_{\mathrm{in}}\times(0,v_+).
		\]
	\end{Lemma}
	\begin{proof}
		At the saddle $\tau(o)$, denote the local stable semiaxis by
		\[
		V^s:=\Phi_{\tau(o)}^{-1}\bigl(\{0\}\times(0,P_{\mathrm{box}})\bigr).
		\]
		By~\eqref{eq:local-linearization}, $P\mapsto e^{-\nu}P$, so $\widetilde f(V^s)\subset V^s$. Since $E_{\mathrm{out}}\subset W^s(\tau(o))$, for each $z\in E_{\mathrm{out}}$ there exists $n(z)\in\mathbb N$ such that $\widetilde f^{n(z)}(z)\in V^s$. By continuity, the same entry time applies on a neighborhood of $z$ in the separatrix. Compactness of $E_{\mathrm{out}}$ gives a finite subcover; let $L$ be the maximum of the corresponding entry times. Forward invariance of $V^s$ implies $\widetilde f^L(E_{\mathrm{out}})\subset V^s$. By compactness, there exist $0<P_-\le P_+<P_{\mathrm{box}}$ such that $\Phi_{\tau(o)}(\widetilde f^L(E_{\mathrm{out}}))\subseteq\{0\}\times[P_-,P_+]$.
		
		Since $E_{\mathrm{out}}$ is mapped into $V^s$, the image of $\overline{B_{\mathrm{out}}}\times\{0\}$ under $\tau^{-1}\circ\widetilde f^L\circ J_{\mathrm{out}}$ lies on the stable semiaxis $\{Q=0,\,P\in[P_-,P_+]\}$ at $o$, where $v=QP=0$. Both the outgoing and incoming phases increase along forward iteration on the branch. Since the derivatives of the charts and iterates are nondegenerate, $b^{\mathrm{raw}}(\theta):=\pi_1\bigl((J_{\mathrm{in}}^{-1}\circ\tau^{-1}\circ\widetilde f^L\circ J_{\mathrm{out}})(\theta,0)\bigr)$ satisfies $(b^{\mathrm{raw}})'>0$. On the overlap, local linearization gives $b^{\mathrm{raw}}(\theta+1)=b^{\mathrm{raw}}(\theta)+1$, so the map extends to the real line by integer translations. Monotonicity and $1<\theta_+-\theta_-<2$ imply
		\[
		1<b^{\mathrm{raw}}(\theta_+)-b^{\mathrm{raw}}(\theta_-)<2.
		\]
		Since $b^{\mathrm{raw}}(\theta_-)\ge-\nu^{-1}\log P_+>-\nu^{-1}\log P_{\mathrm{box}}$, we may choose $\delta_{\mathrm{in}}>0$ sufficiently small so that
		\[
		B_{\mathrm{in}}:=(\theta_{\mathrm{in}}^-,\theta_{\mathrm{in}}^+)
		:=\bigl(b^{\mathrm{raw}}(\theta_-)-\delta_{\mathrm{in}},\,
		b^{\mathrm{raw}}(\theta_+)+\delta_{\mathrm{in}}\bigr)
		\]
		satisfies $e^{-\nu\theta_{\mathrm{in}}^-}<P_{\mathrm{box}}$, $1<|B_{\mathrm{in}}|<2$, and $b^{\mathrm{raw}}(\overline{B_{\mathrm{out}}})\Subset B_{\mathrm{in}}$.
		
		Choose
		\[
		0<v_+<\min\bigl\{P_{\mathrm{box}}e^{\nu\theta_-},\,
		Q_{\mathrm{box}}e^{-\nu(\theta_{\mathrm{in}}^++1)}\bigr\}.
		\]
		By the definitions of $J_{\mathrm{out}}$ and $J_{\mathrm{in}}$, together with~\eqref{eq:theta_plus_bound} and $e^{-\nu\theta_{\mathrm{in}}^-}<P_{\mathrm{box}}$, the images of $\overline{B_{\mathrm{out}}}\times[0,v_+]$ under $J_{\mathrm{out}}$ and of $\overline{B_{\mathrm{in}}}\times[0,v_+]$ under $J_{\mathrm{in}}$, as well as their one-step forward images under $\widetilde f$, lie in $U_{\mathrm{box}}(o)$.
		
		Since $\widehat{\mathcal P}(\theta,0)=\bigl(b^{\mathrm{raw}}(\theta),0\bigr)$ and $b^{\mathrm{raw}}(\overline{B_{\mathrm{out}}})\Subset B_{\mathrm{in}}$, continuity and compactness of $\overline{B_{\mathrm{out}}}$ give $v_0\in(0,v_+]$ such that the image of $\overline{B_{\mathrm{out}}}\times[0,v_0]$ under $\tau^{-1}\circ\widetilde f^L\circ J_{\mathrm{out}}$ remains in the local chart at $o$ and satisfies $0<P<P_{\mathrm{box}}$ and $|Q|<Q_{\mathrm{box}}$, with incoming phase in $B_{\mathrm{in}}$. Since $f$ is isotopic to the identity, it preserves the orientation of the cylinder and the two sides of the invariant separatrix; thus, for $0<v\le v_0$, the image in the incoming chart lies on the chosen side $Q>0$ and satisfies $\widehat{\mathcal P}_2(\theta,v)=QP>0$. Since $\widehat{\mathcal P}_2(\theta,0)=0$, continuity on $\overline{B_{\mathrm{out}}}\times[0,v_0]$ allows us to choose
		\[
		0<v_*<\min\bigl\{v_0,\,v_+,\,e^{-\nu(\theta_{\mathrm{in}}^+-\theta_-)}\bigr\}
		\]
		sufficiently small so that $0<\widehat{\mathcal P}_2(\theta,v)<v_+$ on $\overline{B_{\mathrm{out}}}\times(0,v_*]$. Thus $\widehat{\mathcal P}$ maps $\overline{B_{\mathrm{out}}}\times[0,v_*]$ into $B_{\mathrm{in}}\times[0,v_+)$ and $\overline{B_{\mathrm{out}}}\times(0,v_*]$ into $B_{\mathrm{in}}\times(0,v_+)$. {
The same composition is defined on an open neighborhood of the closed rectangle: the logarithmic charts extend across $v=0$ while $Q>0$ in the outgoing chart and $P>0$ in the incoming chart, and compactness gives a common neighborhood on which all chart inclusions persist. Each factor is a diffeomorphism onto its image. Thus the restriction has a smooth local inverse also at $v=0$, as required below.
}
	\end{proof}
	
	On the overlap $O_{\mathrm{out}} \times [0, v_*)$, the identities $J_{\mathrm{out}} \circ T_1 = \widetilde{f} \circ J_{\mathrm{out}}$ and $J_{\mathrm{in}} \circ T_1 = \widetilde{f} \circ J_{\mathrm{in}}$ give
	\begin{equation}\label{eq:transition-overlap-identity}
		\begin{aligned}
			\widehat{\mathcal P}\circ T_1
			&=J_{\mathrm{in}}^{-1}\circ\tau^{-1}\circ\widetilde f^L\circ J_{\mathrm{out}}\circ T_1\\
			&=J_{\mathrm{in}}^{-1}\circ\widetilde f\circ\tau^{-1}\circ\widetilde f^L\circ J_{\mathrm{out}}
			=T_1\circ\widehat{\mathcal P}.
		\end{aligned}
	\end{equation}
	
	Set $\mathcal{U}_n := (B_{\mathrm{out}} + n) \times [0, v_*)$. Since $|B_{\mathrm{out}}| > 1$, these sets cover $\mathbb{R}\times[0,v_*)$. Define the local maps $\widehat{\mathcal{P}}^{(n)} := T_n \circ \widehat{\mathcal{P}} \circ T_{-n}$. By~\eqref{eq:transition-overlap-identity}, the maps $\widehat{\mathcal P}^{(n)}$ and $\widehat{\mathcal P}^{(n+1)}$ agree on $\mathcal{U}_n \cap \mathcal{U}_{n+1}$. The equality holds on open overlaps, so all derivatives agree there as well. Hence $\widehat{\mathcal{P}}$ extends uniquely to a smooth map on $\mathbb{R} \times [0, v_*)$ satisfying, for every $n \in \mathbb{Z}$,
	\begin{equation*}
		\widehat{\mathcal{P}}(\theta + n, v) = T_n\bigl(\widehat{\mathcal P}(\theta,v)\bigr) \qquad (n \in \mathbb{Z}).
	\end{equation*}
	
	Both charts pull back the area form to the same constant multiple, and $\widetilde f$ preserves area, so $\det D\widehat{\mathcal P}=1$. By Lemma~\ref{lem:finite_transition}, the map $\widehat{\mathcal P}$ is a local diffeomorphism at every point. After reducing $v_*$, it is also injective. Suppose, to the contrary, that no $v_*$ makes $\widehat{\mathcal{P}}$ globally injective. Then for every $k \in \mathbb{N}$ there are distinct points $z_k = (\theta_k, v_k) \neq w_k = (\theta'_k, v'_k)$ satisfying $0 \le v_k, v'_k < 1/k$ and $\widehat{\mathcal{P}}(z_k) = \widehat{\mathcal{P}}(w_k)$. By $\widehat{\mathcal{P}} \circ T_n = T_n \circ \widehat{\mathcal{P}}$, a common integer phase translation allows us to assume $\theta_k \in [0, 1]$. As $k \to \infty$, we have $v_k, v'_k \to 0$. Compactness of $[0, 1]$ gives a convergent subsequence with $\theta_k \to \theta_\infty \in [0, 1]$. Since $\widehat{\mathcal P}_1(\theta,v)-b^{\mathrm{raw}}(\theta)$ is $1$-periodic in $\theta$, continuity implies that it tends to zero uniformly in $\theta\in\mathbb R$ as $v\to0$. Equality of the images then yields $b^{\mathrm{raw}}(\theta'_k) \to b^{\mathrm{raw}}(\theta_\infty)$. Global monotonicity and injectivity of $b^{\mathrm{raw}}$ on $\mathbb{R}$ imply $\theta'_k \to \theta_\infty$, so both $z_k$ and $w_k$ converge to $z_\infty := (\theta_\infty, 0)$. The inverse function theorem, however, makes $\widehat{\mathcal{P}}$ injective in a neighborhood of $z_\infty$. For sufficiently large $k$, both points lie in that neighborhood, forcing $z_k=w_k$, a contradiction. Thus $\widehat{\mathcal{P}}$ is a global smooth embedding.
	
	The definition of the unnormalized finite-step transition map $\widehat{\mathcal{P}}_L := J_{\mathrm{in}}^{-1} \circ \tau^{-1} \circ \widetilde{f}^L \circ J_{\mathrm{out}}$ explicitly depends on the entry time $L \in \mathbb{N}$. If we increase this time by a positive integer $k \in \mathbb{N}$, replacing $L$ by $L' = L + k$ and reducing $v_*$ if necessary, then the identity $J_{\mathrm{in}}^{-1}\circ\widetilde f^k=T_k\circ J_{\mathrm{in}}^{-1}$ gives, on the common domain,
	\begin{equation}
		\widehat{\mathcal{P}}_{L+k} = J_{\mathrm{in}}^{-1} \circ \tau^{-1} \circ \widetilde{f}^{L+k} \circ J_{\mathrm{out}} = \big(J_{\mathrm{in}}^{-1} \circ \widetilde{f}^k\big) \circ \tau^{-1} \circ \widetilde{f}^L \circ J_{\mathrm{out}} = T_k \circ \widehat{\mathcal{P}}_L
	\end{equation}
	To remove the integer phase translation corresponding to the entry time, define the normalized transition map
	\begin{equation}\label{eq:normalized-transition-map}
		\mathcal{P} := T_{-L} \circ \widehat{\mathcal{P}}_L \colon \mathbb{R} \times [0, v_*) \longrightarrow \mathbb{R} \times [0, v_+).
	\end{equation}
	For every admissible $L' = L + k$,
	\begin{equation*}
		T_{-(L+k)} \circ \widehat{\mathcal{P}}_{L+k} = (T_{-L} \circ T_{-k}) \circ (T_k \circ \widehat{\mathcal{P}}_L) = T_{-L} \circ \widehat{\mathcal{P}}_L
	\end{equation*}
	Thus $\mathcal P$ is unchanged on the common domain when the entry time is increased. The phase translation $T_1$ corresponds to one local iterate, and $T_{-L}$ induces the identity on the quotient annulus; see \cite[Sections 4--5]{KK21}. This normalization leaves the action of the return map on the quotient annulus unchanged and removes only the integer displacement of the lift.
	
	In particular, restricting to the invariant boundary $v = 0$ gives the following explicit expression for the boundary transition map introduced in Section~\ref{sec:criterion}:
	\begin{equation*}
		b(\theta) := \pi_1(\mathcal{P}(\theta, 0)) = \pi_1\Big((J_{\mathrm{in}}^{-1}\circ\tau^{-1}\circ\widetilde{f}^L\circ J_{\mathrm{out}})(\theta,0)\Big) - L
	\end{equation*}
	Here $\pi_1$ denotes projection onto the first logarithmic phase coordinate. The construction gives
	\begin{equation*}
		b(\theta + 1) = b(\theta) + 1, \qquad b'(\theta) > 0, \qquad \mathcal{P}_2(\theta, 0) = 0.
	\end{equation*}
	When $f$ is an iterate $G=f_0^M$ and the one-step restriction of $\widetilde f_0$ to the separatrix is close to a translation in a suitable coordinate, the $C^4$ distance of $b$ from the translations can be bounded over $n=ML$ steps as follows.
	
	\begin{Lemma}\label{lem:one-step-transition}
		Let $G=f_0^M$ ($M\in\mathbb N$) satisfy \textbf{(G1)--(G2)}, with chosen lift $\widetilde G=\widetilde f_0^M$ and boundary transition map $b_G$. Suppose that $\widetilde f_0$ preserves $\widetilde\Sigma$, and let $g$ be its restriction to $\widetilde\Sigma\setminus\{o,\tau(o)\}$. For a closed interval $I\subset B_{\mathrm{out}}$ of length one, set $X:=\{J_{\mathrm{out}}(\theta,0):\theta\in I\}$, and fix $L\in\mathbb N$ as in Lemma~\ref{lem:finite_transition} so that $g^n(X)$ lies in the incoming chart at $\tau(o)$ with $n:=ML$.
		
		Suppose an increasing $C^4$ coordinate $t=h(x)$, defined on an open subarc containing $\bigcup_{j=0}^n g^j(X)$, agrees with $a_-\theta_{\mathrm{out}}+c_-$ near $X$ and with $a_+\theta_{\mathrm{in}}+c_+$ near $g^n(X)$ (where $\theta_{\mathrm{in}}$ is measured at $\tau(o)$) for constants $a_-,a_+>0$ and $c_-,c_+\in\mathbb R$. Write $H:=h\circ g\circ h^{-1}$, and assume that for some $d\in\mathbb R$ and $\varepsilon>0$,
		\[
		H(t)=t+d+e(t),\qquad
		\max_{0\le j<n}\|e\|_{C^4(h(g^j(X)))}\le\varepsilon,\qquad
		n\varepsilon\le1.
		\]
		Then there exists a numerical constant $C_4>0$ such that
		\begin{equation}\label{eq:one-step-transition}
			\inf_{c\in\mathbb R}\|b_G-(\operatorname{id}+c)\|_{C^4(\mathbb T)}
			\le C_4\,\frac{a_-^4}{a_+}\,n\varepsilon .
		\end{equation}
	\end{Lemma}
	
	\begin{proof}
		Set $U_j:=H^j$ on $h(X)$. By the chain rule, $\|U_j'\|_\infty\le(1+\varepsilon)^j\le e$ for $0\le j\le n$. For $2\le r\le4$, differentiating $U_{j+1}=U_j+d+e\circ U_j$ via~\eqref{eq:faa-di-bruno} gives a recurrence in which the leading term is $(1+e'\circ U_j)U_j^{(r)}$ and all remaining terms are products of derivatives of $e$ with polynomials in $U_j',\dots,U_j^{(r-1)}$. By induction on $r$, $D_{j,r}:=\|U_j^{(r)}\|_\infty$ satisfies
		\[
		D_{j+1,r}\le(1+\varepsilon)D_{j,r}+C_r\varepsilon,\qquad D_{0,r}=0,
		\]
		whence $D_{j,r}\le C_r j\varepsilon e^{j\varepsilon}\le C_r e\,n\varepsilon$ for $0\le j\le n$ and $2\le r\le 4$. On $I$, the endpoint expressions for $h$ and the definition of $b_G$ yield
		\[
		a_+\bigl(b_G(\theta)+L\bigr)+c_+ = H^n(a_-\theta+c_-).
		\]
		Differentiating four times with respect to $\theta$ gives $\|b_G^{(4)}\|_{L^\infty(I)}\le (a_-^4/a_+)D_{n,4}$. Since $I$ has length one, \eqref{eq:transition-distance-translations} yields~\eqref{eq:one-step-transition}.
	\end{proof}
	\subsection{The Poincar\'e return lift}\label{subsec:auto-3-3}
	Define the outgoing basic strip and the incoming target region in $U_{\mathrm{box}}(o)$ by
	\begin{equation}\label{eq:Sigma-out-window}
		\Sigma_{\mathrm{out}}
		:=J_{\mathrm{out}}\bigl(B_{\mathrm{out}}\times(0,v_*)\bigr),\qquad
		D_{\mathrm{in}}
		:=J_{\mathrm{in}}\bigl(B_{\mathrm{in}}\times(0,v_+)\bigr)
		\subset U_{\mathrm{box}}(o),
	\end{equation}
	The corresponding incoming region at the adjacent saddle is $\tau(D_{\mathrm{in}})\subset U_{\mathrm{box}}(\tau(o))$. For an initial point $z_0=J_{\mathrm{in}}(\theta,v)$ in the subdomain $J_{\mathrm{in}}\bigl(B_{\mathrm{in}}\times(0,v_*)\bigr)\subset D_{\mathrm{in}}$, where $\theta\in B_{\mathrm{in}}$ and $0<v<v_*$, as long as the orbit remains in the linearizing box,
	\begin{equation}\label{eq:local-orbit-evolution}
		\widetilde f^j(z_0)=\Phi_o^{-1}(ve^{\nu(\theta+j)},e^{-\nu(\theta+j)}).
	\end{equation}
	By Lemma~\ref{lem:finite_transition}, $\theta+\nu^{-1}\log v<\theta_-$. Set
	\begin{equation}\label{eq:discrete-N-exact}
		N(\theta,v):=\left\lfloor\theta_- - \theta-\frac{\log v}{\nu}\right\rfloor+1\ge1.
	\end{equation}
	Then $N(\theta,v)$ is the smallest positive integer for which the outgoing phase exceeds $\theta_-$, and
	\[
	\theta_-<\theta_{\mathrm{out}}:=\theta+\nu^{-1}\log v+N(\theta,v)
	\le\theta_-+1<\theta_+.
	\]
	For $0\le j\le N(\theta,v)$, the vertical coordinate decreases and is bounded above by its initial value, while the horizontal coordinate increases and is bounded above by its terminal value $e^{\nu(\theta_-+1)}$. Since both endpoints lie in the box, the local linearization formula shows inductively that the whole orbit segment remains in the box. For $j<N(\theta,v)$ the outgoing phase is at most $\theta_-$, so these points do not belong to $\Sigma_{\mathrm{out}}$. Hence $N(\theta,v)$ is the first entry time into $\Sigma_{\mathrm{out}}$.
	
	For any $(\theta, v) \in B_{\mathrm{in}} \times (0, v_*)$, let $Q_0 = v e^{\nu\theta}$ denote the $Q$-coordinate of the initial point $z_0 = J_{\mathrm{in}}(\theta, v)$ in the $(Q, P)$ chart. After $N = N(\theta, v)$ iterates, the $Q$-coordinate of $\widetilde f^N(z_0) \in \Sigma_{\mathrm{out}}$ is given by
	\[
	Q_N = e^{\nu\theta_{\mathrm{out}}} = v e^{\nu(\theta + N)},
	\]
	while $Q_N P_N = v$. Hence the outgoing logarithmic phase at the entry time $N(\theta, v)$ is $\theta_{\mathrm{out}} = \theta + \nu^{-1}\log v + N(\theta, v)$, and $\widetilde f^{N(\theta, v)}\bigl(J_{\mathrm{in}}(\theta, v)\bigr) = J_{\mathrm{out}}(\theta_{\mathrm{out}}, v)$. In logarithmic coordinates, this local passage from the incoming strip into $\Sigma_{\mathrm{out}}$ is therefore represented by the map
	\begin{equation}\label{eq:C-hat-nu-def}
		\begin{gathered}
			\begin{aligned}
				\widehat{\mathcal{C}}_\nu(\theta, v)
				&:= \left(\theta + \frac{\log v}{\nu} + N(\theta, v),\; v\right)\\
				&= T_{N(\theta, v)}\bigl(\mathcal{C}_\nu(\theta, v)\bigr),
			\end{aligned}\\
			(\theta, v) \in B_{\mathrm{in}} \times (0, v_*).
		\end{gathered}
	\end{equation}
	The definition of $N(\theta, v)$ gives $\widehat{\mathcal{C}}_\nu(B_{\mathrm{in}} \times (0, v_*)) \subset B_{\mathrm{out}} \times (0, v_*)$. By Lemma~\ref{lem:finite_transition}, the composition of the unnormalized transition map $\widehat{\mathcal{P}}$ with $\widehat{\mathcal{C}}_\nu$ is well defined; denote it by
	\begin{equation}\label{eq:closed-loop-return}
		\widehat{\mathcal R}_{\mathrm{raw}}
		:=\widehat{\mathcal P}\circ\widehat{\mathcal C}_\nu:
		B_{\mathrm{in}}\times(0,v_*)\longrightarrow B_{\mathrm{in}}\times(0,v_+),
	\end{equation}
	Then
	\[
	(\tau\circ J_{\mathrm{in}}\circ\widehat{\mathcal R}_{\mathrm{raw}})(\theta,v)
	=\widetilde f^{N(\theta,v)+L}\bigl(J_{\mathrm{in}}(\theta,v)\bigr).
	\]
	Here $N+L$ is the time needed to reach the incoming basic strip at the adjacent saddle along the prescribed route. Since $\vert{}B_{\mathrm{in}}\vert{} > 1$, the overlap $O_{\mathrm{in}} = B_{\mathrm{in}} \cap (B_{\mathrm{in}} - 1)$ is nonempty. For $\theta \in O_{\mathrm{in}}$, the one-step image $\widetilde{f}(J_{\mathrm{in}}(\theta, v)) = J_{\mathrm{in}}(\theta + 1, v)$ still lies in $J_{\mathrm{in}}(B_{\mathrm{in}} \times (0, v_*))$, so $N(\theta, v) + L$ is not the first return time to the incoming basic strip. Moreover, $\widehat{\mathcal{R}}_{\mathrm{raw}}$ takes values in $B_{\mathrm{in}} \times (0, v_+)$ and need not map its domain $B_{\mathrm{in}} \times (0, v_*)$ into itself. Figure~\ref{fig:raw-return-lift} shows this unnormalized return passage in the universal cover and its coordinate representation.
	\begin{figure}[tbp]
\centering
\begingroup
\setlength{\unitlength}{1pt}
\normalsize\color{black}
\scalebox{0.88}{%
\begin{picture}(450,366)
\linethickness{0.45pt}
\put(0,51){\begin{picture}(450,0)
\put(10,300){\makebox(0,0)[l]{(a) Return passage in the universal cover}}
\linethickness{0.9pt}
\qbezier(15,211)(45,193)(72,173)
\qbezier(72,173)(135,208)(187,244)
\qbezier(187,244)(225,263)(263,244)
\qbezier(263,244)(320,208)(378,173)
\qbezier(378,173)(406,191)(438,211)
\linethickness{0.45pt}
\qbezier(39,214)(55,189)(78,189)
\qbezier(78,189)(116,198)(186,260)
\qbezier(186,260)(225,286)(264,260)
\qbezier(264,260)(329,205)(365,195)
\put(44.99,205.81){\circle*{2.7}}
\put(57.53,194.76){\circle*{2.7}}
\put(69.95,189.81){\circle*{2.7}}
\put(92.72,193.96){\circle*{2.7}}
\put(126.17,212.69){\circle*{2.7}}
\put(161.84,239.4){\circle*{2.7}}
\put(200.04,267.68){\circle*{2.7}}
\put(226.56,272.98){\circle*{2.7}}
\put(249.96,267.68){\circle*{2.7}}
\put(286.46,241.66){\circle*{2.7}}
\put(323.76,214.97){\circle*{2.7}}
\put(351.1,200.06){\circle*{2.7}}
\put(72,173){\circle*{3.3}}
\put(378,173){\circle*{3.3}}
\put(72,159){\makebox(0,0)[c]{$o$}}
\put(378,159){\makebox(0,0)[c]{$\tau(o)$}}
\qbezier[18](20,211)(46,192)(72,173)
\qbezier[18](72,173)(99,192)(126,211)
\qbezier[18](126,211)(99,230)(72,249)
\qbezier[18](72,249)(46,230)(20,211)
\qbezier[18](326,211)(352,192)(378,173)
\qbezier[18](378,173)(404.5,192)(431,211)
\qbezier[18](431,211)(404.5,230)(378,249)
\qbezier[18](378,249)(352,230)(326,211)
\put(16,259){\makebox(0,0)[l]{$U_{\mathrm{box}}(o)$}}
\put(433,259){\makebox(0,0)[r]{$U_{\mathrm{box}}(\tau(o))$}}
\put(20,230){\makebox(0,0)[l]{$D_{\mathrm{in}}$}}
\put(120,222){\makebox(0,0)[c]{$\Sigma_{\mathrm{out}}$}}
\put(418,228){\makebox(0,0)[c]{$\tau(D_{\mathrm{in}})$}}
\put(39,214){\circle*{3}}
\put(78,189){\circle*{3}}
\put(365,195){\circle*{3}}
\put(29,209){\makebox(0,0)[r]{$z_0$}}
\put(93,183){\makebox(0,0)[c]{$z_N$}}
\put(371,191){\makebox(0,0)[l]{$z_{N+L}$}}
\put(226,229){\makebox(0,0)[c]{$\widetilde\Sigma$}}
\put(100,200){\vector(2,1){15}}
\put(143,188){\makebox(0,0)[c]{$\widetilde f^{\,N(\theta,v)}$}}
\put(282,244){\vector(2,-1){15}}
\put(310,258){\makebox(0,0)[c]{$\widetilde f^{\,L}$}}
\qbezier[14](365,202)(365,203.5)(365,205)
\qbezier[14](365,208)(365,209.5)(365,211)
\qbezier[14](365,214)(365,215.5)(365,217)
\qbezier[14](365,220)(365,221.5)(365,223)
\qbezier[14](365,226)(365,227.5)(365,229)
\qbezier[14](365,232)(365,233.5)(365,235)
\qbezier[14](365,238)(365,239.5)(365,241)
\qbezier[14](365,244)(365,245.5)(365,247)
\qbezier[14](365,250)(365,251.5)(365,253)
\qbezier[14](365,256)(365,257.5)(365,259)
\qbezier[14](365,262)(365,263.5)(365,265)
\qbezier[14](365,268)(365,269.5)(365,271)
\qbezier[14](365,274)(365,275.5)(365,277)
\qbezier[14](365,280)(365,281.5)(365,283)
\qbezier[14](365,283)(363.5,283)(362,283)
\qbezier[14](359,283)(357.5,283)(356,283)
\qbezier[14](353,283)(351.5,283)(350,283)
\qbezier[14](347,283)(345.5,283)(344,283)
\qbezier[14](341,283)(339.5,283)(338,283)
\qbezier[14](335,283)(333.5,283)(332,283)
\qbezier[14](329,283)(327.5,283)(326,283)
\qbezier[14](323,283)(321.5,283)(320,283)
\qbezier[14](317,283)(315.5,283)(314,283)
\qbezier[14](311,283)(309.5,283)(308,283)
\qbezier[14](305,283)(303.5,283)(302,283)
\qbezier[14](299,283)(297.5,283)(296,283)
\qbezier[14](293,283)(291.5,283)(290,283)
\qbezier[14](287,283)(285.5,283)(284,283)
\qbezier[14](281,283)(279.5,283)(278,283)
\qbezier[14](275,283)(273.5,283)(272,283)
\qbezier[14](269,283)(267.5,283)(266,283)
\qbezier[14](263,283)(261.5,283)(260,283)
\qbezier[14](257,283)(255.5,283)(254,283)
\qbezier[14](251,283)(249.5,283)(248,283)
\qbezier[14](245,283)(243.5,283)(242,283)
\qbezier[14](239,283)(237.5,283)(236,283)
\qbezier[14](233,283)(231.5,283)(230,283)
\qbezier[14](227,283)(225.5,283)(224,283)
\qbezier[14](221,283)(219.5,283)(218,283)
\qbezier[14](215,283)(213.5,283)(212,283)
\qbezier[14](209,283)(207.5,283)(206,283)
\qbezier[14](203,283)(201.5,283)(200,283)
\qbezier[14](197,283)(195.5,283)(194,283)
\qbezier[14](191,283)(189.5,283)(188,283)
\qbezier[14](185,283)(183.5,283)(182,283)
\qbezier[14](179,283)(177.5,283)(176,283)
\qbezier[14](173,283)(171.5,283)(170,283)
\qbezier[14](167,283)(165.5,283)(164,283)
\qbezier[14](161,283)(159.5,283)(158,283)
\qbezier[14](155,283)(153.5,283)(152,283)
\qbezier[14](149,283)(147.5,283)(146,283)
\qbezier[14](143,283)(141.5,283)(140,283)
\qbezier[14](137,283)(135.5,283)(134,283)
\qbezier[14](131,283)(129.5,283)(128,283)
\qbezier[14](125,283)(123.5,283)(122,283)
\qbezier[14](119,283)(117.5,283)(116,283)
\qbezier[14](113,283)(111.5,283)(110,283)
\qbezier[14](107,283)(105.5,283)(104,283)
\qbezier[14](101,283)(99.5,283)(98,283)
\qbezier[14](95,283)(93.5,283)(92,283)
\qbezier[14](89,283)(87.5,283)(86,283)
\qbezier[14](83,283)(81.5,283)(80,283)
\qbezier[14](77,283)(75.5,283)(74,283)
\qbezier[14](71,283)(69.5,283)(68,283)
\qbezier[14](65,283)(63.5,283)(62,283)
\qbezier[14](59,283)(59,281.5)(59,280)
\qbezier[14](59,277)(59,275.5)(59,274)
\qbezier[14](59,271)(59,269.5)(59,268)
\qbezier[14](59,265)(59,263.5)(59,262)
\qbezier[14](59,259)(59,257.5)(59,256)
\qbezier[14](59,253)(59,251.5)(59,250)
\qbezier[14](59,247)(59,245.5)(59,244)
\qbezier[14](59,241)(59,239.5)(59,238)
\qbezier[14](59,235)(59,233.5)(59,232)
\qbezier[14](59,229)(59,227.5)(59,226)
\qbezier[14](59,223)(59,221.5)(59,220)
\qbezier[14](59,217)(59,215.5)(59,214)
\qbezier[14](59,211)(59,209.5)(59,208)
\qbezier[5](59,205)(59,204.5)(59,204)
\put(59,207){\vector(0,-1){10}}
\put(225,294){\makebox(0,0)[c]{$\tau^{-1}$}}
\put(59,195){\circle*{3}}
\put(51,200){\makebox(0,0)[r]{$\bar z$}}
\put(27,154){\makebox(0,0)[c]{$\tau^{-1}(\widetilde\Sigma)$}}
\put(426,154){\makebox(0,0)[c]{$\tau(\widetilde\Sigma)$}}
\put(10,139){\makebox(0,0)[l]{(b) The corresponding coordinate maps}}
\put(55,110){\makebox(0,0)[c]{$(\theta,v)$}}
\put(225,110){\makebox(0,0)[c]{$(\theta_{\mathrm{out}},v)$}}
\put(402,110){\makebox(0,0)[c]{$(\theta_1,v_1)$}}
\put(100,110){\vector(1,0){76}}
\put(279,110){\vector(1,0){77}}
\put(139,124){\makebox(0,0)[c]{$\widehat{\mathcal C}_\nu$}}
\put(318,124){\makebox(0,0)[c]{$\widehat{\mathcal P}$}}
\put(55,96){\vector(0,-1){34}}
\put(225,96){\vector(0,-1){34}}
\put(402,96){\vector(0,-1){34}}
\put(43,79){\makebox(0,0)[r]{$J_{\mathrm{in}}$}}
\put(213,79){\makebox(0,0)[r]{$J_{\mathrm{out}}$}}
\put(390,79){\makebox(0,0)[r]{$J_{\mathrm{in}}$}}
\put(55,49){\makebox(0,0)[c]{$z_0$}}
\put(225,49){\makebox(0,0)[c]{$z_N$}}
\put(402,49){\makebox(0,0)[c]{$\bar z=\tau^{-1}z_{N+L}$}}
\put(100,49){\vector(1,0){76}}
\put(279,49){\vector(1,0){77}}
\put(139,35){\makebox(0,0)[c]{$\widetilde f^{\,N}$}}
\put(318,35){\makebox(0,0)[c]{$\tau^{-1}\circ\widetilde f^{\,L}$}}
\put(225,8){\makebox(0,0)[c]{$\widehat{\mathcal R}_{\mathrm{raw}}=\widehat{\mathcal P}\circ\widehat{\mathcal C}_\nu$}}
\end{picture}}
\put(225,40){\makebox(0,0)[c]{$\theta_{\mathrm{out}}=\theta+\nu^{-1}\log v+N(\theta,v)$}}
\put(225,23){\makebox(0,0)[c]{$\widehat{\mathcal P}=J_{\mathrm{in}}^{-1}\circ\tau^{-1}\circ\widetilde f^{\,L}\circ J_{\mathrm{out}}$}}
\put(225,6){\makebox(0,0)[c]{$D_{\mathrm{in}}=J_{\mathrm{in}}(B_{\mathrm{in}}\times(0,v_+))$}}
\end{picture}%
}
\endgroup

\caption{The unnormalized return passage. (a) The orbit reaches $\tau(D_{\mathrm{in}})$ after $N$ local iterates and $L$ transition iterates; the dashed arrow represents the translation $\tau^{-1}(x,y)=(x-1,y)$. The dotted quadrilaterals indicate the local boxes. (b) The same passage in logarithmic coordinates. Here $z_j=\widetilde f^{\,j}(z_0)$ and $(\theta_1,v_1)=\widehat{\mathcal R}_{\mathrm{raw}}(\theta,v)$.}
\label{fig:raw-return-lift}
\end{figure}
	To remove the integer translation $N(\theta, v)$, which depends on $(\theta, v)$, and the entry time $L$, use the smooth coordinate change $\mathcal{C}_\nu(\theta, v) = (\theta + \nu^{-1}\log v, v)$ and the extended normalized transition map $\mathcal{P} \colon \mathbb{R} \times [0, v_*) \to \mathbb{R} \times [0, v_+)$ to define the Poincar\'e return lift
	\begin{equation}\label{eq:return-lift-def}
		\mathcal R:=\mathcal P\circ\mathcal C_\nu:
		\mathbb R\times(0,v_*)\longrightarrow\mathbb R\times(0,v_+).
	\end{equation}
	Since $\mathcal{P}$ and $\mathcal{C}_\nu$ are smooth embeddings on their domains and satisfy $\mathcal{P} \circ T_1 = T_1 \circ \mathcal{P}$ and $\mathcal{C}_\nu \circ T_1 = T_1 \circ \mathcal{C}_\nu$, their composition $\mathcal{R}$ is a smooth embedding on $\mathbb{R} \times (0, v_*)$ satisfying $\mathcal{R} \circ T_1 = T_1 \circ \mathcal{R}$. For every $(\theta, v) \in B_{\mathrm{in}} \times (0, v_*)$, the identities $\widehat{\mathcal{P}} = T_L \circ \mathcal{P}$, $\widehat{\mathcal{C}}_\nu(\theta, v) = T_{N(\theta, v)}(\mathcal{C}_\nu(\theta, v))$, and $\mathcal P\circ T_1=T_1\circ\mathcal P$ give
	\begin{equation*}
		\begin{aligned}
			\widehat{\mathcal{R}}_{\mathrm{raw}}(\theta, v) &= (\widehat{\mathcal{P}} \circ \widehat{\mathcal{C}}_\nu)(\theta, v) \\
			&= (T_L \circ \mathcal{P})\bigl(T_{N(\theta, v)}(\mathcal{C}_\nu(\theta, v))\bigr) \\
			&= T_{L + N(\theta, v)}\bigl(\mathcal{R}(\theta, v)\bigr).
		\end{aligned}
	\end{equation*}
	Thus, on the local orbit quotient $\mathbb{T} \times (0, v_*)$ by integer translations $T_1$, the maps $\mathcal{R}$ and $\widehat{\mathcal{R}}_{\mathrm{raw}}$ induce the same quotient map on the corresponding domains, while the lift $\mathcal{R}$ is smooth on the universal cover $\mathbb{R} \times (0, v_*)$. In Section~\ref{sec:kam}, we apply Herman's $C^3$ invariant curve theorem to $\mathcal{R}$ on a subannulus $\mathbb{R} \times \eta_k\mathcal{I}_0 \subset \mathbb{R} \times (0, v_*)$. That theorem requires only that $\mathcal{R}$ be an exact symplectic embedding from its domain into $\mathbb{R} \times (0, v_+)$, not that $\mathcal{R}$ map $\mathbb{R} \times (0, v_*)$ into itself.
	
	\begin{Lemma}\label{lem:return_lift_exact}
		The Poincar\'e return lift $\mathcal{R} \colon \mathbb{R} \times (0, v_*) \to \mathbb{R} \times (0, v_+)$ is a smooth exact symplectic embedding satisfying
		\begin{equation}\label{eq:return-lift-exact-identity}
			\mathcal{R}^*(v\,\mathrm{d}\theta) - v\,\mathrm{d}\theta = \mathrm{d}S_{\mathcal{R}}, \qquad S_{\mathcal{R}} := S_{\mathcal{P}} \circ \mathcal{C}_\nu + \frac{v}{\nu},
		\end{equation}
		where $S_{\mathcal{P}}$ is a smooth primitive on $\mathbb{R}\times[0,v_*)$, $1$-periodic in $\theta$, satisfying $\mathcal{P}^*(v\,\mathrm{d}\theta) - v\,\mathrm{d}\theta = \mathrm{d}S_{\mathcal{P}}$.
	\end{Lemma}
	
	\begin{proof}
		Set $\eta_{\mathcal{P}} := \mathcal{P}^*(v\,\mathrm{d}\theta) - v\,\mathrm{d}\theta = A_{\mathcal{P}}(\theta, v)\,\mathrm{d}\theta + B_{\mathcal{P}}(\theta, v)\,\mathrm{d}v$. Since $\mathcal{P}$ preserves $\mathrm{d}\theta \wedge \mathrm{d}v$, the identity $\mathrm{d}\eta_{\mathcal{P}} = 0$ gives the compatibility condition $\partial_v A_{\mathcal{P}} = \partial_\theta B_{\mathcal{P}}$ for partial derivatives. Boundary invariance $\mathcal{P}_2(\theta, 0) = 0$ implies $A_{\mathcal{P}}(\theta, 0) = 0$. Define the $1$-periodic smooth function on $\mathbb{R} \times [0, v_*)$ by
		\begin{equation*}
			S_{\mathcal{P}}(\theta, v) := \int_0^v B_{\mathcal{P}}(\theta, s)\,\mathrm{d}s,
		\end{equation*}
		Its partial derivatives satisfy $\partial_v S_{\mathcal{P}} = B_{\mathcal{P}}$ and, by the fundamental theorem of calculus,
		\begin{equation*}
			\begin{aligned}
				\partial_\theta S_{\mathcal{P}}(\theta, v)
				&= \int_0^v \partial_\theta B_{\mathcal{P}}(\theta, s)\,\mathrm{d}s
				= \int_0^v \partial_s A_{\mathcal{P}}(\theta, s)\,\mathrm{d}s\\
				&= A_{\mathcal{P}}(\theta, v) - A_{\mathcal{P}}(\theta, 0)
				= A_{\mathcal{P}}(\theta, v),
			\end{aligned}
		\end{equation*}
		that is, $\mathrm{d}S_{\mathcal{P}} = \eta_{\mathcal{P}}$. For the local coordinate change $\mathcal{C}_\nu(\theta, v) = (\theta + \nu^{-1}\log v, v)$, direct differentiation gives
		\begin{equation*}
			\mathcal{C}_\nu^*(v\,\mathrm{d}\theta) - v\,\mathrm{d}\theta = v\,\mathrm{d}\left(\frac{\log v}{\nu}\right) = \frac{\mathrm{d}v}{\nu}.
		\end{equation*}
		Combining $\mathcal{R} = \mathcal{P} \circ \mathcal{C}_\nu$ with $(\mathcal P\circ\mathcal C_\nu)^*=\mathcal C_\nu^*\mathcal P^*$ yields
		\begin{align*}
			\mathcal{R}^*(v\,\mathrm{d}\theta) - v\,\mathrm{d}\theta
			&= \mathcal{C}_\nu^*\big(\mathcal{P}^*(v\,\mathrm{d}\theta) - v\,\mathrm{d}\theta\big) + \big(\mathcal{C}_\nu^*(v\,\mathrm{d}\theta) - v\,\mathrm{d}\theta\big) \\
			&= \mathrm{d}(S_{\mathcal{P}} \circ \mathcal{C}_\nu) + \frac{\mathrm{d}v}{\nu} \\
			&= \mathrm{d}\left( S_{\mathcal{P}} \circ \mathcal{C}_\nu + \frac{v}{\nu} \right) = \mathrm{d}S_{\mathcal{R}},
		\end{align*}
		Thus $\mathcal{R}$ is exact symplectic.
	\end{proof}
	
	In Section~\ref{sec:reconstruction} we construct a map $\iota$ and prove $\iota\circ\mathcal R=\tau^{-1}\circ\iota$.
	
	\section{Normalization of the return map and invariant curves}\label{sec:kam}
	
	In Section~\ref{sec:transition} we defined the Poincar\'e return lift on $\mathbb{R}\times(0,v_*)$ by
	\[
	\mathcal{R}(\theta, v) = (\mathcal{P} \circ \mathcal{C}_\nu)(\theta, v) = \mathcal{P}\!\left(\theta + \frac{\log v}{\nu},\; v\right).
	\]
	On a sequence of subannuli approaching the boundary $v = 0$, we use a contraction of the vertical coordinate to transform $\mathcal{R}$ into a nearly integrable exact symplectic map. We then apply Herman's $C^3$ invariant curve theorem to construct smooth invariant curves.
	
	\subsection{Normalization and remainder estimates}\label{subsec:auto-4-1}
	
	We expand the normalized transition map $\mathcal{P}$ in the vertical coordinate $v$ at $v = 0$, and then use a vertical contraction to separate the boundary terms from the higher-order remainder.
	
	Write $\mathcal{P} = (\mathcal{P}_1, \mathcal{P}_2)$. By Section~\ref{subsec:auto-3-2}, on the boundary $v = 0$,
	\[
	\mathcal{P}_1(\theta, 0) = b(\theta), \qquad \mathcal{P}_2(\theta, 0) = 0,
	\]
	Hence $\partial_\theta \mathcal{P}_1(\theta, 0) = b'(\theta) > 0$ and $\partial_\theta \mathcal{P}_2(\theta, 0) = 0$. Area preservation $\det D\mathcal{P}(\theta, v) = 1$ gives, at $v = 0$,
	\[
	b'(\theta)\,\partial_v \mathcal{P}_2(\theta, 0) = 1, \qquad
	\partial_v \mathcal{P}_2(\theta, 0) = \frac{1}{b'(\theta)}.
	\]
	Taylor expansion with integral remainder in $v$ on $\mathbb{R} \times [0, v_*)$ yields
	\begin{equation}\label{eq:P-taylor}
		\mathcal{P}_1(\theta, v) = b(\theta) + v A_{\mathrm{tr}}(\theta, v), \qquad
		\mathcal{P}_2(\theta, v) = \frac{v}{b'(\theta)} + v^2 B_{\mathrm{tr}}(\theta, v),
	\end{equation}
	where the remainder coefficients are
	\[
	A_{\mathrm{tr}}(\theta, v) := \int_0^1 \partial_v \mathcal{P}_1(\theta, tv)\,\mathrm{d}t, \qquad
	B_{\mathrm{tr}}(\theta, v) := \int_0^1 (1-t)\,\partial_v^2 \mathcal{P}_2(\theta, tv)\,\mathrm{d}t.
	\]
	By $\mathcal{P}_1(\theta + 1, v) = \mathcal{P}_1(\theta, v) + 1$ and $\mathcal{P}_2(\theta + 1, v) = \mathcal{P}_2(\theta, v)$, the functions $A_{\mathrm{tr}}$ and $B_{\mathrm{tr}}$ are smooth on $\mathbb{R} \times [0, v_*)$ and $1$-periodic in $\theta$. Differentiating under the integral sign gives
	\[
	\begin{aligned}
		\partial_\theta^i \partial_v^j A_{\mathrm{tr}}(\theta, v) &= \int_0^1 t^j \partial_\theta^i \partial_v^{j+1} \mathcal{P}_1(\theta, tv)\,\mathrm{d}t,\\
		\partial_\theta^i \partial_v^j B_{\mathrm{tr}}(\theta, v) &= \int_0^1 (1-t) t^j \partial_\theta^i \partial_v^{j+2} \mathcal{P}_2(\theta, tv)\,\mathrm{d}t.
	\end{aligned}
	\]
	Since $A_{\mathrm{tr}}$ and $B_{\mathrm{tr}}$ are $1$-periodic in $\theta$, for each fixed map $f$, on the compact subset $[0, 1] \times [0, v_*/2] \subset \mathbb{R} \times [0, v_*)$ we have, for $i + j \le 3$,
	\begin{equation}\label{eq:Atr-Btr-bounds}
		\|\partial_\theta^i \partial_v^j A_{\mathrm{tr}}\|_\infty \le \frac{\|\mathcal{P}_1 - \theta\|_{C^4} + 1}{j+1}, \qquad
		\|\partial_\theta^i \partial_v^j B_{\mathrm{tr}}\|_\infty \le \frac{\|\mathcal{P}_2\|_{C^5}}{(j+1)(j+2)}.
	\end{equation}
	The norms on the right are taken on the compact set $[0, 1] \times [0, v_*/2]$ and are bounded by a finite constant depending only on the fixed map $f$.
	
	Substituting $b(\theta) = \theta + \sigma + \zeta(\theta)$ into the Taylor expansion of $\mathcal{P}$ shows that the coefficient $1/b'(\theta) = (1 + \zeta'(\theta))^{-1}$ of the first-order term in the second component depends only on the phase variable $\theta$ and remains unchanged as $v \to 0^+$. Its deviation from $1$ is
	\begin{equation*}
		\frac{1}{b'(\theta)} - 1 = -\frac{\zeta'(\theta)}{1 + \zeta'(\theta)},
	\end{equation*}
	The absolute value and derivatives of this term are controlled by derivatives of $\zeta$ and do not decay with $v \to 0^+$.
	Substituting $\mathcal{C}_\nu(\theta, v) = (\theta + \nu^{-1}\log v, v)$ into~\eqref{eq:P-taylor} gives the expression for $\mathcal{R} = \mathcal{P} \circ \mathcal{C}_\nu$ on $\mathbb{R} \times (0, v_*)$:
	\begin{equation}\label{eq:R-explicit}
		\mathcal{R}(\theta, v) =
		\begin{pmatrix}
			b\!\left(\theta + \dfrac{\log v}{\nu}\right) + v A_{\mathrm{tr}}\!\left(\theta + \dfrac{\log v}{\nu},\; v\right) \\[8pt]
			\dfrac{v}{b'\!\left(\theta + \dfrac{\log v}{\nu}\right)} + v^2 B_{\mathrm{tr}}\!\left(\theta + \dfrac{\log v}{\nu},\; v\right)
		\end{pmatrix}.
	\end{equation}
	As $v \to 0^+$, we have $\nu^{-1}\log v \to -\infty$ and the partial derivative $\partial_v(\nu^{-1}\log v) = (\nu v)^{-1}$ is unbounded. To analyze the map for small $v$ in a bounded coordinate interval, fix
	\[
	\mathcal{I} := [1/2, 1] \Subset \mathcal{I}_0 := (1/4, 2).
	\]
	All norms in the $(\theta,s)$ coordinates below are taken on the compact annulus $\mathbb T\times\mathcal I$. For each positive integer $k \in \mathbb{N}$, set
	\begin{equation}\label{eq:layer-scaling}
		\eta_k := e^{-k\nu}, \qquad v = \eta_k s,
	\end{equation}
	and define $D_k(\theta,s):=(\theta,\eta_ks)$. For $k$ large enough that $2\eta_k = 2e^{-k\nu} \le v_*/2$, we have $D_k(\mathbb R\times\mathcal I_0)\subset\mathbb R\times(0,v_*/2)$.
	
	On $\mathbb{R} \times \mathcal{I}_0$, define the conjugated map $\mathcal{R}_k$ and the normalized return map $\mathcal{F}_k$ obtained by composing with the integer translation $T_k(\theta, s) = (\theta + k, s)$:
	\begin{equation}\label{eq:Rk-Fk-def}
		\mathcal{R}_k := D_k^{-1} \circ \mathcal{R} \circ D_k, \qquad \mathcal{F}_k := T_k \circ \mathcal{R}_k \colon \mathbb{R} \times \mathcal{I}_0 \longrightarrow \mathbb{R} \times (0, \eta_k^{-1}v_+).
	\end{equation}
	By $D_k^*(v\,\mathrm{d}\theta) = \eta_k s\,\mathrm{d}\theta$ and~\eqref{eq:return-lift-exact-identity} in Lemma~\ref{lem:return_lift_exact},
	\begin{equation}\label{eq:Fk-exactness}
		\mathcal{R}_k^*(s\,\mathrm{d}\theta) - s\,\mathrm{d}\theta
		= \eta_k^{-1} D_k^*\!\left(\mathcal{R}^*(v\,\mathrm{d}\theta) - v\,\mathrm{d}\theta\right)
		= \mathrm{d}\!\left(\eta_k^{-1} S_{\mathcal{R}} \circ D_k\right).
	\end{equation}
	Since $T_k^*(s\,\mathrm{d}\theta) = s\,\mathrm{d}\theta$, the right-hand side also equals $\mathcal{F}_k^*(s\,\mathrm{d}\theta) - s\,\mathrm{d}\theta$. Hence both $\mathcal{R}_k$ and $\mathcal{F}_k$ are exact symplectic embeddings preserving $\mathrm{d}\theta \wedge \mathrm{d}s$.
	
	Write
	\begin{equation}\label{eq:z-Lambda-def}
		z := \theta + \frac{\log s}{\nu}, \qquad \Lambda_\nu(s) := \sigma + \frac{\log s}{\nu}.
	\end{equation}
	Substitute $v = \eta_k s = e^{-k\nu}s$ into~\eqref{eq:R-explicit} and use
	\[
	\theta + \frac{\log(\eta_k s)}{\nu} = z - k, \qquad
	b(z - k) = b(z) - k = z + \sigma + \zeta(z) - k,
	\]
	together with the $1$-periodicity of $b'$, $A_{\mathrm{tr}}$, and $B_{\mathrm{tr}}$ in their first argument:
	\[
	\begin{gathered}
		b'(z - k) = b'(z),\\
		A_{\mathrm{tr}}(z - k, \eta_k s) = A_{\mathrm{tr}}(z, \eta_k s),\qquad
		B_{\mathrm{tr}}(z - k, \eta_k s) = B_{\mathrm{tr}}(z, \eta_k s),
	\end{gathered}
	\]
	In the first component, left composition with $T_k$ cancels the integer $-k$ arising from $\eta_k = e^{-k\nu}$; in the second, left composition with $D_k^{-1}$ cancels the common factor $\eta_k$. Thus
	\begin{equation}\label{eq:Fk-explicit}
		\mathcal{F}_k(\theta, s) =
		\begin{pmatrix}
			\theta + \Lambda_\nu(s) + \zeta(z) + \eta_k s A_{\mathrm{tr}}(z, \eta_k s) \\[6pt]
			\dfrac{s}{b'(z)} + \eta_k s^2 B_{\mathrm{tr}}(z, \eta_k s)
		\end{pmatrix}.
	\end{equation}
	
	Denote the integrable twist map by $\mathcal{T}_{\Lambda_\nu}(\theta, s) := (\theta + \Lambda_\nu(s), s)$. On $\mathcal{I} = [1/2, 1]$ its twist derivative satisfies
	\begin{equation}\label{eq:Lambda-twist}
		\Lambda_\nu'(s) = \frac{1}{\nu s} > 0, \qquad
		\nu^{-1} \le \Lambda_\nu'(s) \le 2\nu^{-1},
	\end{equation}
	These bounds are independent of $k$.
	
	To estimate the deviation of $\mathcal{F}_k$ from the integrable map $\mathcal{T}_{\Lambda_\nu}$, first consider the coefficient $h := (b')^{-1} = (1 + \zeta')^{-1}$ in the second component of $\mathcal{F}_k$. Under the assumptions of Theorem~\ref{thm:boundary-transition-criterion}, we may choose the threshold so that $\delta_{\nu, \gamma_\xi} \le 1/2$. If $\|\zeta\|_{C^4(\mathbb{T})} \le 1/2$, then $\|\zeta'\|_\infty \le 1/2$, and the denominator has the uniform positive lower bound $b' = 1 + \zeta' \ge 1/2$. Hence $\|h\|_\infty \le 2$ and $\|h - 1\|_\infty \le 2\|\zeta'\|_\infty$. Differentiating $(1 + \zeta')h = 1$ a total of $n$ times, with $1 \le n \le 3$, gives
	\[
	h^{(n)} = -\frac{1}{1 + \zeta'} \sum_{j=1}^n \binom{n}{j} \zeta^{(j+1)} h^{(n-j)}.
	\]
	Using the lower bound $1 + \zeta' \ge 1/2$ and $\|\zeta\|_{C^4(\mathbb{T})} \le 1/2$, induction on $n$ gives a constant $C > 0$ such that
	\begin{equation}\label{eq:h-minus-1-C3}
		\|h - 1\|_{C^3(\mathbb{T})} \le C\|\zeta\|_{C^4(\mathbb{T})}.
	\end{equation}
	For $s\in \mathcal{I}$ and $1 \le j \le 3$,
	\[
	\partial_s^j z = \frac{(-1)^{j-1}(j-1)!}{\nu s^j}, \qquad
	\partial_\theta z = 1, \qquad
	\partial_s(\eta_k s) = \eta_k \le 1.
	\]
	Repeatedly applying the product and chain rules (Appendix~\ref{subsec:chain-rule}) to the terms in~\eqref{eq:Fk-explicit}, and using~\eqref{eq:Atr-Btr-bounds} and~\eqref{eq:h-minus-1-C3}, yields
	\[
	\begin{aligned}
		\|\zeta(z)\|_{C^3} + \|s(h(z) - 1)\|_{C^3} &\le C_\nu \|\zeta\|_{C^4(\mathbb{T})},\\
		\|\eta_k s A_{\mathrm{tr}}(z, \eta_k s)\|_{C^3} + \|\eta_k s^2 B_{\mathrm{tr}}(z, \eta_k s)\|_{C^3} &\le C_\nu \eta_k \bigl(\|A_{\mathrm{tr}}\|_{C^3} + \|B_{\mathrm{tr}}\|_{C^3}\bigr).
	\end{aligned}
	\]
	The norms of $A_{\mathrm{tr}}$ and $B_{\mathrm{tr}}$ on the right are taken on the corresponding compact sets. Every term involving $A_{\mathrm{tr}}$ or $B_{\mathrm{tr}}$ has a factor $\eta_k = e^{-k\nu}$, whereas the first estimate gives an $\zeta$ term independent of $k$. We obtain the $C^3$ estimate
	\begin{equation}\label{eq:Fk-C3-bound}
		\|\mathcal{F}_k - \mathcal{T}_{\Lambda_\nu}\|_{C^3} \le C_\nu \|\zeta\|_{C^4(\mathbb{T})} + C_{f,\nu} e^{-k\nu},
	\end{equation}
	Here $C_{f,\nu}$ is independent of $k$, and $C_\nu$ can be chosen uniformly when $\nu$ ranges over a compact subset of $(0, \infty)$.

{
The same estimate holds on any fixed compact interval $J\Subset(0,\infty)$ in place of $\mathcal I$, once $k$ is large enough that $e^{-k\nu}\sup J\le v_*/2$. Indeed, for $s\in J$,
\[
\left|\partial_s^j\left(\theta+\nu^{-1}\log s\right)\right|
\le\frac{(j-1)!}{\nu(\inf J)^j}\qquad(1\le j\le3),
\]
and every transverse remainder still carries the factor $e^{-k\nu}$. The constants may depend on $J$ and are uniform for $\nu$ in compact subsets of $(0,\infty)$. We use this observation on $J=[1/4,2]$ in Section~\ref{subsec:island-proof}.
}
	
	\subsection{A \texorpdfstring{$C^3$}{C3} invariant curve theorem}\label{subsec:auto-4-2}
	
	We use Herman's $C^3$ invariant curve theorem in \cite[Chapter V, \S 7.10, p.~78]{Her86}, which does not require the embedding to be defined on the whole cylinder.
	
	\begin{Lemma}[Herman's $C^3$ theorem]\label{lem:herman-C3}
		Fix an irrational number $\alpha$ of constant type, the arithmetic lower bound $\gamma > 0$ in~\eqref{eq:constant-type}, a constant $\delta_{\mathrm{H}} > 0$, and $p_S \in (3/2, 2)$. Let
		\[
		H_0(\theta, r) = (\theta + \alpha + t(r),\; r), \qquad t(0) = 0,
		\]
		where $t$ is smooth and satisfies $B^{-1} \le t'(r) \le B$ for some $B > 1$. For every $\epsilon_{\mathrm{g}} > 0$ there exists $\epsilon_{\mathrm{H}} > 0$, depending only on these unperturbed data and the arithmetic constant, with the following property. If a map $H$ is defined on a neighborhood of $\mathbb R\times[-2\delta_{\mathrm H},2\delta_{\mathrm H}]$, commutes with $T_1$, induces an exact symplectic $C^3$ embedding, and satisfies
		\[
		\|H - H_0\|_{C^3} < \epsilon_{\mathrm{H}},
		\]
		then it has an invariant graph $r = \psi_{\mathrm{H}}(\theta)$ of rotation number $\alpha$, with $\|\psi_{\mathrm{H}}\|_{C^0} < \delta_{\mathrm{H}}$ and
		\[
		\psi_{\mathrm{H}} \in W^{2, p_S}(\mathbb{T}), \qquad \|\psi_{\mathrm{H}}\|_{W^{2, p_S}} < \epsilon_{\mathrm{g}}.
		\]
		If $H$ is smooth, then the graph function $\psi_{\mathrm{H}}$ is smooth.
	\end{Lemma}
	

	The space $W^{2, p_S}(\mathbb{T})$ consists of periodic functions whose weak derivatives through order two belong to $L^{p_S}(\mathbb{T})$. In one dimension, the fundamental theorem of calculus and H\"older's inequality control the oscillations of $\psi_{\mathrm{H}}$ and $\psi_{\mathrm{H}}'$ by the $L^{p_S}$ norms of their next derivatives, while their means are controlled by their own $L^{p_S}$ norms. Thus there exists $C_{p_S} > 0$ such that
	\begin{equation}\label{eq:sobolev-C1}
		\|\psi_{\mathrm{H}}\|_{C^1(\mathbb{T})} \le C_{p_S}\|\psi_{\mathrm{H}}\|_{W^{2, p_S}(\mathbb{T})}.
	\end{equation}
	Choosing $\epsilon_{\mathrm{g}}$ appropriately therefore makes the invariant graph arbitrarily small in the $C^1(\mathbb{T})$ norm.
	
	\paragraph{The local model and exactness.}
	Let $a \in \operatorname{int}\mathcal{I}$ and $\alpha = \Lambda_\nu(a)$. Introduce the vertical translation $S_a(\theta, r) := (\theta, a + r)$ to center the vertical coordinate at $s = a$. The centered integrable map is
	\begin{equation}\label{eq:centered-H0}
		\begin{gathered}
			H_0(\theta, r) := (S_a^{-1} \circ \mathcal{T}_{\Lambda_\nu} \circ S_a)(\theta, r)
			= \bigl(\theta + \alpha + t_{\nu, a}(r),\; r\bigr),\\
			t_{\nu, a}(r) := \frac{1}{\nu}\log\!\left(1 + \frac{r}{a}\right).
		\end{gathered}
	\end{equation}
	Direct calculation gives $t_{\nu, a}(0) = 0$ and, on a local interval containing $[-2\delta_{\mathrm{H}}, 2\delta_{\mathrm{H}}]$, where $2\delta_{\mathrm{H}} < \operatorname{dist}(a, \partial\mathcal{I})$,
	\[
	t_{\nu, a}'(r) = \frac{1}{\nu(a + r)} > 0.
	\]
	Using a smooth cutoff as in Appendix~\ref{subsec:cutoffs}, extend $t_{\nu, a}'$ outside this interval to a smooth function on $\mathbb{R}$ with positive upper and lower bounds. Integration with the initial value $0$ gives a globally smooth function $\widetilde{t}_{\nu, a}$, without changing the map on the local neighborhood $\mathbb{T} \times [-2\delta_{\mathrm{H}}, 2\delta_{\mathrm{H}}]$.
	
	For the centered map $H := S_a^{-1} \circ \mathcal{F}_k \circ S_a$, write $H = (H_1, H_2)$. Using $\mathcal{F}_k^*(s\,\mathrm{d}\theta) - s\,\mathrm{d}\theta = \mathrm{d}S_k$ from~\eqref{eq:Fk-exactness}, where $S_k := \eta_k^{-1} S_{\mathcal{R}} \circ D_k$, and $S_a^*(s\,\mathrm{d}\theta) = (a + r)\,\mathrm{d}\theta$, direct calculation of the pullback gives
	\begin{equation}\label{eq:H-exactness}
		\begin{aligned}
			H^*(r\,\mathrm{d}\theta) - r\,\mathrm{d}\theta
			&= H^*\bigl(S_a^*(s\,\mathrm{d}\theta) - a\,\mathrm{d}\theta\bigr) - r\,\mathrm{d}\theta \\
			&= S_a^*\bigl(\mathcal{F}_k^*(s\,\mathrm{d}\theta)\bigr) - a\,\mathrm{d}H_1 - r\,\mathrm{d}\theta \\
			&= (a + r)\,\mathrm{d}\theta + \mathrm{d}(S_k \circ S_a) - a\,\mathrm{d}H_1 - r\,\mathrm{d}\theta \\
			&= \mathrm{d}\bigl(S_k \circ S_a - a(H_1 - \theta)\bigr).
		\end{aligned}
	\end{equation}
	Since $H_1(\theta + 1, r) = H_1(\theta, r) + 1$, the function $H_1(\theta, r) - \theta$ is $1$-periodic in $\theta$. The vertical translation thus preserves exactness; its primitive changes only by a periodic term induced by the closed form $a\,\mathrm{d}\theta$. The construction of the normalized transition map in Section~\ref{sec:transition} also ensures that, for sufficiently large $k$, the map $\mathcal{F}_k$, and hence $H$, is a smooth embedding on the required local neighborhood. By the choice of $\delta_{\mathrm H}$, this neighborhood can be taken inside $S_a^{-1}(\mathbb R\times\mathcal I_0)$. Lemma~\ref{lem:herman-C3} therefore applies directly to $H$.
	
	\paragraph{Uniformity in the parameters.}
	Below, $\nu$ ranges over a fixed compact subset of $(0, \infty)$, and $a$ over a fixed compact subinterval of $\operatorname{int}\mathcal{I}$. Choose $\delta_{\mathrm H}>0$ independently of $\nu,a$ so that $2\delta_{\mathrm H}<\operatorname{dist}(a,\partial\mathcal I)$. On a slightly larger vertical interval, the derivative formulas
	\[
	t_{\nu, a}^{(j)}(r) = \frac{(-1)^{j-1}(j-1)!}{\nu(a + r)^j} \qquad (j \ge 1)
	\]
	give uniform upper bounds for all derivatives of $t_{\nu,a}$ and a uniform positive lower bound for $t_{\nu,a}'$. Together with the reduction in \cite[Chapter V, \S 4.1 and \S\S 7.2--7.5]{Her86}, this allows $\epsilon_{\mathrm H}$ in Lemma~\ref{lem:herman-C3} to be chosen uniformly. More explicitly, write the centered map as
	\[
	H(\theta, r) = \bigl(\theta + \alpha + t(r) + \Delta_1(\theta, r),\; r + \Delta_2(\theta, r)\bigr),
	\]
	where $t = t_{\nu, a}$, and set $z_a(\theta, r) := \theta + \nu^{-1}\log(a + r)$. The component perturbations are then
	\[
	\begin{aligned}
		\Delta_1(\theta, r) &= \zeta\bigl(z_a(\theta, r)\bigr) + \eta_k(a + r)A_{\mathrm{tr}}\bigl(z_a(\theta, r),\, \eta_k(a + r)\bigr),\\
		\Delta_2(\theta, r) &= (a + r)\bigl(h(z_a(\theta, r)) - 1\bigr) + \eta_k(a + r)^2 B_{\mathrm{tr}}\bigl(z_a(\theta, r),\, \eta_k(a + r)\bigr).
	\end{aligned}
	\]
	Introduce the coordinate change $G_H(\theta, r) := (\theta,\, t(r) + \Delta_1(\theta, r))$. For $\|\Delta_1\|_{C^1}$ sufficiently small, $\partial_r(t + \Delta_1)$ has a fixed positive lower bound. The implicit function theorem gives a unique smooth inverse $r = \chi(\theta, R)$ determined by $t(\chi) + \Delta_1(\theta, \chi) = R$. Differentiating this equation three times shows that $\chi$ has a uniform $C^3$ bound. Direct substitution in the coordinates $(\theta, R)$ gives the conjugated map
	\[
	(G_H \circ H \circ G_H^{-1})(\theta, R) = \bigl(\theta + \alpha + R,\; R + \phi(\theta, R)\bigr),
	\]
	{
We perform this calculation on a fixed smaller strip whose closure lies in the domain of $H$. Taking $\|\Delta_2\|_{C^0}$ smaller than the vertical margin ensures that both $r$ and $r+\Delta_2(\theta,r)$ remain in the common domain. The inverse coordinate change is defined on a fixed strip $|R|\le\delta_0$, uniformly in the parameters.
}

	Setting $r = \chi(\theta, R)$ and $\Delta_2 = \Delta_2(\theta, r)$, the remainder $\phi$ is
	\[
	\phi(\theta, R) = t(r + \Delta_2) - t(r) + \Delta_1(\theta + \alpha + R,\, r + \Delta_2) - \Delta_1(\theta, r).
	\]
	Using $t(r + \Delta_2) - t(r) = \Delta_2\int_0^1 t'(r + s \Delta_2)\,\mathrm{d}s$ and the chain rule gives
	\begin{equation}\label{eq:phi-C3-bound}
		\|\phi\|_{C^3} \le C\bigl(\|\Delta_1\|_{C^3} + \|\Delta_2\|_{C^3}\bigr),
	\end{equation}
	Here $C > 0$ depends only on the common domain, an upper bound for $(t')^{-1}$, and $\|t\|_{C^4}$, and is therefore uniformly bounded for this family of logarithmic functions. Although $G_H$ need not be symplectic, exactness of $H$ implies the intersection property, which is preserved by conjugacy. The reduced map thus belongs to the class required by Herman's reduction theorem \cite[Chapter V, \S\S 4.7, 5.5]{Her86}. In the reduction theorem of \cite[Chapter V, \S 7.2, p.~74]{Her86}, the perturbation threshold on a common narrower strip $\{|R| \le \delta_0\}$ is $c_{p_S}\gamma\min\{\gamma, \delta_0\}$, and the extension steps in \cite[Chapter V, \S 4.9 and \S 7.5]{Her86} preserve continuous $C^3$ dependence on the local data. For fixed $\gamma>0$ and $\delta_0>0$, estimate~\eqref{eq:phi-C3-bound} thus gives a common positive threshold $\epsilon_{\mathrm{H}} > 0$, which may also incorporate any prescribed upper bound on the invariant graph norm. Adding an integer horizontal translation does not affect this conclusion.
	
	\begin{Remark}\label{rem:derivative-counts}
		We use a $C^3$ invariant curve theorem; Herman also allows $C^{2+\mathrm{Lip}}$ perturbations in \cite[Chapter V, \S 7.11]{Her86}. In~\eqref{eq:Fk-explicit}, taking three derivatives of $1/b'$ requires four derivatives of $\zeta$, which accounts for the $C^4$ norm condition on $\zeta$. The integral formulas for $A_{\mathrm{tr}}$ and $B_{\mathrm{tr}}$ involve, respectively, the first and second vertical partial derivatives of $\mathcal{P}_1$ and $\mathcal{P}_2$. A finite $C^5$ norm of the normalized transition map therefore controls the $C^3$ remainder in~\eqref{eq:Fk-C3-bound}. Since $k$ is chosen after fixing the map, this $C^5$ norm need not be uniformly bounded in an external parameter.
	\end{Remark}
	
	\subsection{Frequency selection and parameterization}\label{subsec:auto-4-3}
	
	Fix the constant-type data $(\xi, \gamma_\xi)$ in~\eqref{eq:fix-constant-type}. All arithmetic estimates depend only on the lower bound $\gamma_\xi$.
	
	On $\mathcal{I} = [1/2, 1]$, the image interval $\Lambda_\nu(\mathcal{I})$ has length $(\log 2)/\nu$. Choose a positive integer $Q_* > 2\nu/(\log 2)$. The mesh size $1/Q_*$ of $\{(j + \xi)/Q_* : j \in \mathbb{Z}\}$ is smaller than the length of the middle half of $\Lambda_\nu(\mathcal{I})$, so this mesh contains a point $\alpha = (j + \xi)/Q_*$ in that middle half. For every $p \in \mathbb{Z}$ and $q \in \mathbb{N}$,
	\begin{equation}\label{eq:alpha-diophantine}
		\left|\alpha - \frac{p}{q}\right|
		= \frac{1}{Q_*}\left|\xi - \frac{Q_* p - jq}{q}\right|
		\ge \frac{\gamma_*}{q^2}, \qquad
		\gamma_* := \frac{\gamma_\xi}{Q_*}.
	\end{equation}
	Moreover, the distance from $\alpha$ to each endpoint of $\Lambda_\nu(\mathcal{I})$ is at least $(\log 2)/(4\nu)$. Set
	\begin{equation}\label{eq:s-alpha-def}
		s_\alpha := \exp\bigl(\nu(\alpha - \sigma)\bigr), \qquad
		2^{-3/4} \le s_\alpha \le 2^{-1/4}.
	\end{equation}
	The distance from $s_\alpha$ to each endpoint of $\mathcal{I} = [1/2, 1]$ is therefore bounded below by the positive number
	\[
	d_s := \min\bigl\{2^{-3/4} - 1/2,\; 1 - 2^{-1/4}\bigr\} > 0.
	\]
	Fix a compact subinterval $[\nu_-, \nu_+]$. For $\nu \in [\nu_-, \nu_+] \Subset (0, \infty)$, the choice $Q_* > 2\nu_+/(\log 2)$ works for all $\nu$ in that interval and every mean displacement $\sigma$, and gives a uniform endpoint-distance lower bound $(\log 2)/(4\nu_+)$ for the frequency. Once $(\xi, \gamma_\xi)$ and $Q_*$ are fixed, the chosen $\alpha$ is independent of $k$, and the dependence of $\delta_{\nu,\gamma_\xi}$ on the arithmetic data is through $\gamma_*=\gamma_\xi/Q_*$.
	
	To obtain bounds for the normalized map and its derivatives that are uniform in the mean displacement $\sigma$, write $j_\sigma := \lfloor \sigma \rfloor$ and temporarily replace $\mathcal{F}_k$ and $\alpha$ by $T_{-j_\sigma} \circ \mathcal{F}_k$ and $\alpha - j_\sigma$, respectively. The integrable part $\Lambda_\nu(s) - j_\sigma$ and its derivatives are then uniformly bounded. After applying Lemma~\ref{lem:herman-C3} to parameterize the invariant curve, use $K_k(t + j_\sigma) = K_k(t) + (j_\sigma, 0)$ to restore the integer part $j_\sigma$. In every subsequent formula recovering the rotation number $\omega_k$ of the original map, $\alpha$ denotes the rotation number of the Poincar\'e return lift with its integer part retained.
	
	In centered coordinates \eqref{eq:centered-H0}, choose $a = s_\alpha$ and $\delta_{\mathrm{H}} = d_s / 8$. The arithmetic lower bound $\gamma_*$ in~\eqref{eq:alpha-diophantine}, the positive number $d_s$, and the uniformity established in Section~\ref{subsec:auto-4-2} allow us to choose a positive threshold $\varepsilon_{\mathrm{KAM}} > 0$ in Lemma~\ref{lem:herman-C3} so that the $C^1$ norm of the invariant graph in~\eqref{eq:sobolev-C1} satisfies~\eqref{eq:psi-k-C1-bound} below. Choose $\delta_{\nu, \gamma_\xi} $ sufficiently small that
	\[
	C_\nu \delta_{\nu, \gamma_\xi} < \frac{\varepsilon_{\mathrm{KAM}}}{2}.
	\]
	For a fixed map $f$ satisfying $\|\zeta\|_{C^4(\mathbb{T})} < \delta_{\nu, \gamma_\xi}$, choose $k_f \in \mathbb{N}$ sufficiently large that
	\[
	C_{f, \nu} e^{-k_f\nu} < \frac{\varepsilon_{\mathrm{KAM}}}{2}, \qquad
	k_f > \alpha + 1, \qquad
	e^{-k_f\nu} \sup\mathcal{I}_0 < \frac{v_*}{2}.
	\]
	Then for every $k \ge k_f$, estimate~\eqref{eq:Fk-C3-bound} and invariance of the norm under translation give
	\[
	\|H - H_0\|_{C^3} \le \|\mathcal{F}_k - \mathcal{T}_{\Lambda_\nu}\|_{C^3}
	\le C_\nu \|\zeta\|_{C^4(\mathbb{T})} + C_{f, \nu} e^{-k\nu}
	< \varepsilon_{\mathrm{KAM}}.
	\]
	By Lemma~\ref{lem:herman-C3} and~\eqref{eq:sobolev-C1}, for every $k \ge k_f$ the map $\mathcal{F}_k$ has a smooth invariant graph of rotation number $\alpha$,
	\[
	\kappa_k^{\mathrm{sc}}(\theta) := \bigl(\theta,\; s_\alpha + \psi_k(\theta)\bigr), \qquad
	\psi_k(\theta + 1) = \psi_k(\theta),
	\]
	satisfying
	\begin{equation}\label{eq:psi-k-C1-bound}
		\|\psi_k\|_{C^1(\mathbb{T})} < \frac{d_s}{2}.
	\end{equation}
	This estimate keeps the invariant graph in the interior of $\mathbb{T} \times \operatorname{int}\mathcal{I}$ at a fixed positive distance from its boundary.
	
	Let $a_k(\theta) := \pi_1\Bigl(\mathcal{F}_k\bigl(\kappa_k^{\mathrm{sc}}(\theta)\bigr)\Bigr)$ be the lift of the circle map induced by $\mathcal{F}_k$ on the invariant graph $\kappa_k^{\mathrm{sc}}$. Differentiating $\mathcal{F}_k \circ \kappa_k^{\mathrm{sc}} = \kappa_k^{\mathrm{sc}} \circ a_k$ gives
	\[
	D\mathcal{F}_k\bigl(\kappa_k^{\mathrm{sc}}(\theta)\bigr)
	\begin{pmatrix}
		1 \\ \psi_k'(\theta)
	\end{pmatrix}
	= a_k'(\theta)
	\begin{pmatrix}
		1 \\ \psi_k'\bigl(a_k(\theta)\bigr)
	\end{pmatrix}.
	\]
	Since $D\mathcal{F}_k$ is invertible everywhere, $a_k'(\theta)$ never vanishes. Together with $a_k(\theta + 1) = a_k(\theta) + 1$ this gives $a_k'(\theta) > 0$. Thus $a_k$ induces a smooth orientation-preserving circle diffeomorphism with constant-type irrational rotation number $\alpha$. The smooth linearization theorem for circle diffeomorphisms (see \cite[Chapter V, \S 5.1, p.~43]{Her86}) gives a smooth diffeomorphism $h_k : \mathbb{R} \to \mathbb{R}$ satisfying $h_k(t + 1) = h_k(t) + 1$ and $h_k'(t) > 0$ such that
	\[
	a_k\bigl(h_k(t)\bigr) = h_k(t + \alpha).
	\]
	Define the smooth regular degree-one parameterization
	\begin{equation}\label{eq:Kk-def}
		K_k(t) := \bigl(h_k(t),\; s_\alpha + \psi_k(h_k(t))\bigr), \qquad
		(K_k)_1'(t) = h_k'(t) > 0,
	\end{equation}
	Then $K_k$ satisfies the conjugacy equation
	\begin{equation}\label{eq:Fk-Kk-conjugacy}
		\mathcal{F}_k\bigl(K_k(t)\bigr) = K_k(t + \alpha).
	\end{equation}
	No bound on the size of $h_k - \operatorname{id}$ or uniform positive lower bound for $h_k'$ is needed here. Finally, set $\widetilde{K}_k := D_k \circ K_k$. By $\mathcal{F}_k = T_k \circ D_k^{-1} \circ \mathcal{R} \circ D_k$ and $K_k(t + 1) = K_k(t) + (1, 0)$,
	\begin{equation}\label{eq:R-tilde-Kk}
		\mathcal{R}\bigl(\widetilde{K}_k(t)\bigr) = \widetilde{K}_k(t + \varrho_k), \qquad
		\varrho_k := \alpha - k < 0, \qquad
		\widetilde{K}_k(t + 1) = \widetilde{K}_k(t) + (1, 0).
	\end{equation}
	
	\section{Proofs of the general criterion and its corollaries}\label{sec:reconstruction}
	We first define the orbit extension map, which maps an invariant curve of the return map back to the original cylinder. We then prove Theorem~\ref{thm:boundary-transition-criterion} and Corollaries~\ref{cor:f0-invariant} and~\ref{cor:graph-corollary}.
	
	\subsection{The orbit extension map}\label{subsec:auto-5-1}
	Fix $0<v_*<v_+$ from Section~\ref{sec:transition}. Since $|B_{\mathrm{in}}|>1$, for any $(\theta,v)\in\mathbb R\times(0,v_+)$ we can choose an integer $n$ such that $\theta-n\in B_{\mathrm{in}}$ and define
	\begin{equation}\label{eq:iota-def}
		\iota(\theta,v):=\widetilde f^n\bigl(J_{\mathrm{in}}(\theta-n,v)\bigr),
		\qquad \theta-n\in B_{\mathrm{in}}.
	\end{equation}
	On the open overlap corresponding to adjacent integers, the identity $J_{\mathrm{in}}(\theta + 1, v) = \widetilde{f}\bigl(J_{\mathrm{in}}(\theta, v)\bigr)$ shows that the expressions for $n$ and $n + 1$ satisfy
	\begin{equation}\label{eq:iota-overlap}
		\widetilde{f}^{n+1}\bigl(J_{\mathrm{in}}(\theta - n - 1, v)\bigr) = \widetilde{f}^n\bigl(J_{\mathrm{in}}(\theta - n, v)\bigr).
	\end{equation}
	This equality holds on an open set, so all derivatives agree as well. Each local expression is a composition of smooth diffeomorphisms. Hence $\iota : \mathbb{R} \times (0, v_+) \to \mathbb{R}^2$ is a smooth local diffeomorphism satisfying
	\begin{equation}\label{eq:time-equivariance}
		\iota(\theta + 1, v) = \widetilde{f}\bigl(\iota(\theta, v)\bigr).
	\end{equation}
	By continuity, $\iota=J_{\mathrm{in}}$ also holds on $\overline{B_{\mathrm{in}}}\times(0,v_+)$.
	
	Treat the outgoing chart $J_{\mathrm{out}}$ in the same way: for any $(\theta,v)\in\mathbb R\times(0,v_*)$, choose an integer $n$ with $\theta-n\in B_{\mathrm{out}}$ and define
	\begin{equation}\label{eq:iota-out-def}
		\iota_{\mathrm{out}}(\theta,v)
		:=\widetilde f^n\bigl(J_{\mathrm{out}}(\theta-n,v)\bigr),
		\qquad \theta-n\in B_{\mathrm{out}}.
	\end{equation}
	On the overlap of adjacent integer translates, the expressions agree and are independent of the choice of $n$. Hence $\iota_{\mathrm{out}} : \mathbb{R} \times (0, v_*) \to \mathbb{R}^2$ is a smooth local diffeomorphism and
	\[
	\iota_{\mathrm{out}}(\theta+1,v)
	=\widetilde f\bigl(\iota_{\mathrm{out}}(\theta,v)\bigr).
	\]
	By continuity, $\iota_{\mathrm{out}} = J_{\mathrm{out}}$ also holds on $\overline{B_{\mathrm{out}}} \times (0, v_*)$. We next verify that $\iota_{\mathrm{out}} \circ \mathcal{C}_\nu = \iota$ throughout $\mathbb{R} \times (0, v_*)$.

	For any $(\theta,v)\in B_{\mathrm{in}}\times(0,v_*)$, set $\vartheta:=\theta+\nu^{-1}\log v$ and choose an integer $n$ such that $\vartheta-n\in B_{\mathrm{out}}$. The condition $v_*<e^{-\nu(\theta_{\mathrm{in}}^+-\theta_-)}$ in Lemma~\ref{lem:finite_transition} gives
	\[
	\vartheta<\theta_-<\vartheta-n,
	\qquad n<0.
	\]
	For integers $0 \le j \le -n$, consider the coordinate points
	\[
	Q_j = e^{j\nu} v e^{\nu\theta}, \qquad P_j = e^{-j\nu} e^{-\nu\theta}, \qquad Q_j P_j = v.
	\]
	The initial and terminal points lie in the chosen incoming and outgoing basic strips, respectively. Since $Q_j$ increases with $j$ and $P_j$ decreases with $j$, all intermediate points remain in the closed linearizing box $\mathcal{B}$. Induction on $j$ using the local linearization formula \eqref{eq:local-linearization} in (G2) gives
	\begin{equation*}
		\widetilde{f}^j\bigl(J_{\mathrm{in}}(\theta, v)\bigr) = \Phi_o^{-1}(Q_j, P_j) \qquad (0 \le j \le -n).
	\end{equation*}
	Taking $j = -n$, the left-hand side is $\widetilde{f}^{-n}\bigl(J_{\mathrm{in}}(\theta, v)\bigr)$ and the right-hand side is $\Phi_o^{-1}(Q_{-n}, P_{-n}) = J_{\mathrm{out}}(\vartheta - n, v)$. Thus
	
	\[
	J_{\mathrm{out}}(\vartheta - n, v) = \widetilde{f}^{-n}\bigl(J_{\mathrm{in}}(\theta, v)\bigr),
	\]
	It follows that
	\[\bigl(\iota_{\mathrm{out}} \circ \mathcal{C}_\nu\bigr)(\theta, v) = \iota_{\mathrm{out}}(\vartheta, v) = \widetilde{f}^n\bigl(J_{\mathrm{out}}(\vartheta - n, v)\bigr) = J_{\mathrm{in}}(\theta, v) = \iota(\theta, v).\]
	
	Extend this identity to the full domain $\mathbb{R} \times (0, v_*)$. For any $(\theta, v) \in \mathbb{R} \times (0, v_*)$, choose an integer $m \in \mathbb{Z}$ such that $\theta - m \in B_{\mathrm{in}}$. Using the commutation relation $\mathcal{C}_\nu \circ T_m = T_m \circ \mathcal{C}_\nu$ together with $\iota_{\mathrm{out}} \circ T_m = \widetilde{f}^m \circ \iota_{\mathrm{out}}$ and $\iota \circ T_m = \widetilde{f}^m \circ \iota$ gives
	\begin{equation*}
		\begin{aligned}
			\bigl(\iota_{\mathrm{out}} \circ \mathcal{C}_\nu\bigr)(\theta, v) &= \bigl(\iota_{\mathrm{out}} \circ \mathcal{C}_\nu \circ T_m\bigr)(\theta - m, v) \\
			&= \bigl(\iota_{\mathrm{out}} \circ T_m \circ \mathcal{C}_\nu\bigr)(\theta - m, v) \\
			&= \widetilde{f}^m\bigl((\iota_{\mathrm{out}} \circ \mathcal{C}_\nu)(\theta - m, v)\bigr) \\
			&= \widetilde{f}^m\bigl(\iota(\theta - m, v)\bigr) = \iota(\theta, v).
		\end{aligned}
	\end{equation*}
	Thus $\iota_{\mathrm{out}} \circ \mathcal{C}_\nu = \iota$ holds throughout $\mathbb{R} \times (0, v_*)$.
	
	By Lemma~\ref{lem:finite_transition}, the image of the unnormalized finite-step transition map satisfies $\widehat{\mathcal{P}}(\overline{B_{\mathrm{out}}} \times (0, v_*)) \subset B_{\mathrm{in}} \times (0, v_+)$. At these points, the orbit extension map $\iota$ agrees with the incoming chart $J_{\mathrm{in}}$. The definition of $\widehat{\mathcal{P}} = J_{\mathrm{in}}^{-1} \circ \tau^{-1} \circ \widetilde{f}^L \circ J_{\mathrm{out}}$ in Section~\ref{subsec:auto-3-2} first gives, on the outgoing basic rectangle,
	\[
	\iota \circ \widehat{\mathcal{P}} = \tau^{-1}\circ\widetilde{f}^L\circ\iota_{\mathrm{out}}.
	\]
	By~\eqref{eq:time-equivariance} and $\widehat{\mathcal P}\circ T_m=T_m\circ\widehat{\mathcal P}$, this identity extends to the entire domain $\mathbb{R} \times (0, v_*)$. Using the decomposition of $\mathcal{R} = \mathcal{P} \circ \mathcal{C}_\nu = T_{-L} \circ \widehat{\mathcal{P}} \circ \mathcal{C}_\nu$ in~\eqref{eq:return-lift-def} of Section~\ref{sec:transition}, we obtain
	\[
	\begin{aligned}
		\iota \circ \mathcal{R}
		&= \iota \circ T_{-L} \circ \widehat{\mathcal{P}} \circ \mathcal{C}_\nu \\
		&= \widetilde{f}^{-L}\circ\iota \circ \widehat{\mathcal{P}} \circ \mathcal{C}_\nu \\
		&= \widetilde{f}^{-L}\circ\tau^{-1}\circ\widetilde{f}^L\circ\iota_{\mathrm{out}} \circ \mathcal{C}_\nu
		= \tau^{-1}\circ\iota.
	\end{aligned}
	\]
	Consequently,
	\begin{equation}\label{eq:spatial-equivariance}
		\iota \circ \mathcal{R} = \tau^{-1}\circ\iota.
	\end{equation}
	Under the orbit extension map $\iota$, $T_1$ corresponds to one iterate of the lift $\widetilde{f}$, whereas the Poincar\'e return lift $\mathcal{R}$ corresponds to the unit horizontal translation $\tau^{-1}(x,y)=(x-1,y)$: if $\iota(\theta, v) = (x, y)$, then $\iota(\mathcal{R}(\theta, v)) = (x - 1, y)$, and these points project to the same point of the cylinder $\mathbb{A}$.
	Equation~\eqref{eq:spatial-equivariance} is well defined on $\mathbb{R} \times (0, v_*)$. By~\eqref{eq:return-lift-def} and Lemma~\ref{lem:finite_transition}, the image of $\mathcal{R}$ lies in the domain of $\iota$:
	\begin{equation*}
		\mathcal{R}\bigl(\mathbb{R} \times (0, v_*)\bigr) \subset \mathbb{R} \times (0, v_+).
	\end{equation*}
	
	For every $k \ge k_f$, the choice of $k_f$ in Section~\ref{subsec:auto-4-3} ensures that the curve $\widetilde{K}_k(\mathbb{R})$ satisfies $0 < v < v_*/2$ and hence lies in the interior of the strip: $\widetilde{K}_k(\mathbb{R}) \subset \mathbb{R} \times (0, v_*/2)$. Moreover, the conjugacy equation~\eqref{eq:R-tilde-Kk} implies that $\widetilde{K}_k(\mathbb{R})$ is invariant under $\mathcal{R}$. Thus~\eqref{eq:spatial-equivariance} is well defined and holds on $\widetilde{K}_k(\mathbb{R})$.
	\subsection{Proof of Theorem~\ref{thm:boundary-transition-criterion}}\label{subsec:auto-5-2}
	We first give the arithmetic lemma needed to estimate the rotation number.
	
	\begin{Lemma}\label{lem:reciprocal-frequency}
		Suppose that $\alpha\in\mathbb R$ satisfies
		\[
		\left|\alpha-\frac pq\right|\ge\frac{\gamma_*}{q^2}
		\qquad(p\in\mathbb Z,\ q\in\mathbb N),
		\]
		where $\gamma_*>0$. If $k\in\mathbb Z$, $M\in\mathbb N$, and $k>\alpha$, and we set $\omega_k=1/[M(k-\alpha)]$, then
		\[
		\left|\omega_k-\frac pq\right|
		\ge\min\left\{\frac{\omega_k}{2},\frac{2\gamma_*}{3M}\right\}\frac1{q^2}.
		\]
		In particular, if $k-\alpha\ge3/(4\gamma_*)$, then
		\[
		\left|\omega_k-\frac pq\right|\ge\frac{\omega_k}{2q^2}
		\qquad(p\in\mathbb Z,\ q\in\mathbb N).
		\]
		This condition applies simultaneously to all positive integers $M$.
	\end{Lemma}
	\begin{proof}
		Set $a := k - \alpha > 0$. By assumption, the irrational number $a$ inherits the constant-type lower bound $\gamma_*$ of $\alpha$: for every integer $m \in \mathbb{Z}$ and positive integer $n \in \mathbb{N}$,
		\[
		\left|a - \frac{m}{n}\right| = \left|\alpha - \frac{k n - m}{n}\right| \ge \frac{\gamma_*}{n^2}.
		\]
		Set $\beta := 1/a > 0$. For any $p \in \mathbb{Z}$ and $q \in \mathbb{N}$, regardless of whether $p/q$ is in lowest terms, consider two cases:
		\begin{enumerate}
			\item[(i)] If $|\beta - p/q| \ge \beta/2$, then $q \ge 1$ immediately gives
			\[
			\left|\beta - \frac{p}{q}\right| \ge \frac{\beta}{2} \ge \frac{\beta}{2q^2}.
			\]
			\item[(ii)] If $|\beta - p/q| < \beta/2$, then
			\[
			\frac{1}{2a} = \frac{\beta}{2} < \frac{p}{q} < \frac{3\beta}{2} = \frac{3}{2a}.
			\]
			In particular, $p \ge 1$ is a positive integer and $a p < 3q/2$. Applying the constant-type lower bound to $a$ with denominator $p \in \mathbb{N}$ gives
			\[
			\left|\beta - \frac{p}{q}\right|
			= \frac{p}{q a}\left|a - \frac{q}{p}\right|
			\ge \frac{p}{q a} \cdot \frac{\gamma_*}{p^2}
			= \frac{\gamma_*}{q a p}
			> \frac{2\gamma_*}{3q^2}.
			\]
		\end{enumerate}
		Combining cases (i) and (ii), for any $p \in \mathbb{Z}$ and $q \in \mathbb{N}$ we have
		\[
		\left|\beta - \frac{p}{q}\right| \ge \min\!\left\{\frac{\beta}{2},\; \frac{2\gamma_*}{3}\right\} \frac{1}{q^2}.
		\]
		
		For $\omega_k=\beta/M$, since $M p \in \mathbb{Z}$ is an integer, replacing $p$ by $M p$ in the estimate for $\beta$ gives
		\[
		\left|\omega_k - \frac{p}{q}\right|
		= \frac{1}{M}\left|\beta - \frac{M p}{q}\right|
		\ge \min\!\left\{\frac{\omega_k}{2},\; \frac{2\gamma_*}{3M}\right\} \frac{1}{q^2}.
		\]
		
		If $k-\alpha\ge3/(4\gamma_*)$, then $\omega_k/2\le2\gamma_*/(3M)$, which yields the last estimate.
	\end{proof}
	
	Taking $p=0,q=1$ in~\eqref{eq:constant-type} shows that the Diophantine constant of the frequency $\omega$ is at most $|\omega|$. Thus the coefficient $\omega_k/2$ in the lower bound above is of the same order as $\omega_k$.
	
	The next lemma converts a one-sided distance estimate into a Hausdorff distance estimate.

	\begin{Lemma}\label{lem:corner-collar}
		Suppose that $f$ satisfies \textup{(G1)}--\textup{(G2)}, and let $\Sigma\subset\mathbb A$ be its separatrix. There exist an open collar neighborhood $U$ of $\Sigma$ in $\mathbb A$, a Lipschitz retraction $\mathfrak p:U\to\Sigma$ satisfying $\mathfrak p|_\Sigma=\operatorname{id}_\Sigma$, and a constant $C_\Sigma\ge1$ such that
		\begin{equation}\label{eq:retraction-bound}
			d_{\mathbb A}\bigl(z,\mathfrak p(z)\bigr)
			\le C_\Sigma\operatorname{dist}(z,\Sigma)
			\qquad(\forall z\in U),
		\end{equation}
		There is also a continuous homotopy $H:U\times[0,1]\to U$ satisfying $H_0=\operatorname{id}_U$, $H_1=\mathfrak p$, and $H_t|_\Sigma=\operatorname{id}_\Sigma$ for $0\le t\le1$. In the local coordinates $(Q, P)$ at the corner $o$, the tangential and transverse coordinates can be taken to be $u=Q-P$ and $\rho=\min(Q,P)$. In particular, $0\le \rho\le\sqrt{QP}=\sqrt v$ in the first quadrant $Q,P\ge0$.
	\end{Lemma}
	
	\begin{proof}
		Define the piecewise linear map
		\[
		F(u,\rho):=(\max(u,0)+\rho,\,\max(-u,0)+\rho).
		\]
		Its inverse is $F^{-1}(Q,P)=(Q-P,\min(Q,P))$, so $F:\mathbb R^2\to\mathbb R^2$ is a global bi-Lipschitz homeomorphism. Clearly, $0\le\rho=\min(Q,P)\le\sqrt{QP}=\sqrt v$ when $Q,P\ge0$.
		
		Fix sufficiently small $u_0\in(0,1/4)$ and $\rho_0>0$ so that the image of $[-2u_0,2u_0]\times[-\rho_0,\rho_0]$ under $F$ is contained in $V$, and define the bi-Lipschitz embedding $\widetilde\Phi_c(u,\rho):=\Phi_o^{-1}(F(u,\rho))$ there. Extending the local expression $\widetilde a(u):=\widetilde\Phi_c(u,0)$ smoothly along $\widetilde\Sigma\setminus\{o,\tau(o)\}$ gives a degree-one piecewise smooth parameterization $\widetilde a:\mathbb R\to\mathbb R^2$ of the separatrix in the universal cover, regular for $u\notin\mathbb Z$, where $\widetilde a'(u)\ne0$, and satisfying $\widetilde a(0)=o$ and $\widetilde a(u+1)=\tau(\widetilde a(u))$. Its quotient map $a:=\pi\circ\widetilde a:\mathbb T\to\Sigma$ is a homeomorphism.
		
		For $0<|u|\le2u_0$, $\widetilde\Phi_c(u,0)$ is smooth on each of $u>0$ and $u<0$, and $\det(\partial_u F(u,0),\partial_\rho F(u,0))=1$ holds on both. Together with $\det D\Phi_o^{-1}\equiv1$, this gives $\det(\widetilde a'(u),\partial_\rho\widetilde\Phi_c(u,0))=1$. On the regular part $|u|\in[u_0/2,1/2]$ away from the corner, the normal vector field $n(u):=(-\widetilde a_2'(u),\widetilde a_1'(u))^T$ is smooth and satisfies $\det(\widetilde a'(u),n(u))=|\widetilde a'(u)|^2>0$. As in the construction of a smooth tubular neighborhood (see \cite[Theorem 6.24]{Lee13}), a smooth partition of unity on the compact set $\{u\in[-1/2,1/2]:|u|\ge u_0/2\}$ gives a smooth $1$-periodic vector field $w(u)$ defined on $\mathbb R\setminus\mathbb Z$ such that $\det(\widetilde a'(u),w(u))>0$ for every $|u|\in[u_0/2,1/2]$ and $w(u)=\partial_\rho\widetilde\Phi_c(u,0)$ on the overlap $u_0\le|u|\le2u_0$. Choose a compactly supported cutoff function $\chi\in C_c^\infty((-2u_0,2u_0),[0,1])$ (see Appendix~\ref{subsec:cutoffs}) that equals $1$ on $[-u_0,u_0]$, and make the following piecewise definition on the basic strip $[-1/2,1/2]\times(-\rho_0,\rho_0)$:
		\[
		\widetilde\Phi(u,\rho):=
		\begin{cases}
			\widetilde\Phi_c(u,\rho),
			& |u|\le u_0,\\[3pt]
			(1-\chi(u))\bigl(\widetilde a(u)+\rho w(u)\bigr)
			+\chi(u)\widetilde\Phi_c(u,\rho),
			& u_0<|u|<2u_0,\\[3pt]
			\widetilde a(u)+\rho w(u),
			& 2u_0\le |u|\le1/2.
		\end{cases}
		\]
		By the support properties of $\chi$, these three expressions agree smoothly near $|u|=u_0$ and $|u|=2u_0$. Since $\widetilde a(u)+\rho w(u)$ is $\tau$-equivariant on $|u|\in[2u_0,1/2]$, the rule $\widetilde\Phi(u+k,\rho):=\tau^k(\widetilde\Phi(u,\rho))$ for $k\in\mathbb Z$ uniquely extends $\widetilde\Phi$ to a map on $\mathbb R\times(-\rho_0,\rho_0)$ and induces a quotient map $\Phi:=\pi\circ\widetilde\Phi:\mathbb T\times(-\rho_0,\rho_0)\to\mathbb A$. On the overlap $u_0<|u|<2u_0$, we have $\partial_\rho\widetilde\Phi(u,0)=w(u)$. The inverse function theorem, applied for $|u|>u_0/2$, and the bi-Lipschitz property of $\widetilde\Phi_c$ in the corner neighborhood $|u|\le u_0$ show that $\Phi$ is a local bi-Lipschitz embedding near every point of $\mathbb T\times\{0\}$. The injectivity of $\Phi|_{\mathbb T\times\{0\}}=a$ and compactness of $\mathbb T$ then give $\rho_*\in(0,\rho_0]$ such that $\Phi$ is injective on $\mathbb T\times(-\rho_*,\rho_*)$.
		
		Shrink $\rho_*$ further if necessary. Compactness of the central circle $\mathbb T\times\{0\}$, nondegeneracy of the derivative on the regular part, and the bi-Lipschitz property of $\widetilde\Phi_c$ near the corner ensure that $\Phi$ is a bi-Lipschitz homeomorphism near every point of $\mathbb T\times(-\rho_*,\rho_*)$. Hence $\Phi$ is a local homeomorphism onto its image. Taking the open strip of half the width, $U:=\Phi\bigl(\mathbb T\times(-\rho_*/2,\rho_*/2)\bigr)\subset\mathbb A$, we obtain an open set $U$ containing $\Sigma$. Define
		\[
		\begin{gathered}
			\mathfrak p\bigl(\Phi(u,\rho)\bigr):=\Phi(u,0)=a(u),\\
			H_t\bigl(\Phi(u,\rho)\bigr):=\Phi\bigl(u,(1-t)\rho\bigr)\quad(0\le t\le1).
		\end{gathered}
		\]
		Set $K:=\mathbb T\times[-\rho_*/2,\rho_*/2]$. The local bi-Lipschitz estimates and a finite covering give uniform bi-Lipschitz estimates for sufficiently close pairs of points in $K$. For pairs at distance at least a fixed positive number, compactness of $K$ and injectivity of $\Phi|_K$ give a positive lower bound for the distance between their images. Together with boundedness of $\Phi(K)$, this shows that $\Phi|_K$ and its inverse are Lipschitz. Thus $\mathfrak p:U\to\Sigma$ is Lipschitz, and $H_t$ is a continuous homotopy fixing $\Sigma$ pointwise. Finally, for any $z\in U$, choose $b\in\Sigma$ such that $d_{\mathbb A}(z,b)=\operatorname{dist}(z,\Sigma)$. By $\mathfrak p(b)=b$,
		\[
		d_{\mathbb A}\bigl(z,\mathfrak p(z)\bigr)
		\le d_{\mathbb A}(z,b)+d_{\mathbb A}\bigl(\mathfrak p(b),\mathfrak p(z)\bigr)
		\le\bigl(1+\operatorname{Lip}(\mathfrak p)\bigr)\operatorname{dist}(z,\Sigma),
		\]
		Taking $C_\Sigma:=1+\operatorname{Lip}(\mathfrak p)$ completes the proof.
	\end{proof}

	\begin{proof}[Proof of Theorem~\ref{thm:boundary-transition-criterion}]
		In addition to the threshold chosen in Section~\ref{subsec:auto-4-3}, require $k_f-\alpha\ge3/(4\gamma_*)$, and increase $k_f$ as needed in the distance estimates below.
		
		\paragraph{Immersion and closure of the extended curve.}
		Fix $k \ge k_f$ and write $\varrho := \varrho_k = \alpha - k < 0$. Define the extended curve in the universal cover $\mathbb{R}^2$ by
		\[
		\ell(t) := \iota\bigl(\widetilde{K}_k(t)\bigr) \qquad (t \in \mathbb{R}).
		\]
		Since $D\iota$ is invertible everywhere and~\eqref{eq:Kk-def} gives
		\[
		\widetilde{K}_k'(t) = DD_k\bigl(K_k(t)\bigr)K_k'(t) \ne (0, 0),
		\]
		\[
		\ell'(t) = D\iota\bigl(\widetilde{K}_k(t)\bigr)\widetilde{K}_k'(t) \ne (0, 0) \qquad (\forall t \in \mathbb{R}),
		\]
		Thus $\ell$ is a smooth regular immersion. Combining~\eqref{eq:R-tilde-Kk} with~\eqref{eq:time-equivariance} and~\eqref{eq:spatial-equivariance}, for every $t \in \mathbb{R}$ we have
		\begin{equation}\label{eq:ell-equivariance}
			\ell(t + 1) = \widetilde{f}\bigl(\ell(t)\bigr), \qquad
			\ell(t + \varrho) = \tau^{-1}\bigl(\ell(t)\bigr).
		\end{equation}
		Set $a := -\varrho = k - \alpha > 0$, and let $c(t) := \pi\bigl(\ell(t)\bigr)$ be the natural projection of $\ell(t)$ to the cylinder $\mathbb{A}$. Equation~\eqref{eq:ell-equivariance} gives
		\begin{equation}\label{eq:c-properties}
			c(t + a) = c(t), \qquad
			c(t + 1) = f\bigl(c(t)\bigr), \qquad
			\ell(t + a) = \tau\bigl(\ell(t)\bigr).
		\end{equation}
		Here $c(t+a)=c(t)$ means that the projected curve is closed, whereas $\ell(t+a)=\tau(\ell(t))$ means that it winds once around the cylinder. We next prove that the quotient map $c : \mathbb{R}/(a\mathbb{Z}) \to \mathbb{A}$ is injective.
		
		\paragraph{Injectivity of the quotient map.}
		Suppose that $t_0, d \in \mathbb{R}$ satisfy $c(t_0) = c(t_0 + d)$. Since $f$ is invertible, the second identity in~\eqref{eq:c-properties} implies that, for every $n \in \mathbb{Z}$,
		\[
		c(t_0 + n) = f^n\bigl(c(t_0)\bigr) = f^n\bigl(c(t_0 + d)\bigr) = c(t_0 + d + n).
		\]
		Since $a=k-\alpha$ is irrational, the set $\{t_0+n-ma:m,n\in\mathbb Z\}$ is dense in the real line. For any $s \in \mathbb{R}$, there is a sequence of integers $\{n_j\}, \{m_j\}$ such that $t_0 + n_j - m_j a \to s$ as $j \to \infty$. Using the $a$-periodicity and continuity of $c$ and taking the limit gives
		\[
		c(s) = \lim_{j \to \infty} c(t_0 + n_j - m_j a)
		= \lim_{j \to \infty} c(t_0 + d + n_j - m_j a) = c(s + d).
		\]
		Hence $d$ is a global period of $c$.

		Consider the set of periods of $c$,
		\[
		H := \bigl\{h \in \mathbb{R} : c(t + h) = c(t) \text{ for every } t \in \mathbb{R}\bigr\}.
		\]
		Continuity of $c$ and the first identity in~\eqref{eq:c-properties} show that $H$ is a closed subgroup of $(\mathbb{R}, +)$ containing $a\mathbb{Z}$. Fix any point $t_* \in \mathbb{R}$. In local coordinates near $c(t_*)$, the condition $c'(t_*) \ne (0, 0)$ implies that one coordinate component is strictly monotone on some small interval $(t_* - \delta, t_* + \delta)$ with $\delta > 0$. Thus
		\[
		H \cap (-\delta, \delta) = \{0\}.
		\]
		The set $H \cap (0, \infty)$ is therefore nonempty and has a positive lower bound, and its infimum satisfies $h_0 := \inf\bigl(H \cap (0, \infty)\bigr) \ge \delta > 0$. Since $H$ is closed, this infimum is attained: $H \cap (0, \infty)$ has a smallest positive element $h_0 > 0$. For any $h \in H$, division with remainder gives a unique integer $m \in \mathbb{Z}$ and a remainder $r \in [0, h_0)$ such that $h = m h_0 + r$. If $r > 0$, then $r \in H \cap (0, \infty)$ and $r<h_0$, contradicting the minimality of $h_0$. Hence $r=0$, and therefore $h = m h_0 \in h_0 \mathbb{Z}$. This proves $H = h_0\mathbb{Z}$. Consequently, there is a unique positive integer $N \ge 1$ such that $a = N h_0$.

		For every $t \in \mathbb{R}$, equality of the projections $c(t + h_0) = c(t)$ gives an integer $d(t) \in \mathbb{Z}$ such that
		\[
		\ell(t + h_0) = \tau^{d(t)}\bigl(\ell(t)\bigr), \qquad
		d(t) = \pi_x\bigl(\ell(t + h_0)\bigr) - \pi_x\bigl(\ell(t)\bigr).
		\]
		Since $\ell$ is continuous, $d(t)$ is a continuous integer-valued function on the real line. Thus $d(t)$ is constant; denote it by $d_0 \in \mathbb{Z}$. Iterating this identity $N$ times and using the third identity in~\eqref{eq:c-properties} gives
		\[
		\tau^{N d_0}\bigl(\ell(t)\bigr) = \ell(t + N h_0) = \ell(t + a) = \tau\bigl(\ell(t)\bigr), \qquad N d_0 = 1.
		\]
		Since $N \ge 1$ and $d_0 \in \mathbb{Z}$, we have $N = d_0 = 1$, and hence $H = a\mathbb{Z}$. Together with the preceding discussion of self-intersections, this proves that $c : \mathbb{R}/(a\mathbb{Z}) \to \mathbb{A}$ is injective.
		\paragraph{Embedding and rotation number.}
		A continuous injection from the compact space $\mathbb{R}/(a\mathbb{Z})$ into the Hausdorff space $\mathbb{A}$ is a homeomorphism onto its image. Since $c$ is a smooth regular immersion, the local inverse function theorem shows that its image $\Gamma_k := c(\mathbb{R})$ is a smoothly embedded circle in the cylinder $\mathbb{A}$. The third identity in~\eqref{eq:c-properties}, namely $\ell(t + a) = \tau\bigl(\ell(t)\bigr)$, shows that the projection degree of $\Gamma_k$ is $1$. Hence $\Gamma_k$ is a smooth essential invariant circle.
		
		Since $c:\mathbb R/(a\mathbb Z)\to\mathbb A$ is a smooth embedding and
		$\ell(t+a)=\tau\bigl(\ell(t)\bigr)=\ell(t)+(1,0)$,
		its lift $\ell:\mathbb R\to\mathbb R^2$ is also a smooth embedding.
		Moreover, $\pi\bigl(\ell(\mathbb R)\bigr)=\Gamma_k$, and
		$\ell(\mathbb R)$ is invariant under every integer horizontal translation $\tau^n$, $n\in\mathbb Z$.
		Therefore,
		\[
		\ell(\mathbb R)=\pi^{-1}(\Gamma_k).
		\]
		Thus the full preimage of $\Gamma_k$ in the universal cover is precisely this smoothly embedded line.
		
		Finally, define the parameterization $\gamma_k(s) := \ell(a s) = \ell(-\varrho s)$. On the universal cover $\mathbb{R}^2$, equations~\eqref{eq:ell-equivariance} and~\eqref{eq:c-properties} directly give
		\[
		\gamma_k(s + 1) = \tau\bigl(\gamma_k(s)\bigr), \qquad
		\widetilde{f}\bigl(\gamma_k(s)\bigr) = \ell(a s + 1) = \gamma_k\!\left(s + \frac{1}{a}\right) = \gamma_k\!\left(s - \frac{1}{\varrho}\right).
		\]
		Thus, with respect to the chosen lift $\widetilde{f}$, the rotation number of $\Gamma_k$ is $\omega_k = -1/\varrho = 1/(k - \alpha) > 0$. This proves~\eqref{eq:rotation-conjugacy} in Theorem~\ref{thm:boundary-transition-criterion}.
		\paragraph{Distance to the separatrix.}
		In the $(\theta, v)$ coordinates, the invariant graph of $\mathcal R$ can be written as
		\[
		\kappa_k(\theta) := D_k\bigl(\kappa_k^{\mathrm{sc}}(\theta)\bigr) = \bigl(\theta,\, V_k(\theta)\bigr), \qquad
		V_k(\theta) := e^{-k\nu}\bigl(s_\alpha + \psi_k(\theta)\bigr),
		\]
		where $V_k$ is $1$-periodic. By~\eqref{eq:s-alpha-def} and~\eqref{eq:psi-k-C1-bound}, there exist positive numbers $c_1 := s_\alpha - d_s/2 > 1/2$ and $c_2 := s_\alpha + d_s/2 < 1$ independent of $k$ such that
		\begin{equation}\label{eq:Vk-bounds}
			c_1 e^{-k\nu} \le V_k(\theta) \le c_2 e^{-k\nu} \qquad (\forall \theta \in \mathbb{R}).
		\end{equation}
		The induced lift $r_k(\theta) := \pi_1\Bigl(\mathcal{R}\bigl(\theta,V_k(\theta)\bigr)\Bigr) = \pi_1\Bigl(\mathcal{F}_k\bigl(\theta,s_\alpha+\psi_k(\theta)\bigr)\Bigr) - k$ defined by $\mathcal{R}\bigl(\kappa_k(\theta)\bigr) = \kappa_k\bigl(r_k(\theta)\bigr)$ is conjugate to the translation $t \mapsto t + \varrho_k$ through the increasing first component $h_k$ of $\widetilde{K}_k$. Hence $r_k$ is strictly increasing on $\mathbb{R}$. By~\eqref{eq:Fk-explicit} and~\eqref{eq:Vk-bounds}, there is a constant $C_0 > 0$ independent of $k$ such that
		\begin{equation}\label{eq:rk-shift-bound}
			\bigl|r_k(\theta) - \theta + k\bigr| \le C_0 \qquad (\forall \theta \in \mathbb{R}).
		\end{equation}
		By~\eqref{eq:rk-shift-bound}, when $k > C_0 + 1$ we have $r_k(\theta) < \theta - 1$ for every $\theta \in \mathbb{R}$. Choose the left endpoint $\theta_{\mathrm{in}}^-$ of the incoming basic phase interval $B_{\mathrm{in}}$ as the reference point, and take $\theta = r_k^{-1}(\theta_{\mathrm{in}}^-)$ in~\eqref{eq:rk-shift-bound} to obtain
		\[
		\theta_{\mathrm{in}}^- + k - C_0 \le r_k^{-1}(\theta_{\mathrm{in}}^-) \le \theta_{\mathrm{in}}^- + k + C_0.
		\]
		Thus the closed interval
		\[
		J_k := \bigl[\theta_{\mathrm{in}}^-,\, r_k^{-1}(\theta_{\mathrm{in}}^-)\bigr]
		\]
		is a basic interval for the action of $r_k$ on $\mathbb{R}$: the intervals $\{r_k^n(J_k)\}_{n \in \mathbb{Z}}$ cover the real line $\mathbb{R}$ and have pairwise disjoint interiors. By~\eqref{eq:Vk-bounds}, the entire invariant graph $\kappa_k(\mathbb{R})$ lies in the substrip $\mathbb{R} \times (0, v_*/2)$, so~\eqref{eq:spatial-equivariance} can be applied repeatedly on this graph. Since the projections of $\iota\bigl(\kappa_k(r_k^n(J_k))\bigr) = \tau^{-n}\Bigl(\iota\bigl(\kappa_k(J_k)\bigr)\Bigr)$ to the cylinder $\mathbb{A}$ all coincide, that is, $\Gamma_k = \pi\bigl(\iota(\kappa_k(J_k))\bigr)$, it suffices to estimate the distance from $\iota\bigl(\kappa_k(J_k)\bigr)$ to the separatrix $\Sigma$ on the finite interval $J_k$.
		
		For any $\theta \in J_k$, write $\theta = \phi + n$, where $\phi \in [\theta_{\mathrm{in}}^-, \theta_{\mathrm{in}}^- + 1)$ and the integer $n$ satisfies $0 \le n \le k + C_1$; take $C_1 := \lceil C_0 \rceil$. By~\eqref{eq:iota-def},
		\[
		\iota\bigl(\kappa_k(\theta)\bigr) = \widetilde{f}^n\Bigl(J_{\mathrm{in}}\bigl(\phi,V_k(\phi)\bigr)\Bigr).
		\]
		Fix $q_* > 0$ so that the basic interval $[q_*, e^\nu q_*] \times \{0\}$ on the local outgoing axis and its image under one further iterate both lie strictly inside the closed linearizing box $\mathcal{B}$. For an incoming point $J_{\mathrm{in}}\bigl(\phi, V_k(\phi)\bigr)$, define its residence time in the linearizing box by
		\[
		n_*(\phi, k) := \left\lceil \frac{\log q_* - \log V_k(\phi)}{\nu} - \phi \right\rceil.
		\]
		For $0 \le j \le n_*$, its local linearizing coordinates are
		\[
		Q_j = e^{j\nu} V_k(\phi) e^{\nu\phi}, \qquad
		P_j = e^{-j\nu} e^{-\nu\phi}, \qquad
		Q_j P_j = V_k(\phi).
		\]
		The definition of the ceiling function and~\eqref{eq:Vk-bounds} give a constant $C_2 > 0$ independent of $k, \phi$ such that
		\[
		q_* \le Q_{n_*} < e^\nu q_*, \qquad
		|n_* - k| \le C_2.
		\]
		For every $0 \le j \le n_*$, since $Q_j \le e^\nu q_*$ and $P_j \le e^{-\nu\theta_{\mathrm{in}}^-}$, the point $(Q_j, P_j)$ remains in the closed linearizing box $\mathcal{B}$. In the $(Q, P)$ coordinates, its Euclidean distance to the two local separatrix half-axes $\{QP = 0,\, Q, P \ge 0\}$ is $\min(Q_j, P_j)$. Thus
		\[
		\begin{aligned}
			\operatorname{dist}\bigl((Q_j, P_j),\, \{QP = 0,\, Q, P \ge 0\}\bigr)
			&= \min(Q_j, P_j)\\
			&\le \sqrt{Q_j P_j}
			= \sqrt{V_k(\phi)} \le \sqrt{c_2}\,e^{-k\nu/2}.
		\end{aligned}
		\]
		Since the derivatives of the local symplectic chart and its inverse are bounded on the compact box $\mathcal{B}$, the same estimate holds, up to a fixed multiplicative constant, for distance on the cylinder $\mathbb{A}$.
		
		If $n > n_*$, then the remaining number of iterates satisfies $0 < n - n_* \le C_1 + C_2$. At time $n_*$, the orbit point has the same horizontal coordinate $Q_{n_*} \in [q_*, e^\nu q_*)$ as the point $(Q_{n_*}, 0)$ on the compact outgoing separatrix segment, and its vertical coordinate satisfies
		\[
		0 < P_{n_*} = \frac{V_k(\phi)}{Q_{n_*}} \le \frac{c_2}{q_*} e^{-k\nu}.
		\]
		Fix an integer $J_* > C_1 + C_2$, and on a sufficiently small compact neighborhood of this compact outgoing separatrix segment set
		\[
		M_* := \max_{0 \le j \le J_*} \sup \|D\widetilde{f}^j\| < \infty.
		\]
		By the mean value inequality, for $1 \le n - n_* \le J_*$ each subsequent iterate is at distance at most $C_f M_* e^{-k\nu}$ from the corresponding iterate on the separatrix, whose projection belongs to $\Sigma$. Combining the cases $0 \le n \le n_*$ and $n > n_*$ gives the one-sided distance estimate
		\begin{equation}\label{eq:one-sided-dist}
			\sup_{z \in \Gamma_k} \operatorname{dist}(z, \Sigma) \le C_f e^{-k\nu/2}.
		\end{equation}
		
		Furthermore, the dynamics on $\Gamma_k$ is conjugate to an irrational rotation, so this circle does not contain the fixed point $o$. If $\Gamma_k$ and $\Sigma \setminus \{o\}$ intersect at a point $z$, its forward orbit satisfies $f^n(z) \to o$ as $n \to +\infty$. Forward invariance and closedness of $\Gamma_k$ would imply $o \in \Gamma_k$, a contradiction. {
Hence $\Gamma_k\cap\Sigma=\varnothing$. The embedded essential curve $\Sigma$ separates the cylinder into two components. The connected curve $\Gamma_k$ meets the chosen component, because $\iota=J_{\mathrm{in}}$ on the incoming basic rectangle and $V_k>0$ corresponds to $Q,P>0$. Consequently, the whole of $\Gamma_k$ lies on the chosen side of $\Sigma$.
}
		
		To estimate the reverse distance, take the neighborhood $U$ and map $\mathfrak p:U\to\Sigma$ from Lemma~\ref{lem:corner-collar}.
		
		For sufficiently large $k$, estimate~\eqref{eq:one-sided-dist} ensures that $\Gamma_k \subset U$. Since $\mathfrak p\simeq\operatorname{id}_U$, this map preserves the winding number of a closed curve. The essential circle $\Gamma_k$ has degree $1$, so the restriction $\mathfrak{p}|_{\Gamma_k} : \Gamma_k \to \Sigma$ has degree $1$ and is therefore surjective. Thus, for every $x \in \Sigma$, there is $z \in \Gamma_k$ such that $\mathfrak{p}(z) = x$. Combining~\eqref{eq:one-sided-dist} and~\eqref{eq:retraction-bound} gives
		\[
		\operatorname{dist}(x, \Gamma_k) \le d_{\mathbb{A}}(x, z) = d_{\mathbb{A}}\bigl(\mathfrak{p}(z), z\bigr)
		\le C_\Sigma \operatorname{dist}(z, \Sigma) \le C_\Sigma C_f e^{-k\nu/2}.
		\]
		Combining this with~\eqref{eq:one-sided-dist} proves the Hausdorff distance estimate
		\begin{equation}\label{eq:hausdorff-bound}
			d_H(\Gamma_k, \Sigma) \le C_f' e^{-k\nu/2} \longrightarrow 0 \qquad (k \to \infty).
		\end{equation}
		Increasing $C_f$ to absorb the constant $C_f'$ in~\eqref{eq:hausdorff-bound} gives~\eqref{eq:Gamma-Sigma-dist}. Taking $M=1$ in Lemma~\ref{lem:reciprocal-frequency} gives~\eqref{eq:omegak-estimate}. The uniform choices of $Q_*$ and $\delta_{\nu,\gamma_\xi}$ were explained in Section~\ref{subsec:auto-4-3}. This completes the proof.
	\end{proof}
	
	\subsection{Proof of Corollary~\ref{cor:f0-invariant}}\label{subsec:auto-5-3}
	\begin{proof}
		Apply Theorem~\ref{thm:boundary-transition-criterion} to the iterate $G=f_0^M$ in place of $f$, fix $k\ge k_G$, and retain $\ell$ and $\varrho=\alpha-k$ from its proof. For $0 \le i < M$, write $\Gamma_{k, i} := f_0^i(\Gamma_k)$. Since $G \circ f_0^i = f_0^i \circ G$, the restrictions $G|_{\Gamma_{k, i}}$ are smoothly conjugate to the same irrational rotation. Thus every $G$-orbit on each $\Gamma_{k, i}$ is dense in that circle. If there is a point $z \in \Gamma_{k, i} \cap \Gamma_{k, j}$, taking orbit closures gives
		\[
		\overline{\{G^n(z) : n \in \mathbb{Z}\}} = \Gamma_{k, i} = \Gamma_{k, j}.
		\]
		Consequently, the distinct members of the finite family $\{\Gamma_{k, i} : 0 \le i < M\}$ are pairwise disjoint.
		
		Every embedded essential circle in the cylinder $\mathbb{A}$ separates its upper and lower ends, so a finite collection of pairwise disjoint essential circles has a natural total order: write $\Gamma < \Gamma'$ if $\Gamma'$ lies in the upper connected component of $\mathbb{A} \setminus \Gamma$. List the distinct members of $\{\Gamma_{k, i} : 0 \le i < M\}$ in this order as $C_1 < C_2 < \dots < C_d$, with $1 \le d \le M$. Since $f_0^M(\Gamma_k) = G(\Gamma_k) = \Gamma_k$, the invertible map $f_0$ induces a permutation of $\{C_1, \dots, C_d\}$. Since $f_0$ is isotopic to the identity, it maps the upper and lower sides of each essential circle to the upper and lower sides of its image, respectively. This permutation therefore strictly preserves the total order $<$. A strictly order-preserving permutation of a finite totally ordered set fixes its least element $C_1$ and, by induction, every element $C_j$. In particular,
		\begin{equation}\label{eq:f0-invariance}
			f_0(\Gamma_k) = \Gamma_k.
		\end{equation}
		
		Since $\Gamma_k$ has degree $1$, Section~\ref{subsec:auto-5-2} shows that its full preimage in the universal cover $\mathbb{R}^2$ is a single connected embedded line $\ell(\mathbb{R})$. Thus, for the chosen fixed-point lift $\widetilde{f}_0$ satisfying $\widetilde{G} = \widetilde{f}_0^M$, we also have $\widetilde{f}_0\bigl(\ell(\mathbb{R})\bigr) = \ell(\mathbb{R})$. Since $\ell : \mathbb{R} \to \ell(\mathbb{R})$ is a homeomorphism, there is a unique continuous function $\mathfrak{a} : \mathbb{R} \to \mathbb{R}$ such that
		\[
		\widetilde{f}_0\bigl(\ell(t)\bigr) = \ell\bigl(\mathfrak{a}(t)\bigr) \qquad (\forall t \in \mathbb{R}).
		\]
		The lift $\widetilde{f}_0$ commutes with both $\widetilde{G}$ and $\tau^{-1}$, so
		
		\[ \ell\bigl(\mathfrak{a}(t+1)\bigr) = \widetilde{f}_0\bigl(\widetilde{G}(\ell(t))\bigr) = \widetilde{G}\bigl(\widetilde{f}_0(\ell(t))\bigr) = \ell\bigl(\mathfrak{a}(t)+1\bigr),\]
		
		\[\ell\bigl(\mathfrak{a}(t+\varrho)\bigr) = \widetilde{f}_0\bigl(\tau^{-1}(\ell(t))\bigr) = \tau^{-1}\bigl(\widetilde{f}_0(\ell(t))\bigr) = \ell\bigl(\mathfrak{a}(t)+\varrho\bigr).\]
		Injectivity of $\ell$ gives
		
		\[
		\mathfrak{a}(t + 1) = \mathfrak{a}(t) + 1, \qquad
		\mathfrak{a}(t + \varrho) = \mathfrak{a}(t) + \varrho.
		\]
		Thus the continuous function $d(t) := \mathfrak{a}(t) - t$ has both the period $1$ and the irrational period $\varrho = \alpha - k$. Since $\{m + n\varrho : m, n \in \mathbb{Z}\}$ is dense in $\mathbb{R}$, continuity gives a constant $c_0 \in \mathbb{R}$ such that $d(t) \equiv c_0$, that is, $\mathfrak{a}(t) = t + c_0$. Iterating $\widetilde{f}_0\bigl(\ell(t)\bigr) = \ell(t + c_0)$ a total of $M$ times and comparing with $\widetilde{G}\bigl(\ell(t)\bigr) = \ell(t + 1)$ gives
		\[
		\widetilde{f}_0^M\bigl(\ell(t)\bigr) = \ell(t + M c_0), \qquad
		\widetilde{G}\bigl(\ell(t)\bigr) = \ell(t + 1).
		\]
		By $\widetilde{G} = \widetilde{f}_0^M$ and injectivity of $\ell$, we have $M c_0 = 1$, or equivalently $c_0 = 1/M$. Substitution into $\gamma_k(s) = \ell(-\varrho s)$ yields
		\[
		\widetilde{f}_0\bigl(\gamma_k(s)\bigr) = \ell\bigl(-\varrho s + \frac{1}{M}\bigr) = \gamma_k\!\left(s - \frac{1}{M\varrho}\right).
		\]
		Hence the rotation number of the chosen lift $\widetilde{f}_0$ of the original map on $\Gamma_k$ is
		\[
		\omega_k = -\frac{1}{M\varrho} = \frac{1}{M(k - \alpha)},
		\]
		This proves the rotation number formula~\eqref{eq:f0-rotation-number} in Corollary~\ref{cor:f0-invariant}.
		Lemma~\ref{lem:reciprocal-frequency} then gives~\eqref{eq:omegak-estimate}. The condition on $k-\alpha$ is independent of $M$, so the threshold $k_G$ from the theorem applied to $G$ can be retained. This proves the corollary.
	\end{proof}
	
	\subsection{Proof of Corollary~\ref{cor:graph-corollary}}\label{subsec:auto-5-4}
	\begin{proof}
		By Theorem~\ref{thm:boundary-transition-criterion} and Corollary~\ref{cor:f0-invariant}, $\Gamma_k$ is a smooth essential invariant circle of $f_0$. The fact stated in Section~\ref{subsec:auto-1-1} shows that $f_0$ is exact symplectic. Together with the positive twist condition, Birkhoff's theorem gives $\Gamma_k=\operatorname{graph}u_k$, where $u_k$ is Lipschitz~\cite{Arn09}.
		
		We now show that this Lipschitz graph is in fact $C^\infty$ smooth. Write a lift of the smooth regular degree-one parameterization of $\Gamma_k$ in component form:
		\[
		\begin{gathered}
			\gamma(t) = \bigl(X(t),\, Y(t)\bigr),\qquad
			\bigl(X'(t),\, Y'(t)\bigr) \ne (0, 0),\\
			X(t + 1) = X(t) + 1,\qquad
			Y(t + 1) = Y(t).
		\end{gathered}
		\]
		Write the invariant Lipschitz graph supplied by Birkhoff's theorem as $y = u(x)$, and set $L_u := \operatorname{Lip}(u) < \infty$. Since points on the curve satisfy $Y(t) = u\bigl(X(t)\bigr)$, for any $t \in \mathbb{R}$ and $h \ne 0$,
		\[
		\left|\frac{Y(t + h) - Y(t)}{h}\right|
		= \frac{\bigl|u(X(t + h)) - u(X(t))\bigr|}{|h|}
		\le L_u \left|\frac{X(t + h) - X(t)}{h}\right|.
		\]
		Letting $h \to 0$ gives $|Y'(t)| \le L_u |X'(t)|$. If $X'(t_0) = 0$ at some point $t_0$, this inequality would give $Y'(t_0) = 0$, contradicting regularity $\bigl(X'(t_0), Y'(t_0)\bigr) \ne (0, 0)$. Thus $X'(t)$ never vanishes on $\mathbb{R}$ and has constant sign. Together with
		\[
		\int_0^1 X'(t)\,\mathrm{d}t = X(1) - X(0) = 1 > 0,
		\]
		this implies that $X'(t) > 0$ for every $t \in \mathbb{R}$. The $1$-periodicity of $X'$ and compactness of $[0, 1]$ give
		\[
		a_\gamma := \min_{t \in [0, 1]} X'(t) > 0.
		\]
		By the inverse function theorem, the strictly increasing smooth function $X : \mathbb{R} \to \mathbb{R}$ has a smooth inverse $h := X^{-1} : \mathbb{R} \to \mathbb{R}$ satisfying
		\[
		h(x + 1) = h(x) + 1, \qquad
		h'(x) = \frac{1}{X'\bigl(h(x)\bigr)}, \qquad
		u(x) = Y\bigl(h(x)\bigr).
		\]
		Therefore $u \in C^\infty(\mathbb{T})$. In particular, its first and second derivatives are
		\[
		u'(x) = \frac{Y'\bigl(h(x)\bigr)}{X'\bigl(h(x)\bigr)}, \qquad
		u''(x) = \frac{Y''X' - Y'X''}{(X')^3}\bigl(h(x)\bigr),
		\]
		Repeated differentiation gives smooth derivatives of every order. The positive lower bound $a_\gamma > 0$ is needed only for each fixed invariant circle and need not be uniform in $k$.
		
		Finally, let $g_u(x) := \pi_x\Bigl(\widetilde{f}_0\bigl(x,u(x)\bigr)\Bigr)$ be the lift of the map induced by $f_0$ on the invariant graph $\Gamma_k = \operatorname{graph} u$. Taking the first component of $\widetilde{f}_0\bigl(\gamma(t)\bigr) = \gamma(t + \omega_k)$ gives
		\[
		g_u\bigl(X(t)\bigr) = X(t + \omega_k).
		\]
		This proves the assertions in Corollary~\ref{cor:graph-corollary} concerning the smooth graph and the smooth conjugacy of the circle dynamics.

		Apply Theorem~\ref{thm:boundary-transition-criterion} to $G=f_0^M$, and denote the constant in the distance estimate by $C_G$. If the separatrix $\Sigma=\operatorname{graph}\psi$ is also a Lipschitz graph, set $L_\psi := \operatorname{Lip}(\psi)$. Since the invariant graph satisfies $V_k > 0$, $\Gamma_k$ intersects the local first quadrant $\{Q > 0, P > 0\}$. Moreover, $\{Q > 0, P > 0\}$ lies above $\Sigma$, $\Gamma_k \cap \Sigma = \varnothing$, and $\Gamma_k$ is connected. Hence $\Gamma_k$ lies entirely above $\Sigma$, that is, $u_k(x) > \psi(x)$ for every $x \in \mathbb{T}$.
		
		For any $x \in \mathbb{T}$, set $z = (x, u_k(x))$ and choose $b = (x_0, \psi(x_0)) \in \Sigma$ such that $d_{\mathbb{A}}(z, b) = \operatorname{dist}(z, \Sigma)$. Lipschitz continuity of $\psi$ and the triangle inequality give
		\[
		\begin{aligned}
			0 < u_k(x) - \psi(x) &\le |u_k(x) - \psi(x_0)| + |\psi(x_0) - \psi(x)| \\
			&\le |u_k(x) - \psi(x_0)| + L_\psi d_{\mathbb{T}}(x, x_0) \\
			&\le (1 + L_\psi)d_{\mathbb{A}}(z, b)
			= (1 + L_\psi)\operatorname{dist}(z, \Sigma)\\
			&\le (1 + L_\psi) C_G e^{-k\nu/2}.
		\end{aligned}
		\]
		This proves~\eqref{eq:graph-distance-estimate} and completes the proof of the corollary.
	\end{proof}
	
\begin{proof}[Proof of Theorem~\ref{thm:intro-criterion}]
Let $f_0$ be the original map and choose $M\in\mathbb N$ so that $G=f_0^M$, with the lift $\widetilde G=\widetilde f_0^M$, satisfies the stated assumptions. Identity~\eqref{eq:transition-distance-translations} converts the assumed $C^4$-closeness to a translation into~\eqref{eq:zeta-C4-assumption}. Theorem~\ref{thm:boundary-transition-criterion} applied to $G$ gives smooth essential invariant circles on the chosen side, together with exponential Hausdorff convergence to the separatrix. Corollary~\ref{cor:f0-invariant} makes these circles invariant under $f_0$ and gives smooth conjugacy to rotations with constant-type frequencies $1/[M(k-\alpha)]$ tending to zero. The same circles retain the Hausdorff estimate. If $f_0$ has positive twist bounded away from zero, Corollary~\ref{cor:graph-corollary}(1) makes them smooth graphs. If $\Sigma=\operatorname{graph}\psi$ is Lipschitz, the argument of part~(2) gives $\|u_k-\psi\|_\infty\le(1+\operatorname{Lip}\psi)d_H(\Gamma_k,\Sigma)$. One-sidedness determines the sign of $u_k-\psi$, so convergence from the chosen side is uniform and exponential. The uniformity of the one-branch threshold follows from the last assertion of Theorem~\ref{thm:boundary-transition-criterion}. The positive-area conclusion follows from Remark~\ref{rem:scope}.
\end{proof}

\subsection{The return map near the heteroclinic cycle}\label{subsec:island-proof}
We now prove Proposition~\ref{prop:island}, which establishes the interior assertion of Theorem~\ref{thm:intro-island}. Work in the universal cover $\mathbb R^2$ of the cylinder, and write $D$ for the island, $o_0=o$ and $o_1=\tau(o_0)$ for its saddle vertices, and $G=f_0^M$. The boundary $\partial D=S_+\cup S_-$ is oriented along the upper arc $S_+$ from $o_0$ to $o_1$ and along the lower arc $S_-$ from $o_1$ to $o_0$; rotation numbers about an interior fixed point are taken positive in this direction. Since the integer translates $\tau^n(D)$ are pairwise disjoint, the cylinder projection $\pi|_D$ is injective.

\paragraph{The two passages.}
On the interior side of either saddle ($Q>0,P<0$ near $o_0$ and $Q<0,P>0$ near $o_1$), set $w=-QP=|Q||P|>0$ and define the outgoing and incoming logarithmic phases by
\[
\theta_{\mathrm{out}}=\nu^{-1}\log|Q|,\qquad
\theta_{\mathrm{in}}=-\nu^{-1}\log|P|.
\]
In each of the four interior charts, $\widetilde G$ is conjugate to the unit translation $T_1(\theta,w)=(\theta+1,w)$, and the local change of coordinates from the incoming chart to the outgoing chart at either vertex is given by the same map $\mathcal C_\nu$ as in Section~\ref{subsec:auto-3-1}:
\[
\mathcal C_\nu(\theta,w)=\left(\theta+\frac{\log w}{\nu},\,w\right).
\]
Index the local charts by $0$ at $o_0$ and by $1$ at $o_1$. Direct differentiation shows that all four inverse charts $J_{0,\mathrm{in}},J_{0,\mathrm{out}},J_{1,\mathrm{in}},J_{1,\mathrm{out}}$ pull back $\mathrm dQ\wedge\mathrm dP$ to $-\nu\,\mathrm d\theta\wedge\mathrm dw$, and that $\mathcal C_\nu^*(w\,\mathrm d\theta)-w\,\mathrm d\theta=\nu^{-1}\mathrm dw$.

Applying Lemma~\ref{lem:finite_transition} and the extension procedure of Section~\ref{subsec:auto-3-2} separately along $S_+$ (from $o_0$ to $o_1$ in $L_+$ iterates of $G$) and along $S_-$ (from $o_1$ to $o_0$ in $L_-$ iterates of $G$), and subtracting $L_+$ and $L_-$ from the respective terminal phases as in~\eqref{eq:normalized-transition-map}, we obtain smooth normalized transition embeddings $\mathcal P_+,\mathcal P_-$ on $\mathbb R\times[0,w_*)$ commuting with $T_1$ and preserving $\mathrm d\theta\wedge\mathrm dw$. On the boundary $w=0$, we have $\mathcal P_i(\theta,0)=(b_i(\theta),0)$ with $b_i'>0$ for $i\in\{+,-\}$. By the proof of Lemma~\ref{lem:return_lift_exact}, each $\mathcal P_i$ is exact symplectic with respect to $w\,\mathrm d\theta$. Define the one-pass return lifts and the full-loop return lift by
\[
\mathcal R_i:=\mathcal P_i\circ\mathcal C_\nu\quad(i=+,-),\qquad
\mathcal R_{\mathrm{in}}:=\mathcal R_-\circ\mathcal R_+ .
\]
Since $\partial_w(\mathcal P_+)_2(\theta,0)=1/b_+'(\theta)$ is positive and $1$-periodic, shrinking $w_*>0$ ensures that $\mathcal R_+$ maps $\mathbb R\times(0,w_*)$ into the domain of $\mathcal R_-$, so $\mathcal R_{\mathrm{in}}$ is a smooth exact symplectic embedding on $\mathbb R\times(0,w_*)$ commuting with $T_1$.

Write $b_i(\theta)=\theta+\sigma_i+\zeta_i(\theta)$ with $\int_0^1\zeta_i(\theta)\,\mathrm d\theta=0$, and set $\delta:=\max_{i\in\{+,-\}}\|\zeta_i\|_{C^4(\mathbb T)}$. For $\eta_k=e^{-k\nu}$ and $D_k(\theta,s)=(\theta,\eta_k s)$, define
\[
\mathcal F_{i,k}:=T_k\circ D_k^{-1}\circ\mathcal R_i\circ D_k,\qquad
\mathcal T_{\Lambda_i}(\theta,s):=\left(\theta+\sigma_i+\frac{\log s}{\nu},\,s\right)\qquad(i=+,-).
\]
Assume $\delta\le 1/2$. On any fixed compact subinterval of $(0,\infty)$, the calculation of Section~\ref{subsec:auto-4-1} (with $w$ in place of $v$ and $z=\theta+\nu^{-1}\log s$) gives
\[
\mathcal F_{i,k}(\theta,s)=
\begin{pmatrix}
	\theta+\sigma_i+\nu^{-1}\log s+\zeta_i(z)+O_{G,\nu}(\eta_k)\\[4pt]
	s/b_i'(z)+O_{G,\nu}(\eta_k)
\end{pmatrix},
\]
and yields the $C^3$ estimate $\|\mathcal F_{i,k}-\mathcal T_{\Lambda_i}\|_{C^3}\le C_\nu\delta+C_{G,\nu}\eta_k$ on $\mathbb T\times[1/4,2]$. For $\delta$ sufficiently small and $k$ sufficiently large, the second component of $\mathcal F_{+,k}$ maps $\mathbb T\times[3/8,5/4]$ into $(1/4,2)$. Since $T_k$ commutes with $\mathcal R_\pm$, the rescaled full-loop return map is the composition
\[
\mathcal F_k^{\mathrm{in}}:=T_{2k}\circ D_k^{-1}\circ\mathcal R_{\mathrm{in}}\circ D_k
=\mathcal F_{-,k}\circ\mathcal F_{+,k}.
\]
Writing $\mathcal F_k^{\mathrm{in}}-\mathcal T_\Lambda=(\mathcal F_{-,k}-\mathcal T_{\Lambda_-})\circ\mathcal F_{+,k}+(\mathcal T_{\Lambda_-}\circ\mathcal F_{+,k}-\mathcal T_{\Lambda_-}\circ\mathcal T_{\Lambda_+})$ and applying the chain rule on $\mathbb T\times[3/8,5/4]$, we obtain
\begin{equation}\label{eq:island-composition-bound}
	\|\mathcal F_k^{\mathrm{in}}-\mathcal T_\Lambda\|_{C^3}
	\le C_\nu\delta+C_{G,\nu}e^{-k\nu},\qquad
	\Lambda(s):=\sigma_++\sigma_-+\frac{2}{\nu}\log s.
\end{equation}
Here $C_\nu$ depends only on bounds for the derivatives of $\nu^{-1}\log s$ through order four on $[1/4,2]$ and can be chosen uniformly for $\nu$ in any compact subset of $(0,\infty)$, whereas $C_{G,\nu}$ may depend on $G$.

On $\mathcal I=[1/2,1]$, the integrable map $\mathcal T_\Lambda$ has twist $\Lambda'(s)=2/(\nu s)=\Lambda_{\nu/2}'(s)>0$, and $\Lambda(\mathcal I)$ has length $(2\log 2)/\nu$. Choosing an integer $Q_*>\nu/\log 2$, we find $\alpha_{\mathrm{in}}=(j+\xi)/Q_*$ in the middle half of $\Lambda(\mathcal I)$ satisfying~\eqref{eq:alpha-diophantine} with $\gamma_*=\gamma_\xi/Q_*$ and reference action $s_{\alpha_{\mathrm{in}}}\in[2^{-3/4},2^{-1/4}]$. Applying Lemma~\ref{lem:herman-C3} and the parameterization of Section~\ref{subsec:auto-4-3} (with $\nu/2$ and $\sigma_++\sigma_-$ in place of $\nu$ and $\sigma$) yields a threshold $\delta^{\mathrm{island}}_{\nu,\gamma_\xi}>0$, uniform on compact $\nu$-intervals, such that whenever $\delta<\delta^{\mathrm{island}}_{\nu,\gamma_\xi}$ and $k$ is sufficiently large, $\mathcal F_k^{\mathrm{in}}$ admits a smooth invariant graph in $\mathbb T\times\operatorname{int}\mathcal I$ with a smooth regular degree-one parameterization $K_k$ satisfying $\mathcal F_k^{\mathrm{in}}(K_k(t))=K_k(t+\alpha_{\mathrm{in}})$. Setting $\widetilde K_k:=D_k\circ K_k$ and $a_k:=2k-\alpha_{\mathrm{in}}>0$, we have
\begin{equation}\label{eq:R-in-tilde-Kk}
	\mathcal R_{\mathrm{in}}\bigl(\widetilde K_k(t)\bigr)=\widetilde K_k(t-a_k),\qquad
	\widetilde K_k(t+1)=\widetilde K_k(t)+(1,0).
\end{equation}

\paragraph{One turn around the interior fixed point.}
Fix a fixed point $p\in D$ of $\widetilde f_0$, and let $\pi_p:\widetilde{\mathcal U}_p\to\mathbb R^2\setminus\{p\}$ be the universal cover of the punctured plane. Let $\Delta:\widetilde{\mathcal U}_p\to\widetilde{\mathcal U}_p$ denote the deck transformation corresponding to one positive loop around $p$ along the oriented Jordan curve $\partial D=S_+\cup S_-$. Choose a lift $\widehat o_0\in\widetilde{\mathcal U}_p$ of $o_0$; lifting $S_+$ from $\widehat o_0$ determines a lift $\widehat o_1$ of $o_1$, and lifting $S_-$ from $\widehat o_1$ terminates at $\Delta(\widehat o_0)$. Let $\widehat f_0:\widetilde{\mathcal U}_p\to\widetilde{\mathcal U}_p$ be the unique lift of $\widetilde f_0|_{\mathbb R^2\setminus\{p\}}$ fixing $\widehat o_0$, and set $\widehat G:=\widehat f_0^M$. Since $\widetilde f_0$ preserves orientation, fixes $p$, and preserves $S_+$ and $S_-$, the lift $\widehat f_0$ commutes with $\Delta$, preserves the lifted arcs, and fixes $\widehat o_1$.

Lift the local charts $J_{0,\mathrm{in}},J_{0,\mathrm{out}}$ to $\widehat o_0$ and $J_{1,\mathrm{in}},J_{1,\mathrm{out}}$ to $\widehat o_1$. Following~\eqref{eq:iota-def} and~\eqref{eq:iota-out-def}, define the orbit extension maps $\iota_0,\iota_{0,\mathrm{out}}$ at $\widehat o_0$ and $\iota_1,\iota_{1,\mathrm{out}}$ at $\widehat o_1$ by
\[
\iota_0(\theta,w):=\widehat G^m\bigl(J_{0,\mathrm{in}}(\theta-m,w)\bigr)\qquad(\theta-m\in B_{0,\mathrm{in}}),
\]
and analogously for the remaining three charts. As in Section~\ref{subsec:auto-5-1}, these maps are smooth local diffeomorphisms satisfying $\iota_0\circ T_1=\widehat G\circ\iota_0$, $\iota_1\circ T_1=\widehat G\circ\iota_1$, $\iota_{0,\mathrm{out}}\circ\mathcal C_\nu=\iota_0$, and $\iota_{1,\mathrm{out}}\circ\mathcal C_\nu=\iota_1$. In $\widetilde{\mathcal U}_p$, the unnormalized finite-step transitions along the lifted upper and lower arcs are given on their basic rectangles by
\[
\widehat{\mathcal P}_+=J_{1,\mathrm{in}}^{-1}\circ\widehat G^{L_+}\circ J_{0,\mathrm{out}},\qquad
\widehat{\mathcal P}_-=J_{0,\mathrm{in}}^{-1}\circ\Delta^{-1}\circ\widehat G^{L_-}\circ J_{1,\mathrm{out}},
\]
since the lower lifted arc from $\widehat o_1$ enters the incoming chart $\Delta\circ J_{0,\mathrm{in}}$ at $\Delta(\widehat o_0)$. Using $\mathcal R_\pm=T_{-L_\pm}\circ\widehat{\mathcal P}_\pm\circ\mathcal C_\nu$, the argument leading to~\eqref{eq:spatial-equivariance} gives
\begin{equation}\label{eq:island-development}
	\iota_1\circ\mathcal R_+=\iota_0,\qquad
	\iota_0\circ\mathcal R_-=\Delta^{-1}\circ\iota_1,\qquad
	\iota_0\circ\mathcal R_{\mathrm{in}}=\Delta^{-1}\circ\iota_0,\qquad
	\iota_0\circ T_1=\widehat G\circ\iota_0 .
\end{equation}

Define $\ell(t):=\iota_0\bigl(\widetilde K_k(t)\bigr)$ in $\widetilde{\mathcal U}_p$ and its projection $c(t):=\pi_p\bigl(\ell(t)\bigr)\in D\setminus\{p\}$. Combining~\eqref{eq:R-in-tilde-Kk} and~\eqref{eq:island-development} yields
\[
\ell(t+1)=\widehat G\bigl(\ell(t)\bigr),\qquad
\ell(t+a_k)=\Delta\bigl(\ell(t)\bigr).
\]
Thus $c$ is a smooth regular immersion satisfying $c(t+a_k)=c(t)$ and $c(t+1)=\widetilde G\bigl(c(t)\bigr)$, and its lift $\ell$ advances by one deck transformation $\Delta$ over each period $a_k$. {
If $c(t_0)=c(t_0+d)$, equivariance under $\widetilde G$ propagates this equality to $t_0+n$ for every $n\in\mathbb Z$. Since $a_k$ is irrational, reduction modulo $a_k$ and continuity give $c(t+d)=c(t)$ for all $t$. Thus every self-coincidence is a global period. As in Section~\ref{subsec:auto-5-2}, regularity makes the period group equal to $h\mathbb Z$ for some $h>0$, and $a_k=Nh$ for an integer $N\ge1$. A lift of the $h$-period advances by $\Delta^{d_0}$ for an integer $d_0$, so
\[
\ell(t+a_k)=\Delta^{Nd_0}\ell(t)=\Delta\ell(t),\qquad Nd_0=1.
\]
The action of $\Delta$ is free; hence $N=d_0=1$. It follows that $c:\mathbb R/(a_k\mathbb Z)\to D\setminus\{p\}$ is injective. Its image $C_k$ is therefore a smooth Jordan curve of winding number one about $p$, with $\pi_p^{-1}(C_k)=\ell(\mathbb R)$, and $\widetilde G|_{C_k}$ is minimal.
}

If $C_k\cap\widetilde f_0(C_k)\ne\varnothing$, minimality of $\widetilde G|_{C_k}$ and $\widetilde f_0\circ\widetilde G=\widetilde G\circ\widetilde f_0$ imply $\widetilde f_0(C_k)=C_k$. If $C_k$ and $\widetilde f_0(C_k)$ were disjoint, then since both are Jordan curves in $D$ enclosing $p$, one of their bounded interiors would be strictly contained in the other, contradicting area preservation of $\widetilde f_0$. Thus $\widetilde f_0(C_k)=C_k$ and $\widehat f_0\bigl(\ell(\mathbb R)\bigr)=\ell(\mathbb R)$. Writing $\widehat f_0\bigl(\ell(t)\bigr)=\ell\bigl(\mathfrak a(t)\bigr)$, commutation of $\widehat f_0$ with $\widehat G$ and $\Delta$ implies that $\mathfrak a(t)-t$ is both $1$-periodic and $a_k$-periodic, hence constant; since $\widehat f_0^M=\widehat G$, this constant is $1/M$. In the parameterization $\gamma_k(s):=\ell(a_k s)$, we have $\gamma_k(s+1)=\Delta\bigl(\gamma_k(s)\bigr)$ and $\widehat f_0\bigl(\gamma_k(s)\bigr)=\gamma_k\bigl(s+1/(Ma_k)\bigr)$. Consequently, the rotation number of $\widetilde f_0|_{C_k}$ about $p$ and its horizontal rotation number satisfy
\begin{equation}\label{eq:island-inner-frequency}
	\rho_{\mathrm{in}}(\widetilde f_0|_{C_k})=\frac{1}{M(2k-\alpha_{\mathrm{in}})},\qquad
	\rho_x(\widetilde f_0|_{C_k})=0,
\end{equation}
where $\rho_{\mathrm{in}}(\widetilde f_0|_{C_k})$ is of constant type by Lemma~\ref{lem:reciprocal-frequency}. Since $C_k\subset D$ and $\pi|_D$ is an embedding onto a topological disk in $\mathbb A$, the projected circle $\pi(C_k)$ is contractible in $\mathbb A$.

\paragraph{Hausdorff distance and positive area.}
Write the unscaled invariant graph of $\mathcal R_{\mathrm{in}}$ as $\kappa_k(\theta)=\bigl(\theta,W_k(\theta)\bigr)$ with $W_k(\theta)\asymp e^{-k\nu}$, and let $r_k(\theta):=\pi_1\bigl(\mathcal R_{\mathrm{in}}(\kappa_k(\theta))\bigr)$ be the induced strictly increasing lift, which satisfies $r_k(\theta)=\theta-2k+O_G(1)$ by~\eqref{eq:island-composition-bound}. Similarly, the intermediate curve $\kappa_k^+(\theta):=\mathcal R_+\bigl(\kappa_k(\theta)\bigr)$ has vertical coordinate $W_k^+(\theta)\asymp e^{-k\nu}$ by~\eqref{eq:P-taylor} and first coordinate $r_{k,+}(\theta):=\pi_1\bigl(\kappa_k^+(\theta)\bigr)=\theta-k+O_G(1)$. Fix the left endpoint $\theta_{0,\mathrm{in}}^-$ of $B_{0,\mathrm{in}}$ and set $J_k:=\bigl[\theta_{0,\mathrm{in}}^-,r_k^{-1}(\theta_{0,\mathrm{in}}^-)\bigr]$, so that $C_k=\pi_p\bigl(\iota_0(\kappa_k(J_k))\bigr)$. Split $J_k$ at $\theta_{0,\mathrm{in}}^-+k$:
\begin{itemize}[topsep=2pt,itemsep=2pt]
	\item For $\theta\in[\theta_{0,\mathrm{in}}^-,\theta_{0,\mathrm{in}}^-+k]$, we have $\pi_p\bigl(\iota_0(\kappa_k(\theta))\bigr)=\widetilde G^j\bigl(J_{0,\mathrm{in}}(\phi,W_k(\phi))\bigr)$ with $\phi\in[\theta_{0,\mathrm{in}}^-,\theta_{0,\mathrm{in}}^-+1)$ and $0\le j\le k$.
	\item For $\theta\in[\theta_{0,\mathrm{in}}^-+k,r_k^{-1}(\theta_{0,\mathrm{in}}^-)]$, the identity
\[
\iota_0(\kappa_k(\theta))=\iota_1\bigl(\mathcal R_+(\kappa_k(\theta))\bigr)=\iota_1\bigl(\kappa_k^+(\theta)\bigr)
\]
represents $\pi_p\bigl(\iota_0(\kappa_k(\theta))\bigr)$ at $o_1$ as $\widetilde G^j\bigl(J_{1,\mathrm{in}}(\phi_1,W_k^+)\bigr)$ with $\phi_1\in[\theta_{1,\mathrm{in}}^-,\theta_{1,\mathrm{in}}^-+1)$ and $-O_G(1)\le j\le k+O_G(1)$.
\end{itemize}
{
The action coordinate in both representations lies between $c_Ge^{-k\nu}$ and $C_Ge^{-k\nu}$ for fixed positive constants. For $j\ge0$, the box residence-time estimate of Section~\ref{subsec:auto-5-2} applies with these bounds and gives distance at most $C_Ge^{-k\nu/2}$ from the corresponding boundary arcs. The possible negative integers in the second representation satisfy $-J_G\le j<0$ for a fixed $J_G$. On the compact incoming basic strip at $o_1$, each map $\widetilde G^{\,j}\circ J_{1,\mathrm{in}}$, $-J_G\le j<0$, has bounded derivative. Its value at $w=0$ belongs to $\partial D$, so the mean value inequality bounds its distance to $\partial D$ by $C_Gw\le C_Ge^{-k\nu}$. Therefore
\[
\sup_{z\in C_k}\operatorname{dist}(z,\partial D)\le C_Ge^{-k\nu/2}.
\]
} Conversely, since $\partial D=S_+\cup S_-$ is a piecewise smooth Jordan curve with transverse corners at $o_0$ and $o_1$, the construction of Lemma~\ref{lem:corner-collar} provides a Lipschitz retraction $\mathfrak p:U\to\partial D$ on a collar neighborhood $U$ of $\partial D$ with $\mathfrak p\simeq\operatorname{id}_U$. For $k$ large enough that $C_k\subset U\setminus\{p\}$, the restriction $\mathfrak p|_{C_k}:C_k\to\partial D$ has degree $1$ around $p$ and is therefore surjective, which gives $d_H(C_k,\partial D)\le C_G e^{-k\nu/2}$.

Finally, given any interior collar of $\partial D$, choose $k$ sufficiently large that $C_k$ lies in the collar, and take a thin tubular neighborhood $U_k\subset D\setminus\{p\}$ of $C_k$ compactly contained in the collar. Since $\widetilde f_0|_{U_k}$ is smooth, preserves area, and leaves the essential circle $C_k$ of $D\setminus\{p\}$ invariant, any closed curve in $U_k$ isotopic to $C_k$ encloses the same area inside $D$ as its $\widetilde f_0$-image and therefore intersects its image. Applying the local theorem of Fayad and Krikorian \cite[Theorem~1; Theorem~3 and Remark~1]{FK09} to $\widetilde f_0$ on $U_k$ yields a positive-area union of smooth quasiperiodic $\widetilde f_0$-invariant circles in $U_k$ surrounding $p$. This completes the proof of Proposition~\ref{prop:island}.

\begin{proof}[Proof of Theorem~\ref{thm:intro-island}]
	Applying Theorem~\ref{thm:intro-criterion} directly to the upper separatrix $\Sigma_+$ yields the upper exterior invariant circles. For the lower separatrix $\Sigma_-$, conjugating $f_0$ by the orientation-preserving involution $(x,y)\mapsto(-x,-y)$ and using $(-Q,-P)$ in the local saddle charts reverses the direction of horizontal translation $\tau$ while preserving the absolute-value logarithmic phases of Definition~\ref{def:island-geometry} up to constant shifts; applying Theorem~\ref{thm:intro-criterion} to the conjugated map and undoing the conjugation yields the lower exterior invariant circles with negative horizontal rotation numbers. The interior assertions follow from Proposition~\ref{prop:island} once the boundary $C^4$-threshold is taken as the minimum of the exterior and interior thresholds.
	
	It remains to verify the sign of the lower and upper horizontal rotation numbers~\eqref{eq:lower-upper-rotation} when $f_0$ satisfies the positive twist condition $\partial_y(\pi_x\circ\widetilde f_0)\ge b_*>0$. {
Birkhoff's graph theorem \cite{Arn09} applies also to the invariant essential separatrix, giving $\Sigma_+=\operatorname{graph}\psi_+$ with $\psi_+$ Lipschitz. For a point $z=(x_0,y_0)$ in the upper exterior, put $d_z=y_0-\psi_+(x_0)>0$. Corollary~\ref{cor:graph-corollary} gives $\Gamma_k^+=\operatorname{graph}u_k^+$ and
\[
0<u_k^+(x)-\psi_+(x)\le(1+\operatorname{Lip}\psi_+)C_Ge^{-k\nu/2}.
\]
Choose $k$ so large that the right-hand side is smaller than $d_z$. Then $y_0>u_k^+(x_0)$.
} Since $\widetilde f_0$ preserves the region above $\Gamma_k^+$, the lifted orbit $(x_j,y_j)=\widetilde f_0^{\,j}(z)$ satisfies $y_j>u_k^+(x_j)$ for all $j\ge 0$. Let $a_k(x):=\pi_x\bigl(\widetilde f_0(x,u_k^+(x))\bigr)$ be the induced circle lift on $\Gamma_k^+$. Since $a_k$ is conjugate to a positive rotation $\omega_k>0$, it has no fixed points on $\mathbb R$, and $1$-periodicity of $a_k-\operatorname{id}$ gives $c_k:=\min_{x\in\mathbb T}\bigl(a_k(x)-x\bigr)>0$. Monotonicity of $\pi_x\circ\widetilde f_0$ in $y$ then yields
	\[
	x_{j+1}-x_j=\pi_x\bigl(\widetilde f_0(x_j,y_j)\bigr)-x_j
	>\pi_x\bigl(\widetilde f_0(x_j,u_k^+(x_j))\bigr)-x_j
	=a_k(x_j)-x_j\ge c_k>0\qquad(j\ge 0),
	\]
	whence $\rho_x^-(z)\ge c_k>0$. The estimate $\rho_x^+(z)<0$ in the lower exterior follows symmetrically from the lower invariant graphs.
\end{proof}

\section{Application to the Herman hyperbolic family}\label{sec:family}
	
	In this section we prove Theorems~\ref{thm:herman-hyperbolic-family} and~\ref{cor:herman-window} for the Herman hyperbolic family $f_\lambda$ of Section~\ref{subsec:auto-1-2}, and then deduce Theorem~\ref{thm:intro-herman}. The vertical rescaling $y=\lambda r$ conjugates $f_\lambda$ to $F_\lambda$. We take the iterate $G_\lambda=F_\lambda^{M_\lambda}$, where $M_\lambda$ and $\nu_\lambda$ are defined in~\eqref{eq:thm1.1-scale-parameters}. We verify conditions (G1)--(G2) for $G_\lambda$ and the uniform estimate
	\[
	\frac12\le\nu_\lambda\le1,\qquad
	\|\zeta_\lambda\|_{C^4(\mathbb T)}\le C\lambda,
	\]
	and then use Theorem~\ref{thm:boundary-transition-criterion} and Corollaries~\ref{cor:f0-invariant} and~\ref{cor:graph-corollary} to recover the conclusions for $f_\lambda$. Unless otherwise stated, the constants $C,C_j$ in this section are independent of $\lambda\in(0,\lambda_0]$, whereas the initial index $k_\lambda$ and the distance constant $C_\lambda$ may depend on $\lambda$.
	
	\subsection{The separatrix and linearizing coordinates}\label{subsec:geometry}
	
	Make the vertical rescaling $y = \lambda r$, and write
	\begin{equation}\label{eq:W-Phi-def}
		W_\lambda(x) := \frac{\psi_\lambda(x)}{\lambda}, \qquad
		\Phi_\lambda(x) := \frac{\varphi_{\lambda^2}(x)}{\lambda^2}.
	\end{equation}
	By~\eqref{eq:C5-interval-norm},~\eqref{eq:positive-drift}, and~\eqref{eq:phi-C4-bound}, we have $\|W_\lambda\|_{C^5([0, 1])} \le C$, $W_\lambda$ has a uniform positive lower bound on every compact subset of $(0, 1)$, and $\Phi_\lambda \in C^\infty(\mathbb{T})$ satisfies $\|\Phi_\lambda\|_{C^j(\mathbb{T})} \le C_j$ for $0 \le j \le 4$. In the coordinates $(x,r)$, $f_\lambda$ is conjugate to the cylinder map $F_\lambda$ with chosen lift
	\begin{equation}\label{eq:F-lambda-def}
		\widetilde F_\lambda(x, r) := \bigl(x + \lambda r,\; r + \lambda\Phi_\lambda(x + \lambda r)\bigr).
	\end{equation}
	A direct calculation gives $\det D\widetilde F_\lambda=1$. Since $\widetilde F_\lambda$ is the composition of the two area-preserving shear diffeomorphisms $(x, r) \mapsto (x + \lambda r, r)$ and $(x, r) \mapsto \bigl(x, r + \lambda\Phi_\lambda(x)\bigr)$, for $0\le t\le1$ the maps
	\[
	(x,r)\longmapsto
	\bigl(x+t\lambda r,\ r+t\lambda\Phi_\lambda(x+t\lambda r)\bigr)
	\]
	give a smooth area-preserving isotopy from the identity to $F_\lambda$. Hence $F_\lambda$ is a smooth area-preserving diffeomorphism isotopic to the identity and preserves the area form $\mathrm{d}x \wedge \mathrm{d}r$.
	
	By definition, $\varphi_{\lambda^2}\bigl(g_\lambda(x)\bigr)=\psi_\lambda\bigl(g_\lambda(x)\bigr)-\psi_\lambda(x)$. Dividing by $\lambda$ gives $\lambda\Phi_\lambda\bigl(g_\lambda(x)\bigr) = W_\lambda\bigl(g_\lambda(x)\bigr) - W_\lambda(x)$. Substitution into~\eqref{eq:F-lambda-def} yields invariance of the graph $\operatorname{graph} W_\lambda$:
	\[
	\widetilde F_\lambda\bigl(x,W_\lambda(x)\bigr) = \bigl(g_\lambda(x),\, W_\lambda(g_\lambda(x))\bigr).
	\]
	Since $g_\lambda(x) > x$ throughout $(0, 1)$ and there are no fixed points there, every interior forward orbit on the arc $\widetilde{\Sigma}_\lambda^{\mathrm{sc}} := \{(x, W_\lambda(x)) : x \in [0, 1]\}$ joining the saddles $o = (0, 0)$ and $\tau(o) = (1, 0)$ in the universal cover $\mathbb{R}^2$ converges monotonically to $\tau(o)$, and every interior backward orbit converges monotonically to $o$. Thus $\widetilde{\Sigma}_\lambda^{\mathrm{sc}}$ with its endpoints removed is both an unstable branch of $o$ and a stable branch of $\tau(o)$. Its projection $\Sigma_\lambda^{\mathrm{sc}} = \operatorname{graph} W_\lambda$ is smooth away from the saddle and satisfies (G1). The zero-flux fact in Section~\ref{subsec:auto-1-1} also shows that $F_\lambda$ is an exact symplectic diffeomorphism.

	We next construct common linear symplectic charts at the integer saddles and verify (G2). Write $c_\lambda := (1 + \lambda)^{-1}$ and $d_\lambda := \sqrt{1 + c_\lambda}$. Near the saddle $\tau^n(o) = (n, 0)$, $n \in \mathbb{Z}$, define the linear coordinate transformation $\Phi_{\tau^n(o)} = \Phi_o \circ \tau^{-n}$ by
	\begin{equation}\label{eq:linear-chart-7-7}
		\begin{gathered}
			x - n = \frac{Q - P}{d_\lambda}, \qquad
			r = \frac{Q + c_\lambda P}{d_\lambda}\\
			\Longleftrightarrow\\
			Q = \frac{c_\lambda(x - n) + r}{d_\lambda}, \qquad
			P = \frac{-(x - n) + r}{d_\lambda}.
		\end{gathered}
	\end{equation}
	Since $\det D\Phi_o = (1 + c_\lambda)/d_\lambda^2 = 1$, we have $\mathrm{d}x \wedge \mathrm{d}r = \mathrm{d}Q \wedge \mathrm{d}P$, and the matrix of this linear transformation and its inverse are uniformly bounded in $\lambda \in (0, \lambda_0]$. If $|x - n| < r_0/4$ and $|\lambda r| < r_0/4$, where $\lambda \le 1$, then $|(x - n) + \lambda r| < r_0/2 \le c_\lambda r_0$. {In this neighborhood,~\eqref{eq:endpoint-affine-structure} gives $\Phi_\lambda(n+t)=c_\lambda t$ for $|t|<c_\lambda r_0$. Applying this with $t=(x-n)+\lambda r$, we see that $\widetilde F_\lambda$ is represented in the local coordinates $(x-n,r)$ by the linear matrix}
	\[
	A_\lambda =
	\begin{pmatrix}
		1 & \lambda \\
		c_\lambda\lambda & 1 + c_\lambda\lambda^2
	\end{pmatrix}.
	\]
	The eigenvectors of $A_\lambda$ associated with $1 + \lambda$ and $c_\lambda$ are $(1, 1)^T$ and $(-1, c_\lambda)^T$, respectively. Thus~\eqref{eq:linear-chart-7-7} diagonalizes the local map as
	\[
	(\Phi_o\circ\widetilde F_\lambda\circ\Phi_o^{-1})(Q,P) = \bigl((1 + \lambda)Q,\; c_\lambda P\bigr).
	\]
	On the two sides of the saddle, the separatrix $\Sigma_\lambda^{\mathrm{sc}}$ satisfies $r = x - n$ for $0 < x - n < r_0/2$, corresponding to $P = 0, Q = d_\lambda(x - n) > 0$, and $r = -c_\lambda(x - n)$ for $-r_0/2 < x - n < 0$, corresponding to $Q = 0, P = -d_\lambda(x - n) > 0$.

	For $M_\lambda = \lfloor \lambda^{-1} \rfloor$ and $\nu_\lambda = M_\lambda \log(1 + \lambda)$, define $G_\lambda := F_\lambda^{M_\lambda}$ and $\widetilde G_\lambda := \widetilde F_\lambda^{M_\lambda}$. The map $G_\lambda$ still has $\Sigma_\lambda^{\mathrm{sc}}$ as a non-splitting separatrix and satisfies \textbf{(G1)}. In this paragraph, the coordinate formulas are evaluated using the globally defined affine transformations in~\eqref{eq:linear-chart-7-7}; their restrictions to $U(\tau^n(o))$ are the local charts appearing in \textbf{(G2)}. By~\eqref{eq:linear-chart-7-7} and $0 < c_\lambda < 1 < d_\lambda$, every $(Q, P) \in \mathbb{R}^2$ satisfies $|x - n| \le |Q| + |P|$ and $|\lambda r| \le |r| \le |Q| + |P|$. Fix $R \in (0, r_0/8)$ independent of $\lambda$, so that the symmetric open square $\mathcal{N} := (-R, R)^2$ satisfies $\Phi_{\tau^n(o)}^{-1}(\mathcal{N}) \subset \{\vert{}x - n\vert{} < r_0/4,\, \vert{}\lambda r\vert{} < r_0/4\}$, on which the one-step linearization holds. Define
	\[
	V := \left(-\frac{R}{e},\, \frac{R}{e}\right)^2, \qquad
	U(o) := \Phi_o^{-1}(V), \qquad
	0 < Q_{\mathrm{box}}, P_{\mathrm{box}} < \frac{R}{2e^2},
	\]
	with $Q_{\mathrm{box}}, P_{\mathrm{box}}$ independent of $\lambda$, and set $\mathcal{B} := [0, Q_{\mathrm{box}}] \times [0, P_{\mathrm{box}}] \Subset V$ and $U_{\mathrm{box}}(o) := \Phi_o^{-1}(\mathcal{B}) \Subset U(o)$. For $0 \le j \le M_\lambda$, we have $1 \le (1 + \lambda)^j \le e^{M_\lambda \lambda} \le e$ and $e^{-1} \le c_\lambda^j \le 1$. For any $(Q, P) \in V$, both the forward candidate points $\bigl((1 + \lambda)^j Q, c_\lambda^j P\bigr)$ and the backward candidate points $\bigl(c_\lambda^j Q, (1 + \lambda)^j P\bigr)$ remain in $\mathcal{N}$ for all $0 \le j \le M_\lambda$. Induction on $j$ yields, throughout $V$,
	\begin{equation}\label{eq:G-lambda-linear}
		(\Phi_o \circ \widetilde G_\lambda^{\,\pm 1} \circ \Phi_o^{-1})(Q, P)
		= \bigl(e^{\pm \nu_\lambda} Q,\; e^{\mp \nu_\lambda} P\bigr), \qquad
		\nu_\lambda \longrightarrow 1 \quad (\lambda \to 0^+).
	\end{equation}
	For $(Q, P) \in \mathcal{B}$, the bound $e^{\nu_\lambda} \le e$ gives $\bigl(e^{\pm \nu_\lambda} Q, e^{\mp \nu_\lambda} P\bigr) \in [0, e Q_{\mathrm{box}}] \times [0, e P_{\mathrm{box}}] \subset [-R/(2e), R/(2e)]^2 \Subset V$, so $\widetilde G_\lambda^{\,\pm 1}\bigl(U_{\mathrm{box}}(o)\bigr) \Subset U(o)$. Consequently, $G_\lambda$ satisfies \textbf{(G2)} with $V$ and $\mathcal{B}$ independent of $\lambda$.

	Moreover, for all sufficiently small $\lambda \in (0, \lambda_0]$,
	\begin{equation}\label{eq:nu-M-bounds}
		\frac{1}{2} \le \nu_\lambda \le 1, \qquad
		M_\lambda \lambda \ge \frac{1}{2}.
	\end{equation}
	
	\subsection{Uniform estimates for the finite-step transition}\label{subsec:auto-6-3}
	
	Fix $\epsilon_0\in(0,1)$, and choose $\theta_-<\theta_+$ independent of $\lambda$ such that
	\[
	\theta_+-\theta_-=1+\epsilon_0,\qquad
	\theta_++1<0,\qquad e^{(\theta_++1)/2}<Q_{\mathrm{box}},
	\]
	Write $B:=[\theta_-,\theta_+]$. By~\eqref{eq:nu-M-bounds}, we have $0 < e^{\nu_\lambda \theta} \le e^{(\theta_+ + 1)/2} < Q_{\mathrm{box}}$ for every $\theta\in B\cup(B+1)$. On the outgoing half-axis $P = 0$, equation~\eqref{eq:linear-chart-7-7} shows that the logarithmic phase coordinate $\theta = \nu_\lambda^{-1}\log Q$ corresponds to the horizontal coordinate $x_\lambda(\theta) = e^{\nu_\lambda \theta}/d_\lambda$. Using $\theta \le \theta_+ < 0$, $\nu_\lambda \in [1/2, 1]$, and $d_\lambda \in [\sqrt{3/2}, \sqrt{2}]$, set
	\[
	a := \frac{e^{\theta_-}}{\sqrt{2}}, \qquad b := \frac{e^{\theta_+/2}}{\sqrt{3/2}},
	\]
	Then $0 < a \le x_\lambda(\theta) \le b < r_0/4$ uniformly for every $\theta \in B$.
	
	Choose $r_{\mathrm{in}} \in (0, r_0)$ independent of $\lambda$ so that the separatrix points satisfying $1 - r_{\mathrm{in}} \le x < 1$ all lie in $\Phi_{\tau(o)}^{-1}\bigl(\{0\} \times (0, P_{\mathrm{box}})\bigr)$. On the compact interval $[a, 1 - r_{\mathrm{in}}] \Subset (0, 1)$, equation~\eqref{eq:positive-drift} gives $c > 0$ such that
	\[
	g_\lambda(x) - x = \lambda W_\lambda(x) \ge c\lambda.
	\]
	The number of $g_\lambda$-iterates required for an orbit starting in $[a, b]$ to enter $[1 - r_{\mathrm{in}}, 1)$ is at most
	\[
	1 + \left\lceil \frac{1 - r_{\mathrm{in}} - a}{c\lambda} \right\rceil \le C_0 \lambda^{-1} \le 2 C_0 M_\lambda.
	\]
	{Fix a positive integer $L\ge2C_0$. The preceding bound shows that $\widetilde G_\lambda^L J_{\mathrm{out}}(B,0)$ lies on the incoming stable-axis segment inside $U_{\mathrm{box}}(\tau(o))$. The same holds for $B+1$, because
\[
J_{\mathrm{out}}(\theta+1,0)=\widetilde G_\lambda J_{\mathrm{out}}(\theta,0)
\]
and that stable-axis segment is forward invariant under $\widetilde G_\lambda$. Thus the same integer $L$, independent of $\lambda$, works on $B\cup(B+1)$.}
	
	On the other hand, $\|W_\lambda'\|_\infty \le C$ and $W_\lambda(1) = 0$ imply $W_\lambda(x) \le C(1 - x)$. When $C\lambda \le 1/2$, using $\log(1 - t) \ge -2t$ for $0 \le t \le 1/2$ gives, for any $x \in [a, b]$ and $0 \le j \le L M_\lambda$,
	\[
	1 - g_\lambda^j(x) \ge (1 - C\lambda)^j(1 - x) \ge e^{-2Cj\lambda}(1 - x) \ge (1 - b)e^{-2CL} =: d > 0.
	\]
	Thus all intermediate points $g_\lambda^j(x)$, $0 \le j \le L M_\lambda$, lie in $[a, 1 - d]$, and the endpoint $g_\lambda^{L M_\lambda}(x)$ lies in $[1 - r_{\mathrm{in}}, 1 - d]$.
	
\begin{figure}[tbp]
\centering
\begingroup
\setlength{\unitlength}{1pt}
\normalsize\color{black}
\scalebox{0.88}{%
\begin{picture}(450,435)
\linethickness{0.45pt}
\put(10,420){\makebox(0,0)[l]{(a) Uniform passage along the separatrix}}
\put(225,402){\makebox(0,0)[c]{$x_\lambda(\theta)=e^{\nu_\lambda\theta}/d_\lambda\in[a,b]\quad(\theta\in B)$}}
\put(225,385){\makebox(0,0)[c]{$\widetilde G_\lambda^L=\widetilde F_\lambda^{LM_\lambda},\qquad L\ge2C_0\ \text{fixed}$}}
\put(24,300){\vector(1,0){414}}
\put(442,300){\makebox(0,0)[l]{$x$}}
\linethickness{0.85pt}
\qbezier(35,300)(110,325)(175,352)
\qbezier(175,352)(225,377)(275,352)
\qbezier(275,352)(340,325)(415,300)
\linethickness{0.45pt}
\put(35,300){\circle*{3}}
\put(415,300){\circle*{3}}
\put(35,288){\makebox(0,0)[c]{$0$}}
\put(415,288){\makebox(0,0)[c]{$1$}}
\put(22,310){\makebox(0,0)[c]{$o$}}
\put(430,310){\makebox(0,0)[c]{$\tau(o)$}}
\qbezier[150](74,297)(74,313.5)(74,330)
\put(74,282){\makebox(0,0)[c]{$a$}}
\qbezier[150](115,297)(115,313.5)(115,330)
\put(115,282){\makebox(0,0)[c]{$b$}}
\qbezier[150](344,297)(344,313.5)(344,330)
\put(344,282){\makebox(0,0)[c]{$1-r_{\mathrm{in}}$}}
\qbezier[150](390,297)(390,313.5)(390,330)
\put(390,282){\makebox(0,0)[c]{$1-d$}}
\qbezier[187](74,300)(94.5,300)(115,300)
\qbezier[110](115,300)(115,312)(115,324)
\qbezier[187](115,324)(94.5,324)(74,324)
\qbezier[110](74,324)(74,312)(74,300)
\qbezier[187](74,306)(94.5,306)(115,306)
\qbezier[187](74,312)(94.5,312)(115,312)
\qbezier[187](74,318)(94.5,318)(115,318)
\qbezier[210](344,300)(367,300)(390,300)
\qbezier[110](390,300)(390,312)(390,324)
\qbezier[210](390,324)(367,324)(344,324)
\qbezier[110](344,324)(344,312)(344,300)
\qbezier[210](344,306)(367,306)(390,306)
\qbezier[210](344,312)(367,312)(390,312)
\qbezier[210](344,318)(367,318)(390,318)
\put(145,317){\vector(1,0){36}}
\put(163,330){\makebox(0,0)[c]{$\ge c\lambda$}}
\put(236,335){\makebox(0,0)[c]{$r=W_\lambda(x)$}}
\qbezier[114](390,338)(402.5,338)(415,338)
\qbezier[37](390,334)(390,338)(390,342)
\qbezier[37](415,334)(415,338)(415,342)
\put(403,350){\makebox(0,0)[c]{$d>0$}}
\put(225,263){\makebox(0,0)[c]{$g_\lambda^{LM_\lambda}([a,b])\subset[1-r_{\mathrm{in}},1-d],\quad d=(1-b)e^{-2CL}$}}
\put(10,243){\makebox(0,0)[l]{(b) Fixed intervals and buffer distances}}
\qbezier[1728](35,220)(225,220)(415,220)
\qbezier[37](35,216)(35,220)(35,224)
\qbezier[37](415,216)(415,220)(415,224)
\put(225,230){\makebox(0,0)[c]{$J\Subset(0,1)$}}
\qbezier[1546](55,195)(225,195)(395,195)
\qbezier[37](55,191)(55,195)(55,199)
\qbezier[37](395,191)(395,195)(395,199)
\put(225,205){\makebox(0,0)[c]{$J_0\Subset J$}}
\qbezier[1346](74,170)(222,170)(370,170)
\qbezier[37](74,166)(74,170)(74,174)
\qbezier[37](370,166)(370,170)(370,174)
\put(225,180){\makebox(0,0)[c]{$[a,1-d]$}}
\qbezier[87](55,170)(64.5,170)(74,170)
\put(64.5,181){\makebox(0,0)[c]{$\ge\mu$}}
\qbezier[114](370,170)(382.5,170)(395,170)
\put(382.5,181){\makebox(0,0)[c]{$\ge\mu$}}
\linethickness{1.2pt}
\qbezier[296](62,149)(94.5,149)(127,149)
\linethickness{0.45pt}
\qbezier[37](62,145)(62,149)(62,153)
\qbezier[37](127,145)(127,149)(127,153)
\put(145,149){\makebox(0,0)[l]{$I_{\mathrm{out}}=[a-\varepsilon,b+\varepsilon]$}}
\put(225,130){\makebox(0,0)[c]{$\varepsilon e^{CL}<\mu/2,\qquad\mu=\operatorname{dist}([a,1-d],\partial J_0)$}}
\put(10,109){\makebox(0,0)[l]{(c) Iterated intervals and Taylor segments}}
\qbezier[855](20,66)(114,66)(208,66)
\qbezier[37](20,62)(20,66)(20,70)
\qbezier[37](208,62)(208,66)(208,70)
\put(204,78){\makebox(0,0)[c]{$J_0$}}
\qbezier[423](62,62)(108.5,62)(155,62)
\qbezier[37](155,62)(155,66)(155,70)
\qbezier[423](155,70)(108.5,70)(62,70)
\qbezier[37](62,70)(62,66)(62,62)
\linethickness{1.2pt}
\qbezier[287](77,66)(108.5,66)(140,66)
\linethickness{0.45pt}
\put(106,94){\makebox(0,0)[c]{$g_\lambda^j(I_{\mathrm{out}})$}}
\put(106,80){\makebox(0,0)[c]{$g_\lambda^j([a,b])$}}
\qbezier[69](62,57)(69.5,57)(77,57)
\qbezier[28](62,54)(62,57)(62,60)
\qbezier[28](77,54)(77,57)(77,60)
\qbezier[69](140,57)(147.5,57)(155,57)
\qbezier[28](140,54)(140,57)(140,60)
\qbezier[28](155,54)(155,57)(155,60)
\put(110,44){\makebox(0,0)[c]{Each gap $\le\varepsilon e^{CL}<\mu/2$}}
\put(110,29){\makebox(0,0)[c]{$0\le j\le LM_\lambda$}}
\qbezier[869](247,66)(342.5,66)(438,66)
\qbezier[37](247,62)(247,66)(247,70)
\qbezier[37](438,62)(438,66)(438,70)
\put(434,78){\makebox(0,0)[c]{$J$}}
\linethickness{0.85pt}
\qbezier[646](270,66)(341,66)(412,66)
\linethickness{1.5pt}
\qbezier[237](330,66)(356,66)(382,66)
\linethickness{0.45pt}
\put(330,66){\circle*{3}}
\put(382,66){\circle*{3}}
\put(286,78){\makebox(0,0)[c]{$J_0$}}
\put(330,52){\makebox(0,0)[c]{$x$}}
\put(390,52){\makebox(0,0)[c]{$g_\lambda(x)$}}
\put(344,94){\makebox(0,0)[c]{$x+s\lambda W_\lambda(x)\in J$}}
\put(341,29){\makebox(0,0)[c]{$x\in J_0,\quad0\le s\le1$}}
\put(225,9){\makebox(0,0)[c]{$\lambda\ \text{fixed}:\quad0\le v<v_*(\lambda)$}}
\end{picture}%
}
\endgroup

\caption{Uniform transition estimates. (a) The initial and terminal horizontal intervals are hatched; the transition time $L$ and the terminal gap $d>0$ are independent of $\lambda$. (b) The fixed nested intervals leave buffer distance $\mu$. (c) The outlined interval is $g_\lambda^j(I_{\mathrm{out}})$; its thick central segment is $g_\lambda^j([a,b])$. The marked end gaps are bounded by $\varepsilon e^{CL}$, and the Taylor segment remains in $J$. }
\label{fig:uniform-finite-transition}
\end{figure}
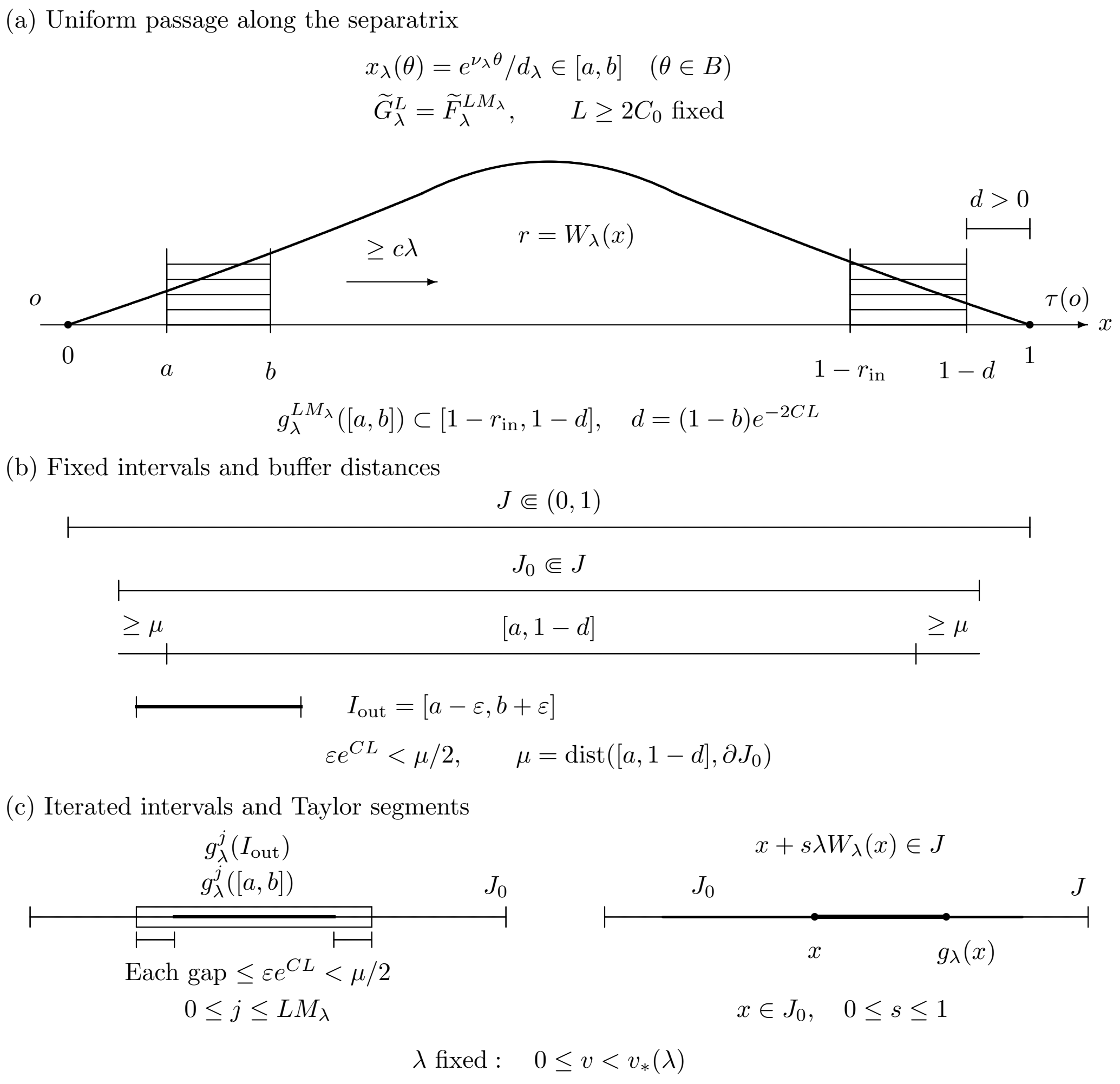

	Fix compact subintervals
	\[
	[a, 1 - d] \Subset J_0 \Subset J \Subset (0, 1), \qquad
	\mu := \operatorname{dist}\bigl([a, 1 - d],\, \partial J_0\bigr) > 0,
	\]
	and take a closed interval $I_{\mathrm{out}} := [a - \varepsilon, b + \varepsilon] \Subset (0,1)$, where $0 < \varepsilon < \frac{1}{2}\mu e^{-CL}$. If necessary, shrink $\varepsilon$ further so that its endpoints remain in $(0,1)$. The mean value theorem and $\|(g_\lambda^j)'\|_\infty \le (1 + C\lambda)^j \le e^{CL}$ for $0 \le j \le L M_\lambda$ give $\sup_{x \in I_{\mathrm{out}}} \operatorname{dist}\bigl(g_\lambda^j(x), g_\lambda^j([a, b])\bigr) < \mu/2$. Hence every iterate of $I_{\mathrm{out}}$ lies in $\operatorname{int} J_0$. Decrease $\lambda_0$ further so that $\|g_\lambda - \operatorname{id}\|_\infty < \frac{1}{2}\operatorname{dist}(J_0, \partial J)$. Then $x + s\lambda W_\lambda(x) \in J$ for every $x \in J_0$ and $0 \le s \le 1$. For each fixed $\lambda \in (0, \lambda_0]$, Lemma~\ref{lem:finite_transition} ensures that the transition map induced by $\widetilde G_\lambda^L$ is well defined on some one-sided collar $v \in [0, v_*(\lambda))$. Figure~\ref{fig:uniform-finite-transition} summarizes the uniform transition time along the separatrix, the uniform separation of the orbit from the endpoints, and the fixed intervals used in the subsequent estimates.
	
	\subsection{Estimates for the boundary transition map}\label{subsec:auto-6-4}
	
	Fix a base point $x_* \in J$ in the compact interval $J$, and define
	\begin{equation}\label{eq:h-lambda-def}
		h_\lambda(x) := \int_{x_*}^x \frac{\mathrm{d}u}{W_\lambda(u)} \qquad (x \in J).
	\end{equation}
	We call $t=h_\lambda(x)$ an \emph{integral coordinate} along this portion of the separatrix: it is obtained by integrating the reciprocal of the normalized displacement $W_\lambda(x)=(g_\lambda(x)-x)/\lambda$. The identity $h_\lambda'W_\lambda=1$ makes the first-order displacement in this coordinate equal to the constant $\lambda$. The conjugated discrete map retains a quadratic remainder, as shown in~\eqref{eq:H-lambda-expansion}.
	Since $0 < c_J \le W_\lambda(x) \le C_J$ on $J$, we have $1/C_J \le h_\lambda'(x) \le 1/c_J$. Differentiating $W_\lambda \cdot (1/W_\lambda) = 1$ a total of $r$ times, for $1 \le r \le 5$, gives
	\[
	\left(\frac{1}{W_\lambda}\right)^{(r)}
	= -\frac{1}{W_\lambda} \sum_{i=1}^r \binom{r}{i} W_\lambda^{(i)} \left(\frac{1}{W_\lambda}\right)^{(r-i)}.
	\]
	Using $\|W_\lambda\|_{C^5(J)} \le C$, induction gives $\|1/W_\lambda\|_{C^5(J)} \le C$, and hence $\|h_\lambda\|_{C^6(J)} \le C$. The inverse differentiation formulas and $h_\lambda' \ge 1/C_J$ show that the inverse map $h_\lambda^{-1}$ also has a uniform $C^6$ norm bound on $h_\lambda(J)$.
	
	Write $x = h_\lambda^{-1}(t)$. Substituting $h_\lambda'(x)W_\lambda(x) = 1$ into the second-order Taylor expansion of $h_\lambda\bigl(x + \lambda W_\lambda(x)\bigr)$ gives the following expansion of the conjugate map $\mathcal{H}_\lambda := h_\lambda \circ g_\lambda \circ h_\lambda^{-1}$:
	\begin{equation}\label{eq:H-lambda-expansion}
		\begin{aligned}
			\mathcal{H}_\lambda(t) &= t + \lambda + \lambda^2 E_\lambda(t),\\
			E_\lambda(t) &:= W_\lambda(x)^2 \int_0^1 (1 - s)h_\lambda''\bigl(x + s\lambda W_\lambda(x)\bigr)\,\mathrm{d}s, \qquad
			\|E_\lambda\|_{C^4(h_\lambda(J_0))} \le C,
		\end{aligned}
	\end{equation}
	Here the $C^4$ bound follows directly from the Leibniz rule,~\eqref{eq:faa-di-bruno}, and $\|h_\lambda\|_{C^6(J)}, \|W_\lambda\|_{C^4(J)} \le C$.
	
	By~\eqref{eq:endpoint-affine-structure}, we have $W_\lambda(x) = x$ in the outgoing affine region at the left endpoint and $W_\lambda(x) = c_\lambda(1 - x)$ in the incoming affine region at the right endpoint. Integrating~\eqref{eq:h-lambda-def} on these two subintervals of $J$ gives
	\begin{equation}\label{eq:h-lambda-endpoints}
		h_\lambda(x) = \log x + A_\lambda^-, \qquad
		h_\lambda(x) = -(1 + \lambda)\log(1 - x) + A_\lambda^+,
	\end{equation}
    where $A_\lambda^-,A_\lambda^+$ are constants.
    On the outgoing and incoming branches, respectively, \eqref{eq:linear-chart-7-7} gives
    $Q=d_\lambda x$ and $P=d_\lambda(1-x)$. Thus \eqref{eq:h-lambda-endpoints}
    is precisely the endpoint matching in Lemma~\ref{lem:one-step-transition}, with
    \[
       a_-=\nu_\lambda,\qquad a_+=(1+\lambda)\nu_\lambda,\qquad
       \frac{a_-^4}{a_+}=\frac{\nu_\lambda^3}{1+\lambda}\le1
    \]
    by \eqref{eq:nu-M-bounds}.

    Apply the lemma to $F_\lambda$ and $G_\lambda=F_\lambda^{M_\lambda}$.
    In the coordinate $x$ along the lifted separatrix the restriction of $F_\lambda$ is $g_\lambda$.
    Section~\ref{subsec:auto-6-3} gives a common transition time $L$ independent of $\lambda$:
    the interval $[\theta_-,\theta_-+1]\subset B$ covers one full phase period, and all intermediate iterates lie in $J_0$ in the $x$ coordinate.
    By \eqref{eq:H-lambda-expansion} the translation step in the $h_\lambda$ coordinate is $d=\lambda$, with defect $e=\lambda^2 E_\lambda$.
    We may therefore take
    \[
       n=LM_\lambda,\qquad
       \varepsilon_\lambda:=\lambda^2\|E_\lambda\|_{C^4(h_\lambda(J_0))}\le C\lambda^2,
       \qquad n\varepsilon_\lambda\le CL\lambda .
    \]
    In particular $n\varepsilon_\lambda\le1$ for all sufficiently small $\lambda$.
    If $b_\lambda$ denotes the boundary transition map of $G_\lambda$ and $\zeta_\lambda$ the zero-mean part of its displacement, then
    \eqref{eq:transition-distance-translations} and \eqref{eq:one-step-transition} yield
    \begin{equation}\label{eq:zeta-lambda-C4}
       \|\zeta_\lambda\|_{C^4(\mathbb T)}
       =\inf_{c\in\mathbb R}\|b_\lambda-(\operatorname{id}+c)\|_{C^4(\mathbb T)}
       \le C_4\frac{\nu_\lambda^3}{1+\lambda}\,n\varepsilon_\lambda
       \le C\lambda .
    \end{equation}

	\subsection{Proof of Theorem~\ref{thm:herman-hyperbolic-family}}\label{subsec:auto-6-5}
	
	\paragraph{Existence and rotation number.}
	Fix the constant-type data $(\xi, \gamma_\xi)$. Since $\nu_\lambda\in[1/2,1]$, the positive integer $Q_*$ and threshold $\delta_* := \inf_{\nu\in[1/2,1]}\delta_{\nu,\gamma_\xi} > 0$ in Theorem~\ref{thm:boundary-transition-criterion} are independent of $\lambda$. Choose $\lambda_*\in(0,\lambda_0]$ so that~\eqref{eq:zeta-lambda-C4} satisfies $C\lambda_*<\delta_*$. For each $0<\lambda<\lambda_*$, applying Theorem~\ref{thm:boundary-transition-criterion} to $G_\lambda$ gives $\alpha_\lambda=(j_\lambda+\xi)/Q_*$ and an integer $k_\lambda$ such that, for every $k\ge k_\lambda$, $G_\lambda$ has a smooth essential invariant circle $\Gamma_{\lambda,k}^{\mathrm{sc}}$ with rotation number $1/(k-\alpha_\lambda)$ relative to $\widetilde G_\lambda$.
	
	By Corollary~\ref{cor:f0-invariant} with $f_0=F_\lambda, M=M_\lambda$, $\Gamma_{\lambda,k}^{\mathrm{sc}}$ is invariant under $F_\lambda$, and its rotation number $\omega_{\lambda,k}=1/[M_\lambda(k-\alpha_\lambda)]$ relative to $\widetilde F_\lambda$ satisfies~\eqref{eq:thm1.1-rotation-number}--\eqref{eq:thm1.1-diophantine}. Moreover, $\partial_r(\pi_x\circ\widetilde F_\lambda)=\lambda>0$, so Corollary~\ref{cor:graph-corollary}(1) shows that $\Gamma_{\lambda,k}^{\mathrm{sc}}$ is a smooth graph. Transforming back by $y=\lambda r$ gives a smooth invariant graph $\Gamma_{\lambda,k}=\operatorname{graph}u_{\lambda,k}$ of $f_\lambda$, on which the induced circle map is smoothly conjugate to $R_{\omega_{\lambda,k}}$.
	
	\paragraph{Accumulation on the separatrix.}
	In the first quadrant $\{Q > 0, P > 0\}$ of the local symplectic chart \eqref{eq:linear-chart-7-7}, the local expressions for the separatrix on the two sides of the saddle, $W_\lambda(x) = x - n$ for $x \ge n$ and $W_\lambda(x) = -c_\lambda(x - n)$ for $x \le n$, give
	\[
	r - W_\lambda(x) =
	\begin{cases}
		d_\lambda P > 0, & x \ge n,\\[4pt]
		d_\lambda Q > 0, & x \le n,
	\end{cases}
	\]
	Thus $\{Q > 0, P > 0\}$ lies above $\Sigma_\lambda^{\mathrm{sc}} = \operatorname{graph} W_\lambda$. Since $\operatorname{Lip}(W_\lambda) \le C$, applying Corollary~\ref{cor:graph-corollary}(2) to $F_\lambda$ and multiplying by the scaling factor $\lambda > 0$ immediately yields~\eqref{eq:thm1.1-u-above-psi} and~\eqref{eq:thm1.1-decay-rate}:
	\[
	0 < u_{\lambda, k}(x) - \psi_\lambda(x) \le C_\lambda e^{-k\nu_\lambda/2} \qquad (\forall x \in \mathbb{T}).
	\]
	
	\paragraph{Uniform slope bound.}
	Estimate~\eqref{eq:thm1.1-slope-bound} follows from the general derivative estimate below.
	
	\begin{Lemma}\label{lem:sqrt-slope}
		Let $\varphi\in C^1(\mathbb T)$ and $\widetilde f_\varphi(x,y)=(x+y,y+\varphi(x+y))$. If $u\in C^1(\mathbb T)$ defines an invariant graph $\Gamma=\operatorname{graph}u$ of $f_\varphi$, and we set $d := \|\varphi'\|_\infty$, then
		\begin{equation}\label{eq:sqrt-slope-bound}
			-\frac{\sqrt{d^2 + 4d} - d}{2} \le u'(x) \le \frac{\sqrt{d^2 + 4d} + d}{2} \qquad (\forall x \in \mathbb{T}).
		\end{equation}
	\end{Lemma}
	
	\begin{proof}
		The maps $f_\varphi$ and $f_\varphi^{-1}$ induce mutually inverse orientation-preserving circle diffeomorphisms on the invariant graph $\Gamma = \operatorname{graph} u$, with lifts
		\begin{equation}\label{eq:gu-def}
			g_u(x) = x + u(x), \qquad
			g_u^{-1}(x) = x - u(x) + \varphi(x),
		\end{equation}
		Hence $g_u'(x) = 1 + u'(x) > 0$, and the second component satisfies the invariance equation
		\begin{equation}\label{eq:u-invariance-eq}
			u\bigl(g_u(x)\bigr) = u(x) + \varphi\bigl(g_u(x)\bigr).
		\end{equation}
		Write $s(x) := u'(x)$, $m_s := \min_{\mathbb{T}} s$, and $M_s := \max_{\mathbb{T}} s$. By $\int_0^1 s(x)\,\mathrm{d}x = 0$ and $1 + s(x) > 0$, we have $-1 < m_s \le 0 \le M_s$. Differentiating~\eqref{eq:u-invariance-eq} and dividing by $1 + s(x)$ gives
		\begin{equation}\label{eq:s-iteration-eq}
			s\bigl(g_u(x)\bigr) = \frac{s(x)}{1 + s(x)} + \varphi'\bigl(g_u(x)\bigr).
		\end{equation}
		Since $T(t) := t/(1 + t)$ is strictly increasing on $(-1, \infty)$, applying~\eqref{eq:s-iteration-eq} at $s(z) = M_s$, with $x = g_u^{-1}(z)$, and at $s(x_-) = m_s$ gives
		\[
		M_s \le \frac{M_s}{1 + M_s} + d, \qquad
		m_s \le s\bigl(g_u(x_-)\bigr) \le \frac{m_s}{1 + m_s} + d.
		\]
		Clearing denominators and rearranging gives $M_s^2 - d M_s - d \le 0$ and $m_s^2 - d m_s - d \le 0$. Solving these quadratic inequalities yields~\eqref{eq:sqrt-slope-bound}.
	\end{proof}
	
	For $f_\lambda$, equation~\eqref{eq:phi-C4-bound} gives $d = \|\varphi_{\lambda^2}'\|_\infty \le C_0 \lambda^2$. When $0 < \lambda \le 1$, Lemma~\ref{lem:sqrt-slope} and $\sqrt{d^2 + 4d} \le d + 2\sqrt{d}$ imply
	\[
	\|u_{\lambda, k}'\|_\infty \le d + \sqrt{d} \le \bigl(C_0 + \sqrt{C_0}\bigr)\lambda,
	\]
	This proves~\eqref{eq:thm1.1-slope-bound} and completes the proof of Theorem~\ref{thm:herman-hyperbolic-family}.
	
\subsection{The two exterior regions and the interior of Herman's island}
\label{subsec:herman-island-application}

In the original coordinates, write the upper and lower separatrix functions and the enclosed island as
\[
\psi_+(x)=g_\lambda(x)-x,\qquad
\psi_-(x)=g_\lambda^{-1}(x)-x,\qquad
D_\lambda=\{(x,y):0<x<1,\ \psi_-(x)<y<\psi_+(x)\},
\]
where $\psi_\pm$ are extended $1$-periodically. On $(0,1)$ we have $\psi_-<0<\psi_+$, and both vanish at the endpoints. In the universal cover, their graphs over $[0,1]$ form a heteroclinic cycle between $o_0=(0,0)$ and $o_1=(1,0)$: $\widetilde f_\lambda$ induces $g_\lambda$ on the upper arc (moving points toward $o_1$) and $g_\lambda^{-1}$ on the lower arc (moving points toward $o_0$). Thus $\widetilde f_\lambda$ preserves $D_\lambda$ and each of its integer translates, and every orbit in $\overline D_\lambda$ has horizontal rotation number zero.

To relate the lower branch to the upper branch, we use the reversing involution of $f_\lambda$.

\begin{Lemma}\label{rem:negative-rotation}
	Define the involution $\widetilde{\mathcal S}_\lambda(x,y):=\bigl(x,\,\varphi_{\lambda^2}(x)-y\bigr)$ on $\mathbb R^2$. Then $\det D\widetilde{\mathcal S}_\lambda=-1$ and
	\begin{equation}\label{eq:reversing-symmetry}
		\widetilde{\mathcal S}_\lambda^2=\operatorname{id},\qquad
		\widetilde{\mathcal S}_\lambda\circ\tau=\tau\circ\widetilde{\mathcal S}_\lambda,\qquad
		\widetilde{\mathcal S}_\lambda\circ\widetilde f_\lambda\circ\widetilde{\mathcal S}_\lambda=\widetilde f_\lambda^{-1}.
	\end{equation}
	Consequently:
	\begin{enumerate}[label=(\roman*), font=\normalfont]
		\item If $\gamma:\mathbb R\to\mathbb R^2$ satisfies $\gamma(t+1)=\tau\bigl(\gamma(t)\bigr)$ and $\widetilde f_\lambda\bigl(\gamma(t)\bigr)=\gamma(t+\omega)$, then $\widetilde\gamma:=\widetilde{\mathcal S}_\lambda\circ\gamma$ satisfies $\widetilde\gamma(t+1)=\tau\bigl(\widetilde\gamma(t)\bigr)$ and $\widetilde f_\lambda\bigl(\widetilde\gamma(t)\bigr)=\widetilde\gamma(t-\omega)$. In particular, $\widetilde{\mathcal S}_\lambda$ maps $\operatorname{graph}u$ to $\operatorname{graph}(\varphi_{\lambda^2}-u)$ and sends the upper separatrix $y=\psi_+(x)$ to the lower separatrix $y=\psi_-(x)$.
		\item In the scaled coordinates $(x,r)=(x,y/\lambda)$ and the local symplectic charts~\eqref{eq:linear-chart-7-7} at the integer saddles, the scaled involution $(x,r)\mapsto\bigl(x,\,\lambda\Phi_\lambda(x)-r\bigr)$ takes the form $(Q,P)\mapsto(-P,-Q)$. In particular, it preserves $QP$, exchanges the outgoing and incoming logarithmic phases $\nu_\lambda^{-1}\log|Q|$ and $-\nu_\lambda^{-1}\log|P|$ with sign reversed, and relates the normalized boundary transition lifts $b_+,b_-$ of $G_\lambda=F_\lambda^{M_\lambda}$ along the upper and lower arcs by
		\begin{equation}\label{eq:b-minus-inverse}
			b_-(s)=-b_+^{-1}(-s).
		\end{equation}
	\end{enumerate}
\end{Lemma}

\begin{proof}
	Writing $\widetilde f_\lambda=K\circ D$ with $D(x,y)=(x+y,y)$ and $K(x,y)=\bigl(x,y+\varphi_{\lambda^2}(x)\bigr)$, direct evaluation gives
	\[
	\bigl(\widetilde{\mathcal S}_\lambda\circ\widetilde f_\lambda\circ\widetilde{\mathcal S}_\lambda\bigr)(x,y)
	=\bigl(x+\varphi_{\lambda^2}(x)-y,\,y-\varphi_{\lambda^2}(x)\bigr)
	=(D^{-1}\circ K^{-1})(x,y)
	=\widetilde f_\lambda^{-1}(x,y),
	\]
	which proves~\eqref{eq:reversing-symmetry}. Assertion~(i) follows from~\eqref{eq:reversing-symmetry} and $\varphi_{\lambda^2}-\psi_+=g_\lambda^{-1}-\operatorname{id}=\psi_-$.
	
	For~(ii), near each saddle $x=n$, the endpoint affine formulas~\eqref{eq:endpoint-affine-structure} give $\lambda\Phi_\lambda(x)=(1-c_\lambda)(x-n)$. Substituting this into~\eqref{eq:linear-chart-7-7} shows that $(x,r)\mapsto\bigl(x,\lambda\Phi_\lambda(x)-r\bigr)$ is represented in $(Q,P)$ coordinates by $(Q,P)\mapsto(-P,-Q)$. If an outgoing phase $t$ on the upper arc at $o_0$ reaches the incoming phase $b_+(t)+L$ at $o_1$ after $L$ iterates of $G_\lambda$, applying the involution gives a forward $G_\lambda^L$-orbit on the lower arc from outgoing phase $-b_+(t)-L$ at $o_1$ to incoming phase $-t$ at $o_0$. Subtracting $L$ gives $b_-\bigl(-b_+(t)-L\bigr)=-t-L$ on a basic interval; since $b_\pm$ commute with unit translation, setting $s=-b_+(t)-L=-b_+(t+L)$ yields~\eqref{eq:b-minus-inverse} on $\mathbb R$.
\end{proof}

\paragraph{The exterior regions.}
By Lemma~\ref{rem:negative-rotation}(i), applying $\widetilde{\mathcal S}_\lambda$ to the upper invariant graphs $u_{\lambda,k}^+:=u_{\lambda,k}>\psi_+$ of Theorem~\ref{thm:herman-hyperbolic-family} yields smooth invariant graphs $u_{\lambda,k}^-:=\varphi_{\lambda^2}-u_{\lambda,k}^+<\psi_-$ below the lower branch, with $\psi_--u_{\lambda,k}^-=u_{\lambda,k}^+-\psi_+>0$ and constant-type rotation numbers $-\omega_{\lambda,k}$ (since $|-\omega-p/q|=|\omega-(-p)/q|$). By~\eqref{eq:phi-C4-bound} and~\eqref{eq:thm1.1-slope-bound}, their slopes also satisfy $\|(u_{\lambda,k}^-)'\|_\infty\le\|\varphi_{\lambda^2}'\|_\infty+\|(u_{\lambda,k}^+)'\|_\infty\le C\lambda$. Since $\det D\widetilde{\mathcal S}_\lambda=-1$, the positive-area families in upper exterior collars (Remark~\ref{rem:scope}) are mapped to positive-area families in the corresponding lower exterior collars.

These graphs also bound the horizontal rotation numbers of arbitrary exterior orbits. For any $z=(x_0,y_0)$ strictly above $\operatorname{graph}\psi_+$, uniform convergence allows us to choose $k$ so that $y_0>u_{\lambda,k}^+(x_0)$. Since $\widetilde f_\lambda$ preserves the region above $\operatorname{graph}u_{\lambda,k}^+$, the orbit $(x_j,y_j)=\widetilde f_\lambda^{\,j}(z)$ satisfies
\[
x_{j+1}-x_j=y_j>u_{\lambda,k}^+(x_j)\ge\min_{\mathbb T}u_{\lambda,k}^+>0\qquad(j\ge0),
\]
so $\rho_x^-(z)>0$. Symmetrically, the negative maxima of $u_{\lambda,k}^-$ imply $\rho_x^+(z)<0$ for every point strictly below $\operatorname{graph}\psi_-$.

\paragraph{The interior of the island.}
In the scaled coordinates $(x,r)$, the local chart~\eqref{eq:linear-chart-7-7} shows that $D_\lambda$ occupies $\{Q>0,P<0\}$ near $o_0$ and $\{Q<0,P>0\}$ near $o_1$, where the exact linearization~\eqref{eq:G-lambda-linear} holds. The finite-step transition along the upper arc in Section~\ref{subsec:auto-6-3} extends smoothly to the interior side of the boundary, with boundary transition map $b_+=b_\lambda$. By~\eqref{eq:zeta-lambda-C4}, $b_+'=1+O(\lambda)$ and $\|b_+^{(j)}\|_\infty=O(\lambda)$ for $2\le j\le4$. {These bounds also hold for the inverse. Indeed, put $a=b_+^{-1}$ and write $b=b_+$. Since $b'\ge1/2$ for small $\lambda$, differentiating $b\circ a=\operatorname{id}$ gives $a'=1/(b'\circ a)=1+O(\lambda)$ and
\[
\begin{aligned}
a''&=-\frac{b''\circ a}{(b'\circ a)^3},\qquad
a'''=\frac{3(b''\circ a)^2}{(b'\circ a)^5}-\frac{b'''\circ a}{(b'\circ a)^4},\\
a^{(4)}&=-\frac{15(b''\circ a)^3}{(b'\circ a)^7}
+\frac{10(b''\circ a)(b'''\circ a)}{(b'\circ a)^6}
-\frac{b^{(4)}\circ a}{(b'\circ a)^5}.
\end{aligned}
\]
Thus $\|a^{(j)}\|_\infty=O(\lambda)$ for $2\le j\le4$. Combining this with~\eqref{eq:b-minus-inverse} and~\eqref{eq:transition-distance-translations} gives}
\[
\max_{j\in\{+,-\}}\inf_{c\in\mathbb R}
\|b_j-(\operatorname{id}+c)\|_{C^4(\mathbb T)}
=\max_{j\in\{+,-\}}\|b_j^{(4)}\|_\infty\le C\lambda .
\]
Since $\nu_\lambda\in[1/2,1]$ by~\eqref{eq:nu-M-bounds}, the boundary transition hypotheses of Proposition~\ref{prop:island} hold for all sufficiently small $\lambda$. Moreover, since $\varphi_{\lambda^2}(x)$ is positive for small $x>0$ and negative for $x<1$ near $1$, it has a zero $x_*\in(0,1)$, which gives an explicit fixed point $(x_*,0)\in D_\lambda$ of $\widetilde f_\lambda$. For any fixed point $p\in D_\lambda$ of $\widetilde f_\lambda$, Proposition~\ref{prop:island} (after rescaling back by $y=\lambda r$) then yields the contractible smooth invariant circles surrounding $p$ and converging to $\partial D_\lambda$ exponentially in Hausdorff distance, their clockwise rotation numbers $1/[M_\lambda(2k-\alpha_{\lambda,\mathrm{in}})]$ (oriented by $o_0\to o_1\to o_0$), and the positive-area families in every interior collar of $\partial D_\lambda$.

\paragraph{The fixed segment for Herman's original cutoff.}
{The original cutoff in \cite[Chapter~II, proof of Proposition~4.2, pp.~65--66]{Her83} gives a fixed interval on which $g_\lambda'=1$. To check this explicitly, denote Herman's fixed cutoff width by $0<\epsilon_{\mathrm c}<1/2$, put
\[
\epsilon_{\mathrm c,\lambda}=\frac12-\frac{1/2-\epsilon_{\mathrm c}}{1+\lambda},
\]
and let $\eta_{\mathrm c}(t)=1$ for $t\le0$ and $\eta_{\mathrm c}(t)=0$ for $t\ge\epsilon_{\mathrm c}/10$. Near $1/2$ his formulas read
\[
g_\lambda'(x)-1=
\begin{cases}
\lambda\eta_{\mathrm c}\bigl((1+\lambda)(x-1/2+\epsilon_{\mathrm c,\lambda})\bigr),&x\le1/2,\\
-\dfrac{\lambda}{1+\lambda}\eta_{\mathrm c}(1/2+\epsilon_{\mathrm c}-x),&x\ge1/2.
\end{cases}
\]
For $|x-1/2|\le\epsilon_{\mathrm c}/2$, each relevant argument is at least $\epsilon_{\mathrm c}/2$. Hence $g_\lambda'=1$ on the fixed interval $I=(1/2-\epsilon_{\mathrm c}/2,1/2+\epsilon_{\mathrm c}/2)$, and $g_\lambda(x)=x+d_\lambda^*$ there with $0<d_\lambda^*=O(\lambda)$.} On a smaller fixed interval $J$ about $1/2$, both $x$ and $x-d_\lambda^*$ lie in $I$, whence
\[
\varphi_{\lambda^2}(x)=(x+d_\lambda^*)+(x-d_\lambda^*)-2x=0\qquad(x\in J).
\]
Thus $J\times\{0\}\subset D_\lambda$ is a segment of fixed points on a neighborhood of which $\widetilde f_\lambda(x,y)=(x+y,y)$ is the unipotent shear. Every interior point $(x_0,0)$ of this segment, including $p=(1/2,0)$, is Lyapunov unstable: if $[x_0-2b,x_0+2b]\subset J$ and $0<y_0<b$, the orbit of $(x_0,y_0)$ follows the shear for $n=\lceil b/y_0\rceil$ steps and moves horizontally by at least $b$.

\subsection{Proofs of Theorems~\ref{cor:herman-window} and~\ref{thm:intro-herman}}\label{subsec:herman-conclusions}

\begin{proof}[Proof of Theorem~\ref{cor:herman-window}]
		Choose $N_*(a, \gamma_\xi)$ large enough that $\lambda_N = N^{-a/2} < \lambda_*$ for all $N \ge N_*$. Define
		\begin{equation*}
			k_N = \max\left\{k_{\lambda_N},\; \left\lfloor \alpha_{\lambda_N} + \frac{1}{M_{\lambda_N}r_N}\right\rfloor + 1\right\},
		\end{equation*}
		Then $k_N - \alpha_{\lambda_N} > 1 / (M_{\lambda_N} r_N)$, and Theorem~\ref{thm:herman-hyperbolic-family} gives
		\begin{equation*}
			0 < \omega_N = \frac{1}{M_{\lambda_N}(k_N - \alpha_{\lambda_N})} < r_N.
		\end{equation*}
		In particular, taking $r_N=\lambda_N N^{-\varepsilon}=N^{-a/2-\varepsilon}$ yields a smooth invariant graph satisfying~\eqref{eq:herman-question}. We first fix $N$ and then choose $k_N$ sufficiently large; no uniform bound on $k_{\lambda_N}$ with respect to $N$ is needed.
	\end{proof}

\begin{proof}[Proof of Theorem~\ref{thm:intro-herman}]
For each fixed $0<\lambda<\lambda_*$, Theorem~\ref{thm:herman-hyperbolic-family} gives the sequence of smooth invariant graphs above $\Sigma_\lambda$, smooth conjugacy to constant-type rotations, exponential uniform convergence, and the slope bound with a constant independent of $\lambda$ and the sequence index. Since the positive rotation numbers tend to zero, each prescribed positive window contains infinitely many of them. For the power-law parameters $\lambda_N=N^{-a/2}$, Theorem~\ref{cor:herman-window} gives the conclusion for every sufficiently large $N$ and every prescribed radius $r_N>0$, with the threshold on $N$ independent of the radii. Taking the radii in~\eqref{eq:herman-question} gives the negative answer to Herman's question. Section~\ref{subsec:herman-island-application} supplies the two positive-area exterior families, the bounds on horizontal rotation, the contractible interior circles and the description of the fixed segment. For the reflected graphs the slope bound follows from $(u_{\lambda,k}^-)'=\varphi_{\lambda^2}'-u_{\lambda,k}'$ and \eqref{eq:phi-C4-bound}.
\end{proof}

\begin{Remark}\label{rem:second-derivative-blowup}
		For a qualitative discussion of the divergence of higher derivatives of invariant curves approaching a separatrix with a corner, see also \cite[Remark 1.1]{KK21}. For the fixed $\lambda \in (0, \lambda_*)$ considered here, the right and left derivatives of $\psi_\lambda$ at the origin are $\lambda$ and $-c_\lambda\lambda$, respectively, with jump $J_\lambda := (1 + c_\lambda)\lambda > 0$. Thus $\psi_\lambda(t) + \psi_\lambda(-t) - 2\psi_\lambda(0) = J_\lambda t$ on $0 < t \le \frac{r_0}{2}$. Set $\epsilon_{\lambda, k} := \|u_{\lambda, k} - \psi_\lambda\|_\infty$. The second-order Taylor expansion gives
		\[
		\begin{aligned}
			J_\lambda t - 4\epsilon_{\lambda, k}
			&\le u_{\lambda, k}(t) + u_{\lambda, k}(-t) - 2u_{\lambda, k}(0)\\
			&= \int_0^t (t - s)\bigl(u_{\lambda, k}''(s) + u_{\lambda, k}''(-s)\bigr)\,\mathrm{d}s\\
			&\le t^2 \|u_{\lambda, k}''\|_\infty.
		\end{aligned}
		\]
		For sufficiently large $k$, take $t = 8\epsilon_{\lambda, k}/J_\lambda \le \frac{r_0}{2}$ and use~\eqref{eq:thm1.1-decay-rate} to obtain
		\[
		\|u_{\lambda, k}''\|_\infty \ge \frac{J_\lambda^2}{16\epsilon_{\lambda, k}} \ge \frac{J_\lambda^2}{16C_\lambda} e^{k\nu_\lambda/2} \longrightarrow \infty \qquad (k \to \infty).
		\]
	\end{Remark}

	\appendix
	\section{Auxiliary results}\label{sec:appendix}
	The parameters are chosen in the following order: first fix $(\xi,\gamma_\xi)$, then $\lambda$, and finally let $k\to\infty$. Constants with subscript $f$ or $\lambda$ may depend on this fixed map. Uniformity in $\lambda$ is stated separately whenever needed.
	
	\subsection{Norms and the chain rule}\label{subsec:chain-rule}
	On a fixed compact coordinate domain $D$, set
	\[
	\|u\|_{C^r(D)}=\max_{|\mathbf a|\le r}\sup_D|\partial^{\mathbf a}u|.
	\]
	For vectors we use the maximum component norm, and for matrices its induced operator norm $\|A\|_\infty=\max_i\sum_j|A_{ij}|$, so that $\|AB\|_\infty\le\|A\|_\infty\|B\|_\infty$. The $C^r$ norm of a vector-valued function is the maximum of its component norms. In fixed dimension this is equivalent to using operator norms of higher derivatives, with the equivalence constants absorbed into the estimates. The zeroth-order norm of a lift is taken on its periodic displacement relative to a reference lift. Higher norms of $g_\lambda$ are used only on $[0,1]$: its degree-one extension is generally only a piecewise smooth circle homeomorphism, so derivatives that do not match at the endpoints cannot be treated as smooth periodic derivatives.
	
	By the Leibniz rule and the multivariable chain rule, for $r\ge0$,
	\[
	\|uv\|_{C^r}\le2^r\|u\|_{C^r}\|v\|_{C^r},\qquad
	\|u\circ h\|_{C^r}\le C_r\|u\|_{C^r}
	\left(1+\max_{1\le j\le r}\|D^jh\|_\infty\right)^r.
	\]
	Here the norm of $u$ is taken on a domain containing the range of $h$, and for $r=0$ the last factor is interpreted as $1$.
	
	The one-dimensional higher-order chain rule can be written as follows. Define
	\[
	\begin{gathered}
		B_{r,\ell}(z_1,\dots,z_{r-\ell+1})\\
		:=\sum_{\substack{m_1+\cdots+m_{r-\ell+1}=\ell\\
				m_1+2m_2+\cdots+(r-\ell+1)m_{r-\ell+1}=r}}
		\frac{r!}{m_1!\cdots m_{r-\ell+1}!}
		\prod_{i=1}^{r-\ell+1}\left(\frac{z_i}{i!}\right)^{m_i},
	\end{gathered}
	\]
	where $1\le\ell\le r$, and the sum ranges over nonnegative integers $m_i$ satisfying the constraints above. Differentiating the one-dimensional composition $A \circ B$ a total of $r \ge 1$ times gives, by induction,
	\begin{equation}\label{eq:faa-di-bruno}
		(A \circ B)^{(r)} = \sum_{\ell=1}^r \bigl(A^{(\ell)} \circ B\bigr)\,B_{r, \ell}\bigl(B', \dots, B^{(r-\ell+1)}\bigr).
	\end{equation}
	In particular, $B_{r, 1}(z_1, \dots, z_r) = z_r$. For $\ell \ge 2$, only derivatives of $B$ of order at most $r - 1$ occur in $B_{r, \ell}$.
	
	\subsection{Smooth cutoffs and interval models}\label{subsec:cutoffs}
	Set
	\[
	E(t)=\begin{cases}e^{-1/t},&t>0,\\0,&t\le0,\end{cases}
	\qquad \chi_0(t)=\frac{E(t)}{E(t)+E(1-t)}.
	\]
	The function $\chi_0$ is smooth, equals $0$ on $(-\infty,0]$ and $1$ on $[1,\infty)$, and takes values in $[0,1]$. For $0<r<1/2$,
	\[
	\chi_r(x)=\chi_0\left(\frac{x-r}{1-2r}\right),\qquad
	W(x)=(1-\chi_r(x))x+\chi_r(x)(1-x)
	\]
	defines a smooth interval function that is strictly positive in $(0,1)$ and equals $x$ and $1-x$ near the two endpoints, respectively. Translations, rescalings, and products of $\chi_0$ give the cutoff functions used in Section~\ref{subsec:auto-4-2} and Lemma~\ref{lem:corner-collar}.
	
	The same $W$ also gives a concrete example satisfying the assumptions of Section~\ref{subsec:auto-1-2}, showing that this class is nonempty. Fix $0<r<1/4$, and set
	\[
	W_\lambda^{\mathrm{mod}}(x):=(1-\chi_r(x))x
	+\frac{\chi_r(x)}{1+\lambda}(1-x),\qquad
	g_\lambda^{\mathrm{mod}}(x):=x+\lambda W_\lambda^{\mathrm{mod}}(x).
	\]
	We have $g_\lambda^{\mathrm{mod}}(x)=(1+\lambda)x$ on $[0,r]$ and $g_\lambda^{\mathrm{mod}}(x)=1+(x-1)/(1+\lambda)$ on $[1-r,1]$. Since the $C^5$ norm of $W_\lambda^{\mathrm{mod}}$ is uniformly bounded, sufficiently small $\lambda>0$ ensures that $(g_\lambda^{\mathrm{mod}})'>0$. Moreover, $g_\lambda^{\mathrm{mod}}(0)=0$ and $g_\lambda^{\mathrm{mod}}(1)=1$, so $g_\lambda^{\mathrm{mod}}$ is a smooth diffeomorphism of $[0,1]$. For every $K\Subset(0,1)$, $W_\lambda^{\mathrm{mod}}$ has a uniform positive lower bound on $K$. Thus this family satisfies~\eqref{eq:endpoint-affine-structure}--\eqref{eq:positive-drift}. For this particular family, $W_\lambda^{\mathrm{mod}}=W+O(\lambda)$ in every fixed $C^m([0,1])$ norm, so $W$ is its leading term as $\lambda\to0^+$.
	
	\section*{AI declaration}
We used OpenAI’s ChatGPT to assist in preparing schematic figures, improving the clarity of certain arguments, and proofreading the manuscript.

\end{document}